\documentclass[11pt]{article}
\usepackage[pagebackref,colorlinks=true,urlcolor=blue,linkcolor=red,citecolor=red]{hyperref}

\usepackage{epsfig,graphics, graphicx}
\usepackage{subcaption, comment, bm}
\usepackage[percent]{overpic}
\usepackage{float}
\usepackage{amssymb,bm}
\usepackage{amsthm}
\usepackage{color}
\usepackage[]{amsmath}
\usepackage[]{amsfonts}
\usepackage[]{fancyhdr}
\usepackage[]{graphicx,wrapfig}
\usepackage{enumitem}
\usepackage{array}
\usepackage{multirow}

\renewcommand{\div}{\nabla\cdot}

\newcommand{\cout}[1]{}

\usepackage{graphicx}
\usepackage[dvipsnames]{xcolor}

\newcommand{\D}{\mathrm{d}}

\newcommand{\Bc}{\mathcal{B}}

\newcommand{\Ic}{\mathcal{I}}
\newcommand{\Jc}{\mathcal{J}}

\newcommand{\Lc}{\mathcal{L}}

\newcommand{\Rc}{\mathcal{R}}
\newcommand{\Sc}{\mathcal{S}}
\newcommand{\Tc}{\mathcal{T}}

\newcommand{\Xc}{\mathcal{X}}

\newcommand{\Rb}{\mathbb{R}}

\newcommand{\vf}{\textbf{\textit{f}}}
\newcommand{\vx}{\textbf{\textit{x}}}
\newcommand{\vy}{\textbf{\textit{y}}}

\newcommand{\vw}{\textbf{\textit{w}}}
\newcommand{\vu}{\textbf{\textit{u}}}
\newcommand{\vv}{\textbf{\textit{v}}}

\newcommand{\vgamma}{\pmb{\gamma}}
\newcommand{\vpsi}{\pmb{\psi}}
\newcommand{\vxi}{\pmb{\xi}}

\usepackage{amsthm}
\usepackage{amsfonts}
\newcommand{\Beq}{\begin{equation}}
\newcommand{\Eeq}{\end{equation}}
\newcommand{\beq}{\begin{equation*}}
\newcommand{\eeq}{\end{equation*}}
\newcommand{\bal}{\begin{align}}
\newcommand{\eal}{\end{align}}

\graphicspath{{./Figs/}}

\newtheorem*{theor}{Theorem}
\newtheorem{thr}{Theorem}
\newtheorem{defn}{Definition}
\newtheorem{rem}{Remark}

\title{\vspace{-1cm} Numerical implementation of generalized V-line transforms\\ 
on 2D vector fields and their inversions}

\author{Gaik Ambartsoumian\thanks{Department of Mathematics, University of Texas at Arlington, Arlington, TX, USA.  \url{gambarts@uta.edu}}\and Mohammad Javad Latifi Jebelli\thanks{Department of Mathematics, Dartmouth College, Hanover, NH, USA. \url{mohammad.javad.latifi.jebelli@dartmouth.edu}} \and Rohit Kumar Mishra\thanks{Mathematics Discipline, Indian Institute of Technology, Gandhinagar, Gujarat, India. \url{rohit.m@iitgn.ac.in}}}

\begin{document}
\date{}
\maketitle
\vspace{-5mm} 
\begin{abstract}
The paper discusses numerical implementations of various inversion schemes for generalized V-line transforms on vector fields introduced in \cite{Gaik_Mohammad_Rohit}. It demonstrates the possibility of efficient recovery of an unknown vector field from five different types of data sets, with and without noise. We examine the performance of the proposed algorithms in a variety of setups, and illustrate our results with numerical simulations on different phantoms.
\end{abstract}
\vspace{-5mm}
\section{Introduction}\label{Introduction}
A peculiar class of generalized Radon transforms have recently attracted considerable interest in integral geometry and its imaging applications \cite{amb-book}. These transforms map functions to their integrals along paths or surfaces with a ``vertex'' inside their support, e.g. along broken rays (also called V-lines) \cite{Ambartsoumian-chapter, amb-lat_2019,  Ambartsoumian_Moon_broken_ray_article, ambartsoumian2016numerical, Florescu-Markel-Schotland, florescu2018, Gouia_Amb_V-line, Kats_Krylov-13, Sherson, walker2019broken} and stars \cite{Amb_Lat_star, ZSM-star-14} in $\mathbb{R}^2$,  or over various conical surfaces \cite{amb-lat_2019, gouia2014analytical, Gouia_Amb_V-line, Palamodov2017, Terz_IP-15} in $\mathbb{R}^3$ and higher dimensions. Such operators appear in mathematical models of several imaging techniques based on scattered particles, including single scattering optical tomography \cite{FMS-PhysRev-10, FMS-09}, single scattering X-ray tomography \cite{Kats_Krylov-15, walker2021iterative}, fluorescence imaging \cite{florescu2018}, Compton scattering emission tomography \cite{nguyen-truong-02}, and Compton camera imaging \cite{Terz-review-18}. The integral geometric formulations of image reconstruction problems in these setups are typically obtained through the Born approximation of the solution of the radiative transport equation (RTE) (e.g. see \cite{FMS-PhysRev-10, FMS-09}). This equation describes the propagation of radiation through a medium by way of a balance relation between the numbers of emitted, transmitted, absorbed and scattered particles in an infinitesimal volume \cite{duderstadt1979transport}.

The mathematical models leading to the generalized Radon transforms described above, neglect the effects of polarization of electromagnetic radiation. While that approach can be justified by the relative simplicity of the resulting models, it has been proposed by various authors that  studying the effects of polarization within the framework of the \textit{vector RTE} (e.g. see \cite{fernandez1999, fernandez1993}) may provide additional information about the inhomogeneities in the system \cite{FMS-09}. The analysis of Born approximation of the solution of the vector RTE is a difficult task, and we are not aware of any rigorous results on that subject. Hence, the derivation of an accurate integral-geometric model properly reflecting the physics of single-scattering of polarized light photons is also an open problem. However, it is clear that in such a model the scalar functions corresponding to the attenuation and scattering coefficients will be replaced by $4\times4$ extinction and phase matrix functions. Therefore, instead of recovering the attenuation coefficient from its integrals along V-lines and stars, one may need to recover the extinction matrix or some of its components from its integral transforms along such trajectories. This has motivated consideration of generalized V-line transforms (VLTs) on vector fields and on tensor fields of higher order. It must be noted, that the generalization of the classical X-ray and Radon transforms to vector fields and tensor fields of higher order has been subject of intense research for many decades (e.g. see \cite{paternain2023geometric, SchusterReview2008, sharafutdinov2012integral} and the references in the next paragraph). Our choice of the transformations studied in this paper is primarily influenced by that research.

In our paper \cite{Gaik_Mohammad_Rohit}, we introduced several generalizations of the aforementioned V-line and star transforms from scalar fields to vector fields in $\mathbb{R}^2$. The list of these new operators included the longitudinal  and transverse V-line transforms, their corresponding first moments, and the vector star transform. The first four concepts were motivated by the analogous generalizations of the classical Radon transform to vector fields (e.g. see \cite{abhishek2019support, denisjuk2006inversion, derevtsov2, GriesmaierRohit2018, Holman2013, katsevich2013exact, kim2020three, krishnan2019momentum, krishnan2020momentum, krishnan2018microlocal, krishnan2019solenoidal, VRS, mishra2020full, MonardRohit2021, Rohit_Suman_2021, Norton, Novikov_Sharafutdinov, Sharafutdinov_TRT_2008}). The vector star transform is a natural extension of the longitudinal and transverse VLTs to the case of trajectories with more than two branches. 
In \cite{Gaik_Mohammad_Rohit} we studied various properties of these transforms and derived several exact inversion formulas for them. 

The goal of the current article is the study of the image reconstruction algorithms ensuing from the theoretical results obtained in \cite{Gaik_Mohammad_Rohit}, discussion of their numerical implementations and analysis of their performance in various setups.
Development of reconstruction algorithms based on exact inversion formulas of generalized Radon transforms and their numerical validation are essential tasks in tomography (e.g. see \cite{ambartsoumian2007thermoacoustic, ambartsoumian2016numerical, kunyansky2001, monard2014numerical}). While such undertakings in vector and tensor tomography utilizing integrals along straight lines have been studied before (e.g. see \cite{Andersson2005, Desbat1995, Bukhgeim2004, SchusterReview2008}), this paper is the first work exploring such algorithms for transforms integrating along trajectories with a vertex. In addition to the standard visualization technique for vector fields using colored images of separate scalar components, 
we present some results of our vector field reconstructions on a single image using the RGB color model.  We also provide a link to a webpage containing implementations of the vector star transform and its inversion, where an interested reader can experiment with the image reconstruction of their own phantoms.

The rest of this article is organized as follows. In Section \ref{sec:theory} we give the formal definitions of five integral transforms acting on vector fields and state the theorems containing explicit formulas for reconstruction of vector fields from those transforms. In Section \ref{sec:Implementations} we provide the numerical schemes of inverting the generalized VLTs, as well as examples of their implementations in Matlab on various phantoms. In Section \ref{sec:Star-Python} we present the numerical implementation of the vector star transform and its inverse in Python. These codes are made available by the authors as an open access notebook in the Google Colab, with options for user customized experiments. We finish the paper with Conclusions in Section \ref{sec:conclusions} and Acknowledgements in Section \ref{sec:acknowledge}.


\section{Theoretical Background}\label{sec:theory}
\subsection{Definitions and notations}\label{subsec: Def and notations}
Let us start with an introduction of some notations and the definitions of the operators discussed in this paper. 
Throughout the article, we use a bold font to denote vectors in $\mathbb{R}^2$ (e.g. $\vx$, $\textbf{\textit{u}}$, $\textbf{\textit{v}}$, $\textbf{\textit{f}}$, etc), and a regular font to denote scalar variables (e.g. $t$, $h$, $f_i$, etc). The usual dot product between vectors $\vx$ and $\vy$ is written as $\vx\cdot\vy$. For a scalar function $V(x_1, x_2)$  and a vector field $\vf =(f_1,f_2)$, we use the notations 
\begin{align}\label{eq: definition of div and curl}
\nabla V := \left(\frac{\partial V}{\partial x_1}, \frac{\partial V}{\partial x_2}\right), \enspace  
\nabla^\perp V := \left(-\frac{\partial V}{\partial x_2}, \frac{\partial V}{\partial x_1}\right), \enspace  
D_{\vu} V := \vu \cdot \nabla V, \enspace \\
\operatorname{div}\vf :=  \frac{\partial f_1}{\partial x_1}+ \frac{\partial f_2}{\partial x_2},
\enspace  \operatorname{curl} \vf :=  \frac{\partial f_2}{\partial x_1}- \frac{\partial f_1}{\partial x_2}.
\end{align}
       



\begin{figure}[ht]
\begin{center}
\begin{subfigure}{.45\textwidth}
\centering
 \includegraphics[height=3.5cm]{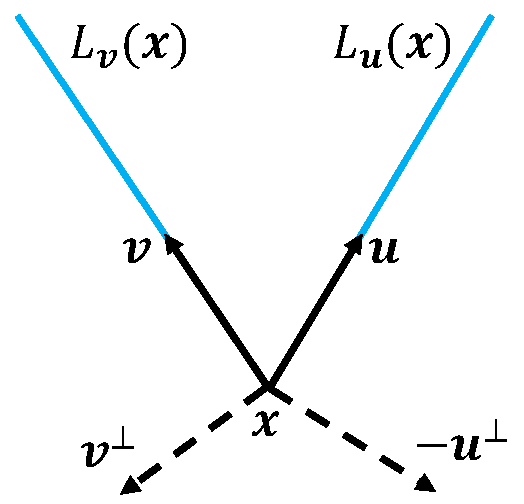} 
 \caption{A V-line with the vertex $\vx$, ray directions $\vu$, $\vv$, and outward normals $-\vu^\perp$, $\vv^\perp$.}
  \label{fig1a}
\end{subfigure} \qquad 
\begin{subfigure}{.45\textwidth}
\centering
\includegraphics[height=3.5cm]{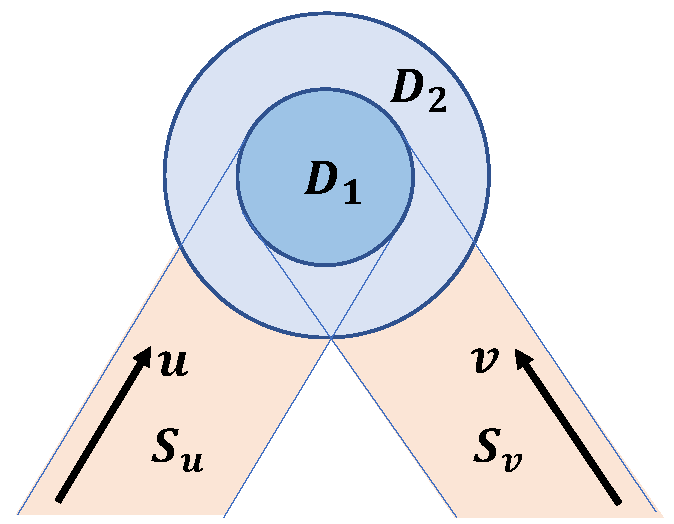}
\caption{A sketch of the compact support of  $\vf$ and the unbounded supports of $\Lc\vf$, $\Tc\vf$, $\Ic\vf$, $\Jc\vf$.} \label{fig1b}
\end{subfigure}
\end{center}
\vspace{-5mm}
\caption{From \cite{Gaik_Mohammad_Rohit}.}
\label{fig1} 
\end{figure}

Let  $\vu,\vv\in S^1$ be a pair of fixed, linearly independent, unit vectors. We denote by  $L_{\vu}(\vx)$ and $L_{\vv}(\vx)$ the rays emanating from $\vx \in \mathbb{R}^2$ in the directions of $\vu$ and $\vv$, i.e.
$$ L_{\vu}(\textbf{\textit{x}}) := \left\{\textbf{\textit{x}} +t \textbf{\textit{u}}: 0 \leq t < \infty\right\} \quad \mbox{ and } \quad  L_{\vv}(\textbf{\textit{x}}) := \left\{\textbf{\textit{x}} +t \textbf{\textit{v}}: 0 \leq t < \infty\right\}.$$
A \textbf{V-line with the vertex $\textbf{\textit{x}}$} is the union of rays $L_{\vu}(\textbf{\textit{x}})$ and $L_{\vv}(\textbf{\textit{x}})$. Note that since  $\textbf{\textit{u}}$ and $\textbf{\textit{v}}$ are fixed, all V-lines have the same ray directions and can be parametrized simply by the coordinates $\textbf{\textit{x}}=(x_1,x_2)$ of their vertices (see Figure \ref{fig1a}). 

\begin{defn}
 The \textbf{divergent beam transform} $\mathcal{X}_{\vu}$ maps a function on $\mathbb{R}^2$ to a set of its integrals along rays, namely
  \begin{equation}\label{def:DivBeam}
   \mathcal{X}_{\vu}h(\vx) :=  \int_{0}^{\infty} h(\vx+t \vu)\,dt.
 \end{equation}
 \end{defn}

The next concept is a natural generalization of the well-known longitudinal ray transform (sometimes also called Doppler transform) \cite{derevtsov2, katsevich2013exact, Thomas_Schuster_2000}, which maps a vector field to the line integrals of its component parallel to the line of integration. If the straight line is substituted by a V-line, then one obtains the following.

\begin{defn} \label{def:definition of V line Doppler transform} Let $\vf = (f_1, f_2)$ be a vector field in $\mathbb{R}^2$ with components $f_i \in C^2_c(\Rb^2)$ for $i =1, 2$. The \textbf{longitudinal V-line transform (LVT)} of $\textbf{\textit{f}}$ is defined as
\begin{align}\label{eq:def V-line transform}
\mathcal{L}_{\vu, \vv}\, \vf\  := -\mathcal{X}_{\vu}  \left(\vf \cdot \vu\right)  + \mathcal{X}_{\vv}  \left(\vf \cdot \vv\right).
\end{align}
\end{defn}

The negative sign in the first term of formula (\ref{eq:def V-line transform}) is due to the direction of traveling along the V-line. It can be interpreted as the path of particles emitted in the direction $-\vu$ at some point outside of the support of the vector field and scattered in the direction $\vv$ at a point $\vx$ inside the support. 

For the next integral transform we need a properly defined notion of the normal unit vector for each branch of the V-line. Given a vector $\vx=(x_1, x_2)$, we denote $\vx^\perp := (-x_2,x_1) $.

\begin{defn}\label{def:definition of transverse V line Doppler transform}
Let $\vf = (f_1, f_2)$ be a vector field in $\mathbb{R}^2$ with components $f_i \in C^2_c(\Rb^2)$ for $i =1, 2$. The \textbf{transverse V-line transform (TVT)} of $\textbf{\textit{f}}$ is defined as
\begin{align}\label{eq:def transverse V-line transform}
\mathcal{T}_{\vu, \vv}\, \vf\  := - \mathcal{X}_{\vu}  \left(\vf \cdot \vu^\perp\right) + \mathcal{X}_{\vv}  \left(\vf \cdot \vv^\perp\right).
\end{align}
\end{defn}
The orientation of normal vectors on each branch of the V-line is chosen towards the same (left) side of the trajectory of the moving particles.
Thus, the transverse V-line transform maps a vector field to the V-line integrals of its ``component'' in the direction of the outward unit normal to the V-line at each point (see Figure \ref{fig1a}). 



\begin{defn}
The \textit{\textbf{first moment  divergent beam transform}} maps a function on $\mathbb{R}^2$ to a set of its weighted integrals along rays, namely
$$\Xc^1_{\vu} h\, (\vx) :=  \int_0^\infty h(\vx + t \vu)\, t\,  dt.  $$ 
\end{defn}

Using the above definition, we generalize the well-known momentum ray transforms mapping vector or tensor fields to weighted integrals of their components along straight lines (e.g. see \cite{abhishek2019support, Andersson2005, kim2020three, krishnan2019momentum, krishnan2020momentum, mishra2020full, Rohit_Suman_2021}) to the case of transforms integrating along V-line paths as follows.

\begin{defn} \label{def:definition of first V line moment  transform} Let $\vf = (f_1, f_2)$ be a vector field in $\mathbb{R}^2$ with components $f_i \in C^2_c(\Rb^2)$ for $i =1, 2$. The \textbf{first moment longitudinal V-line transform (LVT1)} of $\textbf{\textit{f}}$ is defined as
\begin{align}\label{eq:first V line moment  transform}
\mathcal{I}_{\vu, \vv}\, \vf\  &:=  - \Xc_{\vu}^1 (\vf \cdot \vu) + \Xc_{\vv}^1(\vf\cdot \vv).
\end{align}
\end{defn}

\begin{defn} \label{def:definition of first V line transverse moment  transform} Let $\vf = (f_1, f_2)$ be a vector field in $\mathbb{R}^2$ with components $f_i \in C^2_c(\Rb^2)$ for $i =1, 2$. The \textbf{first moment transverse V-line transform (TVT1)} of $\textbf{\textit{f}}$ is defined as
\begin{align}\label{eq:first V line moment  transform2}
\mathcal{J}_{\vu, \vv}\, \vf\  &:=  - \Xc_{\vu}^1 \left(\vf \cdot \vu^\perp \right) + \Xc_{\vv}^1\left(\vf\cdot \vv^\perp \right).
\end{align}
\end{defn}
\vspace{2mm}


\begin{rem}
One can easily verify that $\mathcal{T}_{\vu,\vv}\, \vf = -\mathcal{L}_{\vu,\vv}\, \vf^{\,\perp}$ and $\Jc_{\vu, \vv}\, \vf = - \Ic_{\vu, \vv}\, \vf^{\,\perp}$.
\end{rem}

\vspace{2mm}
\begin{rem}
Since the unit vectors $\vu$ and $\vv$ are fixed, in the rest of the paper we drop the indices $\vu,\vv$ and refer to $\mathcal{T}_{\textbf{\textit{u}}, \textbf{\textit{v}}}$, $\mathcal{L}_{\textbf{\textit{u}},  \textbf{\textit{v}}}$, $\mathcal{I}_{\vu, \vv}$, and $\Jc_{\vu, \vv}$ simply as $\mathcal{T}$, $\mathcal{L}$, $\mathcal{I}$, and $\Jc$.
\end{rem}


Let us assume that  $\operatorname{supp}\vf\subseteq D_1$, where $D_1$ is an open disc of radius $r_1$ centered at the origin. Then $\Lc\vf$, $\Tc\vf$, $\Ic\vf$ and $\Jc \vf$ are supported inside an unbounded domain $D_2\cup S_{\vu} \cup S_{\vv}$, where $D_2$ is a disc of some finite radius $r_2>r_1$ centered at the origin, while $S_{\vu}$ and $S_{\vv}$ are semi-infinite strips (outside of $D_2$) stretching in the direction of $-\vu$ and $-\vv$, respectively (see Figure \ref{fig1b}). It is easy to notice that all three transforms $\Lc\vf$, $\Tc\vf$, $\Ic\vf$ and $\Jc \vf$ are constant along the rays in the directions of $-\vu$ and $-\vv$ inside the corresponding strips $S_{\vu}$ and $S_{\vv}$. In other words, the restrictions of $\Lc\vf$, $\Tc\vf$,  $\Ic\vf$ and $\Jc \vf$ to $\overline{D_2}$ completely define them in $\mathbb{R}^2$. 

\begin{rem}\label{rem:support}
Throughout the paper we assume that  the vector field $\vf$ is supported in $D_1$, and the transforms $\Lc\vf\,(\vx)$, $\Tc\vf\,(\vx)$, $\Ic\vf\,{(\vx)}$,  $\Jc\vf\,{(\vx)}$ are known for all $\vx\in \overline{D_2}$.
\end{rem}

\subsection{Recovery of $\vf$ using $\mathcal{T}\vf$, $\mathcal{L}\vf$, $\mathcal{I}\vf$, and $\Jc \vf$}\label{sec:LVT-TVT}
The following theorems follow directly from the results proven in \cite{Gaik_Mohammad_Rohit}.

\begin{thr}\label{th: Recovery of potential and solenoidal part vector field} Consider a vector field $\vf$ with components in $C_c^2(\mathbb{R}^2)$.
\begin{itemize}
    \item If $\vf$ is a potential vector field, i.e. $\vf =  \nabla V$ for some scalar function $V$ supported in $D_1$, then $V$ can be explicitly reconstructed  from $\Tc \vf$ by solving the following Dirichlet boundary value problem: 
    \begin{align*}
    \left\{\begin{array}{lll}
    &\Delta V(\vx)  =    \displaystyle -\frac{1}{\det(\vv, \vu) }D_{\vu}D_{\vv}\, \Tc \vf\,(\vx)&  \mbox{\textup{ in }} D_1,  \\
      &V(\vx) = 0  &\mbox{ \textup{on} } \partial D_1.
    \end{array}\right.
\end{align*}
    \item If $\vf$ is a solenoidal vector field, i.e. $\vf =  \nabla^\perp W$ for some scalar function $W$ supported in $D_1$, then $W$ can be explicitly reconstructed from $\Lc \vf$ by solving the following Dirichlet boundary value problem: 
\begin{align*}
    \left\{\begin{array}{lll}
    &\Delta W(\vx)  =    \displaystyle \frac{1}{\det(\vv, \vu) }D_{\vu}D_{\vv}\, \Lc \vf\,(\vx)&  \mbox{\textup{ in }} D_1,  \\
      & W (\vx) = 0  &\mbox{ \textup{on} } \partial D_1.
    \end{array}\right.
\end{align*}
\end{itemize}
\end{thr}
\begin{thr}\label{th: Laplacian of components of f}
Consider a vector field $\vf$ with components in $C_c^2(\mathbb{R}^2)$. If $\Lc \vf$ and $\Tc \vf$ are known, then the Laplacian of each component of $\vf$ can be explicitly recovered using the following formulas:
\begin{align}
   \Delta f_1 =  -\frac{1}{\det(\vv, \vu)}D_{\vv}D_{\vu}\left\{\nabla \cdot \begin{pmatrix}
   \Tc \vf\\
   \Lc \vf
   \end{pmatrix}\right\},\label{eq:Laplacian recovery f1}\\
    \Delta f_2 =  \frac{1}{\det(\vv, \vu)}D_{\vv}D_{\vu}\left\{\nabla^\perp \cdot \begin{pmatrix}
   \Tc \vf\\
   \Lc \vf
\end{pmatrix}\right\}.\label{eq:Laplacian recovery f2}
\end{align}
\end{thr}
Therefore, one can reconstruct the entire vector field by solving for $f_1$ and $f_2$ the Dirichlet boundary value problems corresponding to equations \eqref{eq:Laplacian recovery f1} and \eqref{eq:Laplacian recovery f2}. 

A key feature in proving the above results is the possibility of expressing
$\operatorname{curl}  \vf$ and $\operatorname{div}  \vf$ in terms of given data $\Lc \vf$ and $\Tc \vf$. More specifically, we have the following identities (see \cite[Theorem 3 and Theorem 4] {Gaik_Mohammad_Rohit} for details): 
\begin{align}
\operatorname{curl} \vf &= \frac{1}{\det(\vv, \vu)}D_{\vu}D_{\vv}\, \Lc \vf, \label{eq:curl-explicit}\\
\operatorname{div} \vf &= -\frac{1}{\det(\vv, \vu)}D_{\vu}D_{\vv} \Tc \vf. \label{eq:div-explicit}
\end{align}
These identities can also be combined with the next theorem to address the problem of reconstructing a vector field using the first moments of longitudinal and transverse V-line transform. 
\begin{thr}\label{th: reconstruction with integral moments}
Consider a vector field $\vf$ with components in $C_c^2(\mathbb{R}^2)$, and let $\vw=(\vv-\vu)/||\vv-\vu||$.
\begin{itemize}
    \item If $\Lc \vf$ and $\Ic \vf$ are known, then each component of $\vf$ can be explicitly recovered using the relation (\ref{eq:curl-explicit}) and the following formulas:
\begin{align}
   f_1(\vx)=  \frac{1}{\|\vv - \vu\|}D_{\vv}D_{\vu}\int_0^\infty\left\{ \frac{\partial \Ic \vf}{\partial x_1} + u_2 \Xc_{\vu}^1( \operatorname{curl}  \vf) -v_2 \Xc_{\vv}^1 ( \operatorname{curl} \vf) \right\}(\vx + t \vw) dt ,\label{eq:recovery f1}\\
  f_2(\vx)=  \frac{1}{\|\vv - \vu\|}D_{\vv}D_{\vu}\int_0^\infty\left\{ \frac{\partial \Ic \vf}{\partial x_2} - u_1 \Xc_{\vu}^1( \operatorname{curl} \vf) + v_1 \Xc_{\vv}^1 ( \operatorname{curl} \vf) \right\}(\vx + t \vw) dt.\label{eq: recovery f2}
\end{align}
    \item If $\Tc \vf$ and $\Jc \vf$ are known, then each component of $\vf$ can be explicitly recovered using the relation (\ref{eq:div-explicit}) and the following formulas:
\begin{align}
   f_1(\vx)&=  \frac{1}{\|\vv - \vu\|}D_{\vv}D_{\vu}\int_0^\infty\left\{ -\frac{\partial \Jc \vf}{\partial x_2} -  u_1 \Xc_{\vu}^1(\operatorname{div} \vf) + v_1 \Xc_{\vv}^1 (\operatorname{div} \vf) \right\}(\vx + t \vw) dt ,\label{eq:recovery f3}\\
  f_2(\vx)&=  \frac{1}{\|\vv - \vu\|}D_{\vv}D_{\vu}\int_0^\infty\left\{ \frac{\partial \Jc \vf}{\partial x_1} - u_2 \Xc_{\vu}^1(\operatorname{div} \vf) + v_2 \Xc_{\vv}^1 (\operatorname{div} \vf) \right\}(\vx + t \vw) dt.\label{eq: recovery f4}
\end{align}
\end{itemize}
\end{thr}

\subsection{Vector star transform and its inversion}\label{sec:star}
The VLTs discussed in the previous section comprise a  difference of two divergent beam transforms. The star transform is composed of an arbitrary linear combination of the corresponding divergent beam transforms. In this section, we give a formal definition of the star transform on vector fields and present its inversion formula derived in \cite[Section 6]{Gaik_Mohammad_Rohit}. 

\begin{defn}\label{def:star_d}
Let $\vgamma_1, \dots, \vgamma_m\in S^1$ be a set of fixed, unit vectors in $\mathbb{R}^2$, and $c_1, \dots , c_m$ be a set of non-zero weights in $\mathbb{R}$. The \textbf{vector star transform} $\mathcal{S}\vf$ of a vector 
field
$\vf$ is defined as 
\begin{equation}\label{def:star}
\mathcal{S}\vf  
:=\sum_{i=1}^m c_i \, \mathcal{X}_{\vgamma_i}  \begin{bmatrix}
  \vf \cdot \vgamma_i  \\
  \vf \cdot \vgamma_i^{\perp}
\end{bmatrix}, 
\end{equation}
where $\mathcal{X}_{\vgamma_i}$ is applied to the vector in the right-hand side of (\ref{def:star}) component-by-component. 
\end{defn}

Note that, in contrast to each VLT discussed in the previous section, the vector star transform data contains integrals of both the longitudinal and the transverse components of the vector field, which suggests the possibility of full recovery of the field from that data. 

\begin{defn}
 We call a star transform $\mathcal{S}$ \textbf{symmetric}, if  $m=2k$ for some $k\in\mathbb{N}$ and (after possible re-indexing) $\vgamma_i=-\vgamma_{k+i}$  with $c_i=-c_{k+i}$ for all $i=1,\ldots, k$. 
\end{defn}

Let $\mathcal{R}h(\vpsi,s)$ denote the (classical) Radon transform of a scalar function $h$ in $\mathbb{R}^2$, along the line normal to the unit vector $\vpsi\in S^1$ and at a signed distance $s\in\mathbb{R}$ from the origin. 

\begin{thr}\label{th:star}
Consider the vector star transform $\mathcal{S}\vf$ with branch directions $\vgamma_1, \dots , \vgamma_m$, and let
\begin{equation}
\vgamma(\vpsi) := -\sum_{i=1}^m  \frac{c_i\, \vgamma_i}{ \vpsi \cdot \vgamma_i } \; \in \mathbb{R}^2, 
\quad \mbox{and}\quad
Q(\vpsi) := \begin{bmatrix}
 \vgamma(\vpsi) \\
 \vgamma(\vpsi)^{\perp}
\end{bmatrix}^{-1} .
\end{equation}
If the unit vector $\vpsi$ is in the domain of $Q(\vpsi)$, then 
\begin{equation}\label{eq:star-inv}
Q(\vpsi) \frac{d}{ds}\mathcal{R}(\mathcal{S}\vf) (\vpsi,s) = \mathcal{R}\vf\,{(\vpsi,s)},
\end{equation}
where $\mathcal{R}\vf$ is the component-wise Radon transform of a vector field in $\mathbb{R}^2$. 
\end{thr}

It was shown in \cite{Gaik_Mohammad_Rohit} that the function $Q(\vpsi)$ is defined for all but finitely many $\vpsi\in S^1$, if and only if $\mathcal{S}$ is not symmetric. In that case, one can recover $\vf$ from $\mathcal{S}\vf$ by applying to the left-hand-side of equation (\ref{eq:star-inv}) any inversion formula of the classical Radon transform.


\begin{rem}
When $m=2$ and $c_1=-c_2=1$, the vector star transform becomes $\mathcal{S}\vf=(\Lc\vf,\Tc\vf)$. Hence, Theorem \ref{th:star} provides another approach to the recovery of the full vector field $\vf$ from its longitudinal and transverse VLTs. In the special case when $\vgamma_1 = -\vgamma_2$ (and only in that case), the matrix $Q(\vpsi)$ is undefined for any $\vpsi$, and the corresponding transform is not invertible. 
\end{rem}


\section{Numerical Implementation}\label{sec:Implementations}
In this section, we provide the numerical schemes of inverting the generalized VLTs, as well as examples of their implementations on various phantoms. In particular, we demonstrate an efficient recovery of the unknown vector field $\vf$ from the following five data sets.
\begin{enumerate}
\item \textbf{Special vector fields:}  Either $\Lc \vf$ or $\Tc \vf$ is used to reconstruct, respectively, a solenoidal or a potential vector field $\vf$ (see Theorem \ref{th: Recovery of potential and solenoidal part vector field}). 
\item $\Lc \vf$ and $\Tc \vf$ are used together to recover the full unknown vector field $\vf$ (see Theorem \ref{th: Laplacian of components of f}). 
\item A combination of $\Lc \vf$ and its first moment $\Ic \vf$ is used to recover $\vf$ (see Theorem \ref{th: reconstruction with integral moments}). 
\item A combination of $\Tc \vf$ and its first moment $\Jc \vf$ is used to recover $\vf$ (see Theorem \ref{th: reconstruction with integral moments}). 
\item $\Sc \vf$ is used to recover $\vf$ (see Theorem \ref{th:star}).
\end{enumerate}

To avoid cumbersome notations of discretized variables, in this section we will denote the components of the vector variable $\vx$ by $\vx=(x,y)$, instead of $\vx=(x_1,x_2)$.


\subsection{Description of phantoms}
\begin{figure}[H]
\centering
\includegraphics[width=0.95\textwidth]{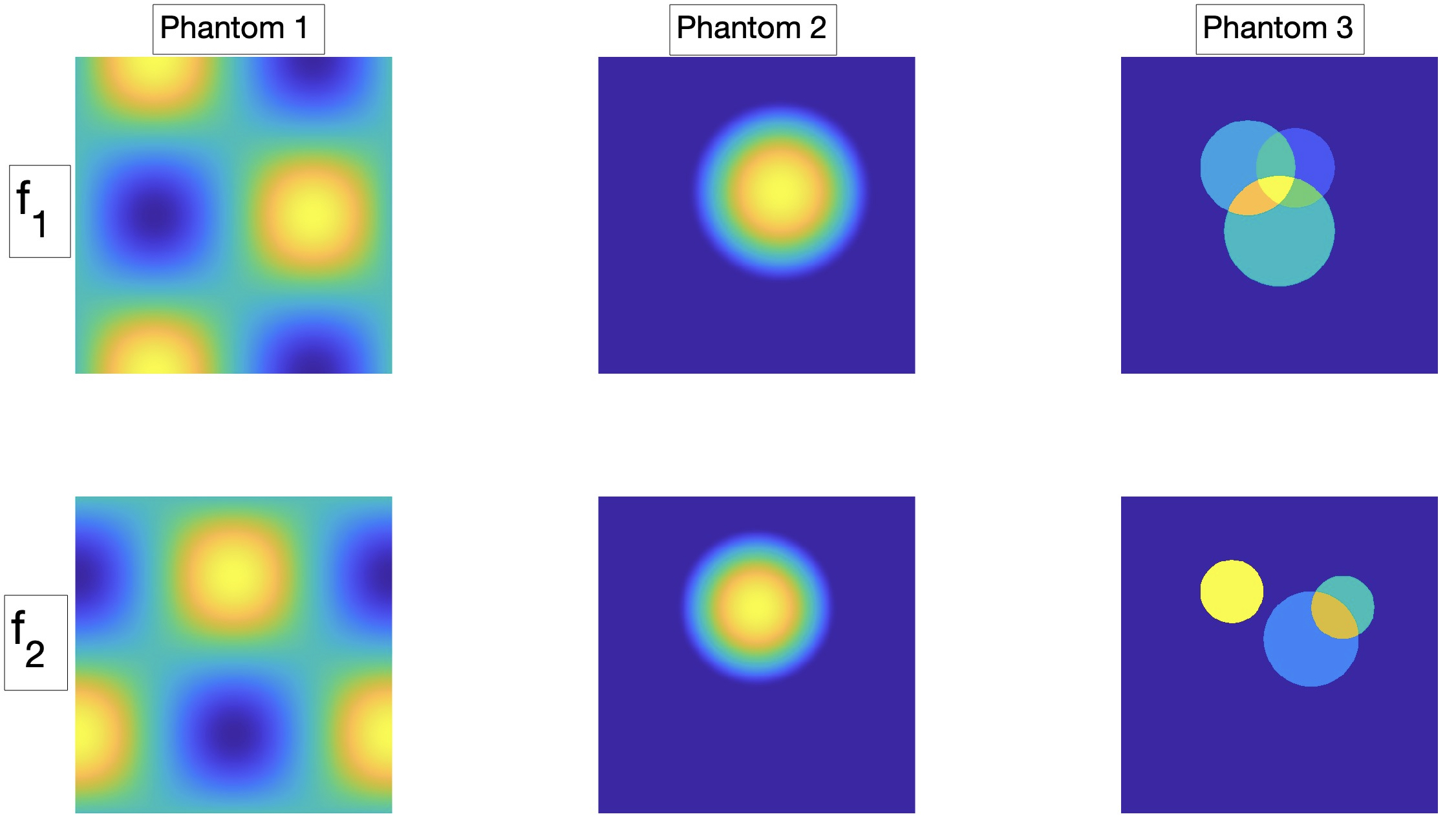}
    \caption{Images of the scalar components of the phantoms used in the numerical simulations. }\label{fig:phantoms}
  \end{figure}

For the five cases described above, we test the performance of the numerical algorithms using various combinations of the following three vector field phantoms defined on $[-1, 1] \times [-1, 1]$ and depicted in Figure \ref{fig:phantoms}.
\begin{itemize}
\item \textbf{Phantom 1:} $\vf(x, y) = (f_1(x, y), f_2(x, y))$, where 
$$f_1(x,y) =   1 + \sin(\pi x)\cos(\pi y), \quad \mbox{ and }\quad  f_2(x,y) = 1 + \sin(\pi y)\cos(\pi x).$$
\item \textbf{Phantom 2:} $\vf(x, y) = (f_1(x, y), f_2(x, y))$, where 
$$f_1(x,y) =  \left\{\begin{array}{cc}
e^{-0.4/\{0.4 - [(x - 0.15)^2+(y -0.15)^2 ]\}},    &  (x - 0.15)^2+(y -0.15)^2 < 0.4 \\
0,   & (x - 0.15)^2+(y -0.15)^2 \geq 0.4
\end{array}\right.$$  and
$$f_2(x,y) = \left\{\begin{array}{cc}
e^{-0.3/\{0.3 - [x^2+(y- 0.3)^2]\}},    &  x^2+(y- 0.3)^2 < 0.3\\
0,   & x^2+(y- 0.3)^2 \geq 0.3.
 \end{array}\right.$$
\item \textbf{Phantom 3:} $\vf(x, y) = (f_1(x, y), f_2(x, y))$, where $f_1$ and $f_2$ are sums of three weighted characteristic functions of disks of different radii $r_j$ and center locations $(x_j,y_j)$. Namely,
\end{itemize}

\begin{center}
\begin{tabular}{ |c|c|c|c|c| } 
 \hline
 $f_1$ & $r$ & $x$ & $y$ & $w$ \\ 
 \hline
 disc 1 & $0.25$ & $0.1$ & $0.3$ & $0.3$ \\ 
 disc 2 & $0.35$ & $0$ & $-0.1$ & $0.9$ \\ 
 disc 3 & $0.3$ & $-0.2$ & $0.3$ & $0.7$ \\
 \hline
\end{tabular}
\hspace{10mm} 
\begin{tabular}{ |c|c|c|c|c| } 
 \hline
 $f_2$ & $r$ & $x$ & $y$ & $w$\\ 
 \hline
 disc 1 & $0.3$ & $0.2$ & $0.1$ & $0.25$\\ 
 disc 2 & $0.2$ & $0.4$ & $0.3$ & $0.45$\\ 
 disc 3 & $0.2$ & $-0.3$ & $0.4$ & $0.9$\\
 \hline
\end{tabular}
\end{center}

\vspace{6mm}

Each of these phantoms has its own specific characteristics examining the pros and cons of the five inversion techniques discussed in the paper. For example, the support of Phantom 1 is not separated away from the boundary of the unit square, which creates difficulties in certain algorithms. Meanwhile, Phantom 3 is piecewise constant, i.e. it lacks the smoothness required in the hypotheses of the inversion results.

\subsection{Data formation}\label{subsec:data formulation}

Unless otherwise specified, in the numerical simulations involving the V-line transforms, the unit vectors defining the V-lines are taken to be $\vu = (\cos \pi/4, \sin \pi/4)$ and $\vv = (\cos 3\pi/4, \sin 3\pi/4)$.{We discuss the effects of the V-line opening angle on reconstructions in Section \ref{Sec:angles-effect}. In the image reconstructions using the vector star transform, we employ stars with three branches defined by angles $\phi_1 = 0$, $\phi_2= 2\pi/3$, and $ \phi_3= 4\pi/3$, and the weights $c_i =  1$ for $i = 1, 2, 3$.

 All integral transforms under consideration are linear combinations of the divergent beam transform and its first moment of various projections of the vector field $\vf$. Therefore, to generate the forward data (corresponding to the V-line transforms and the vector star transform) one needs to have numerical algorithms for computing the divergent beam transform and its first moment of a given scalar function of two variables. We discuss below the process of computing those transforms for a pixelized image $F$. \\


\noindent \textbf{Numerical implementation of the divergent beam transform.} We start with an $m \times m$  pixelized image $F$ defined on $[-1, 1] \times [-1, 1]$. The divergent beam transform of $F$ will also be of the same size $m \times m$, as the rays are parametrized by the coordinates of their vertices, and we consider only the rays emanating from the centers of pixels. To compute the divergent beam transform of $F$ at a vertex $\vx = (x, y)$ in the direction $\vu = (\cos \phi, \sin \phi)$, we first find the intersections of the ray emanating from $\vx$ in the direction $\vu$ and the boundaries of square pixels appearing on the path of this ray. Then, for each such pixel $(i,j)$ we take the product of $F(i,j)$ and the length of the line segment of the ray inside the pixel $(i,j)$. Summing up these products over all such pixels yields the divergent beam transform of $F$ at $\vx$ in the direction $\vu$.\\

\noindent \textbf{Numerical implementation of the first moment of the divergent beam transform.} We use a similar approach to compute this weighted integral. The only difference is that here each term of the sum described above is a product of three quantities. We first multiply $F(i,j)$ by the distance between the center of the pixel $(i, j)$ and the vertex $\vx$ of the ray, and then by the length of the line segment of the ray inside the pixel $(i,j)$. Notice that this method of computing the first moment of the divergent beam transform is not exact, since we use the same constant as the distance between the vertex and any point of the ray inside the pixel. One can easily modify the procedure to account for the variable distance too, but the difference in the generated forward data is negligible for a reasonably fine discretization of the image. \\ 

To generate $\Lc \vf$ and $\Tc \vf$, we evaluate numerically the divergent beam transforms $\mathcal{X}_{\vu}$ of projections $\langle\vf, \vu\rangle$, $\langle\vf, \vu^\perp\rangle$, and $\mathcal{X}_{\vv}$ of $\langle\vf, \vv\rangle$, $\langle\vf, \vv^\perp\rangle$, and combine these quantities according to formulas \eqref{eq:def V-line transform} and \eqref{eq:def transverse V-line transform}. 
The data for $\Ic \vf$ and $\Jc \vf$ are generated in a similar fashion by numerical evaluation of $\mathcal{X}_{\vu}^1$ and $\mathcal{X}_{\vv}^1$ of the appropriate projections and combining the resulting quantities according to formulas (\ref{eq:first V line moment  transform}) and (\ref{eq:first V line moment  transform2}).
 Finally, to obtain the first and the second components of $\mathcal{S}\vf$, we combine the divergent beam transforms $\mathcal{X}_{\vgamma_i}$ of $\langle \vf, \vgamma_i\rangle$ and $\langle \vf, \vgamma_i^\perp\rangle$ respectively, using all $i=1,\ldots, m$ and formula (\ref{def:star}). 

\begin{rem}
    In many of our numerical experiments we add $5\%$, $10 \%$, and $20 \%$ noise to the integral transforms data before applying the inversion procedures.
\end{rem}

\begin{rem}
    The image reconstruction from LVT and TVT data involves solving a Laplace equation, which requires an inversion of an $m^2 \times m^2$ matrix. To curb the computational time, in the numerical implementation of inverting the LVT and TVT we use images with a resolution of $160 \times 160$ pixels. In the problems of recovering a vector field from the other sets of integral transforms, we use images with a resolution of $512 \times 512$ pixels.
\end{rem}

\subsection{Recovery of solenoidal and potential vector fields}\label{subsec:Recovery of solenoidal and potential vector field}
It was shown in Theorem \ref{th: Recovery of potential and solenoidal part vector field} that the solenoidal and the potential vector fields can be recovered just from the knowledge of their longitudinal and transverse V-line transforms, respectively. In this subsection we present such reconstructions using only one of the transformations.

Each reconstruction requires numerically solving a boundary value problem for the Laplace equation, which we achieve through the finite difference method discussed below. 
We present the implementation details for a solenoidal vector field $\vf$ recovered from the knowledge of its longitudinal V-line transform $\Lc\vf$; the recovery of a potential vector field $\vf$ from its transverse V-line transform $\Tc\vf$ follows similarly with obvious changes. 

Let $\vf =  \nabla^\perp W$ be the unknown solenoidal vector field, and $\Lc \vf$ be the given data.  Recall from the second part of Theorem \ref{th: Recovery of potential and solenoidal part vector field} that the function $W$ should satisfy the following relations:
\begin{align*}
    \left\{\begin{array}{lll}
    \Delta W(\vx)  &=    \displaystyle \frac{1}{\det(\vv, \vu) }D_{\vu}D_{\vv} \Lc \vf\,(\vx)&  \mbox{\textup{ in }} D_1,  \\
       W (\vx) &=\  0  &\mbox{ \textup{on} } \partial D_1.
    \end{array}\right.
\end{align*}

First, we compute the gradient $\left(\partial_{x} \Lc \vf,  \partial_{y} \Lc \vf\right)$ with the help of the Matlab function \texttt{gradient}.  Then the directional derivative $D_{\vv} \Lc \vf$ is obtained at every grid point by the direct computation:
$$ D_{\vv} \Lc \vf (x_i,y_j) = -\frac{1}{\sqrt{2}}\partial_{x} \Lc \vf (x_i,y_j) + \frac{1}{\sqrt{2}}\partial_{y} \Lc \vf (x_i,y_j). $$
Applying the same process we get $D_{\vu}D_{\vv} \Lc \vf\,(x_i,y_j)=-\Delta W (x_i,y_j)$, since ${\det (\vv, \vu)} = - 1$ for our choice of $\vu$ and $\vv$. 

Next, we describe our method for numerically solving the boundary value problem for $W$, which will complete the reconstruction of the solenoidal vector field $\vf$ from $\Lc\vf$. In fact, we discuss the numerical solution of a general Dirichlet boundary value problem for the Poisson equation, as it will also appear with different source terms in other places of our paper. In particular, we write the numerical scheme for the following problem:
\begin{align}\label{eq:Poisson equation}
    \left\{\begin{array}{cll}
    - \Delta u &=    f&  \mbox{\textup{ in }} \Omega = [-1, 1] \times [-1, 1],  \\
       u &=\  g  &\mbox{ \textup{on} } \partial \Omega.
    \end{array}\right.
\end{align}
Dividing $\Omega$ into $N \times N$ uniform pixels with the pixel size $h\times h$, we write the central difference approximation for the second-order derivatives at an interior grid point $(x_i, y_j) $ as:
\begin{align}
    \left(\frac{\partial^2 u}{\partial x^2} \right) (x_i, y_j) = \frac{u_{i+1,j}-2u_{i,j} + u_{i-1,j} }{h^2}, \quad  \quad 
    \left(\frac{\partial^2 u}{\partial y^2}\right)(x_i, y_j) = \frac{u_{i,j+1}-2u_{i,j} + u_{i,j-1} }{h^2} 
\end{align}
where $u_{i,j} =  u(x_i, y_j)$. 
Then an approximation of the Laplace operator at an interior grid point $(x_i, y_j)$ can be written as:
\begin{equation*}\label{eq: Finite difference Laplace}
 - (\Delta_h u)_{i,j}  
=   \frac{4u_{i,j}-u_{i-1,j}-u_{i+1,j}-u_{i,j-1}-u_{i,j+1}}{h^2}.
\end{equation*}
Consequently, a finite difference version of the Poisson equation \eqref{eq:Poisson equation} is given by 
\begin{align}\label{eq:finite difference Laplace}
         - (\Delta_h u)_{i,j}   =   f_{i,j} \quad \mbox{\textup{at the interior grid points}}.
\end{align}
We write the interior $(N -2)  \times (N-2) $ grid points in one row using a single index $k = 1$ to $(N-2)^2$ for $u_k = u_{i(k),j(k)}$ and $f_k = f_{i(k),j(k)}$. We use the index map $(i, j) \rightarrow k = (N-2)(i-2)+(j-1)$ for $2 \leq i, j \leq N-1$. With this choice of indexing, equation \eqref{eq:finite difference Laplace} can be written as a matrix equation
\begin{align}
    A U = F,
\end{align}
where $A$ is an $(N-2)^2 \times (N-2)^2$ matrix of the following block tridiagonal structure:
\begin{align*}
   A = -
   \begin{pmatrix}
        B & I & 0 &  \cdots & 0 & 0 & 0\\
        I & B & I &  \cdots & 0 & 0 & 0\\
       \vdots & \vdots &  \vdots &  \ddots & \vdots & \vdots & \vdots \\
       0 & 0 &  0 & \cdots & I & B & I\\
       0 & 0 &  0 & \cdots &  0 & I & B
    \end{pmatrix},
    \mbox{ where }  B =   \begin{pmatrix}
        -4 & 1 & 0 &  \cdots & 0 & 0 & 0\\
        1 & -4 & 1 &  \cdots & 0 & 0 & 0\\
       \vdots & \vdots &  \vdots &  \ddots & \vdots & \vdots & \vdots \\
       0 & 0 &  0 & \cdots & 1 & -4 & 1\\
       0 & 0 &  0 & \cdots &  0 & 1 & -4
    \end{pmatrix}_{(N-2) \times (N-2)}
\end{align*}
and $I$ is the identity matrix of order $(N-2) \times (N-2)$, 
$U = (u_k)_{1 \leq k \leq (N-2)^2}$ and $F = h^2(\tilde{f}_k)_{1 \leq k \leq (N-2)^2}$. Here $\tilde{f}$ represents the modified source term, which satisfies $\tilde{f}_{i,j} =  f_{i,j}$, for $3 \leq i, j \leq N-2$, and involves the boundary terms (i.e. the given data $g(i,j)$) for other indices. More specifically, we have
\begin{align*}
    \begin{array}{lll}
      \tilde{f}_{2,2}  = f_{2,2} + \frac{1}{h^2}\left(g_{1, 2}+ g_{2, 1}\right), &  & \tilde{f}_{N-1,2}  = f_{N-1,2} + \frac{1}{h^2}\left(g_{N, 2}+ g_{N-1, 1}\right), \\ 
      \tilde{f}_{2,j}  =   f_{2,j} + \frac{1}{h^2}g_{1, j}, \, \mbox{for } 3 \leq j \leq N-2, & & \tilde{f}_{N-1,j} =   f_{N-1,j} + \frac{1}{h^2}g_{N, j}, \, \mbox{for } 3 \leq j \leq N-2, \\
      \tilde{f}_{i,2}  =   f_{i,2} + \frac{1}{h^2}g_{i, 1}, \, \mbox{for } 3 \leq i \leq N-2, &  &\tilde{f}_{i,N-1} =   f_{i,N-1} + \frac{1}{h^2}g_{i, N}, \, \mbox{for } 3 \leq i  \leq N-2,\\
      \tilde{f}_{2,N-1}  = f_{2,N-1} + \frac{1}{h^2}\left(g_{1, N-1}+ g_{2, N}\right), & &  \tilde{f}_{N-1,N-1}  = f_{N-1,N-1} + \frac{1}{h^2}\left(g_{N, N-1}+ g_{N-1, N}\right).
    \end{array}
\end{align*}
Finally, we solve the system of linear equations $A U = F$ to get $U$ as a numerical approximation of the solution $u$ of the required boundary value problem \eqref{eq:Poisson equation}. \\

In a set of figures below, we present the reconstructions of a pair of scalar phantoms $W$ from the TVT of the potential vector field $\nabla W$ and from the LVT of the solenoidal vector field $\nabla^\perp W$ with various levels of additive Gaussian noise. The relative errors of these and other reconstructions in the paper are computed using the formula
\[
\operatorname{rel\_error} = \frac{||FTrue-FRec||_2}{||FTrue||_2}\times 100\, \%
\]
and are summarized in tables presented at the end of the corresponding sections. 

\begin{figure}[H]
\centering
{\includegraphics[width=0.85\textwidth]{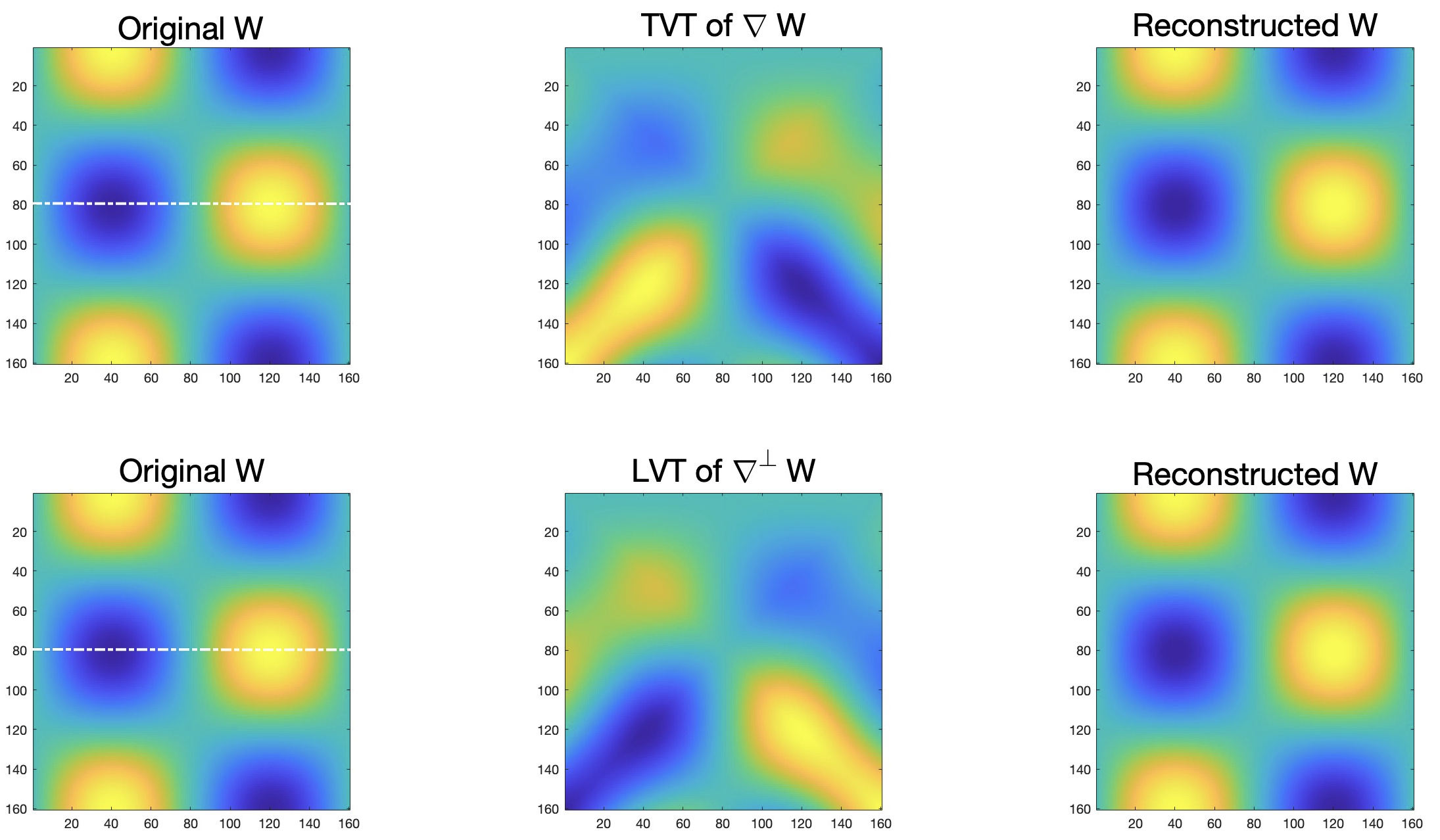}}
    \caption{Reconstruction of function $W$ from $\Tc (\nabla W)$ (top row) and $\Lc (\nabla^\perp W)$ (bottom row). }
    \label{phantom1}
  \end{figure}
\begin{figure}[H]
\centering
{\includegraphics[width=0.95\textwidth]{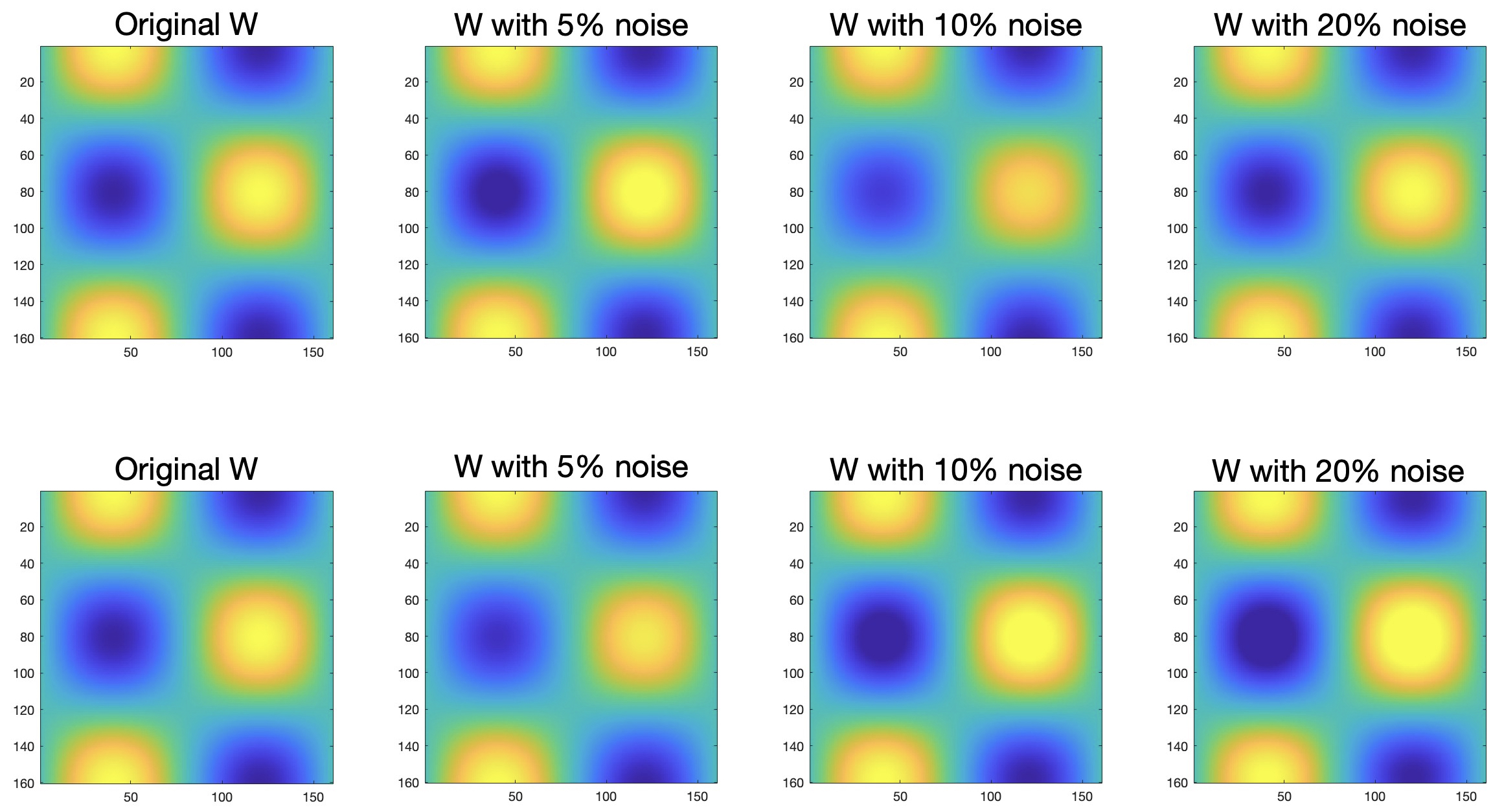}}
    \caption{Reconstructions with  $5\%$, $10 \%$, and $20 \%$ noise. The top row represents the reconstruction of $W$ from $\Tc (\nabla W)$, and the bottom row represents the reconstruction of $W$ from $\Lc (\nabla^\perp W)$.}
  \end{figure}

\begin{figure}[H]
\centering
{\includegraphics[width=0.84\textwidth]{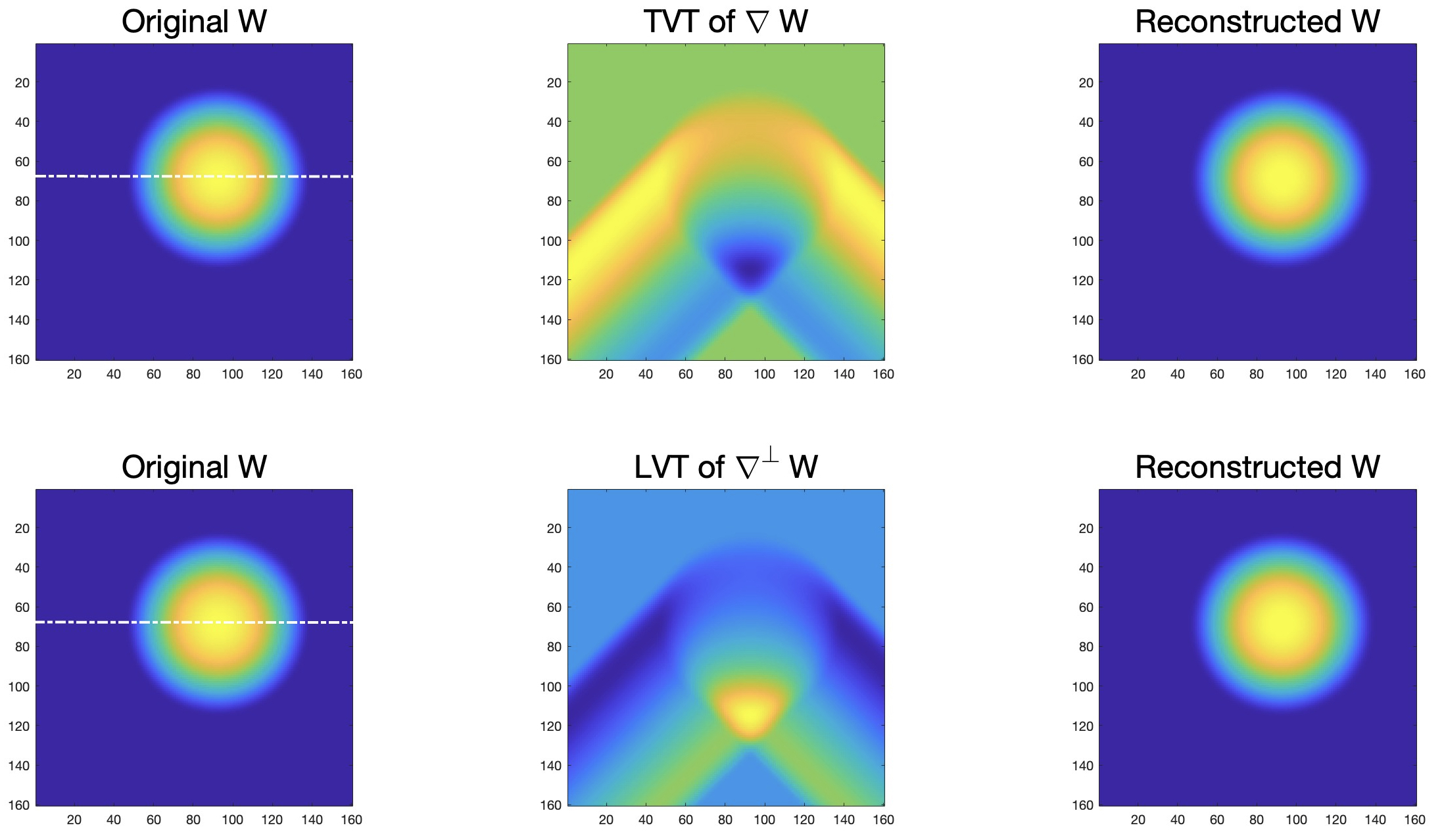}}
    \caption{Reconstruction of function $W$ from $\Tc (\nabla W)$ (top row) and $\Lc (\nabla^\perp W)$ (bottom row). }
  \end{figure}

\begin{figure}[H]
\centering
{\includegraphics[width=0.94\textwidth]{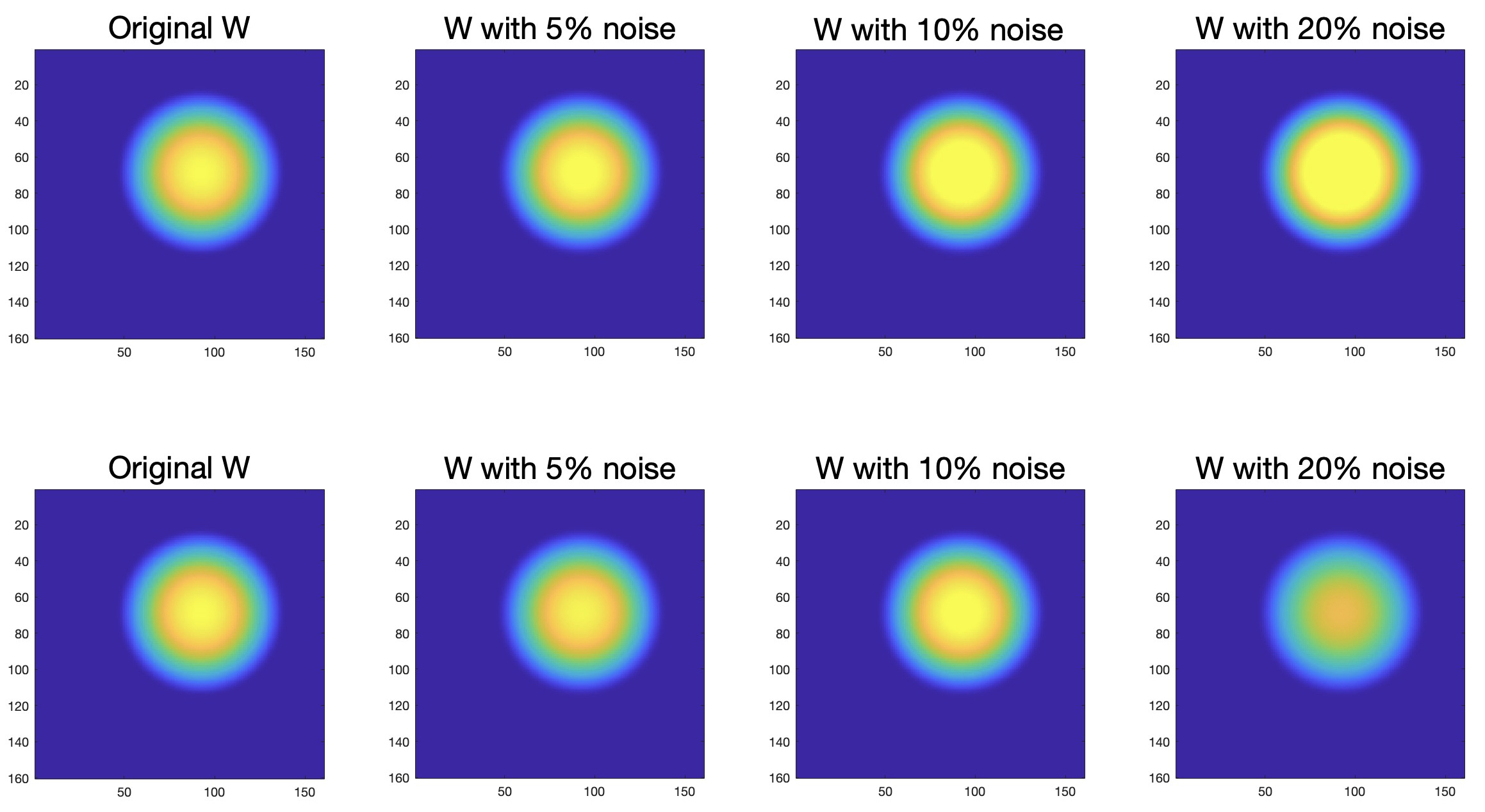}}
    \caption{Reconstructions with  $5\%$, $10 \%$ and $20 \%$ noise. The top row represents the reconstruction of $W$ from $\Tc (\nabla W)$, and the bottom row represents the reconstruction of $W$ from $\Lc (\nabla^\perp W)$.}
  \end{figure}

\begin{rem}
    The dashed white lines on the images of the original phantoms in Figure 3 (as well as in other Figures throughout the paper) are manually added to mark the lines along which the profile of the phantom is compared to that of the reconstructions (see Figure \ref{fig:prof1}).
\end{rem}

\begin{figure}[H]
\centering
{\includegraphics[width=0.95\textwidth]{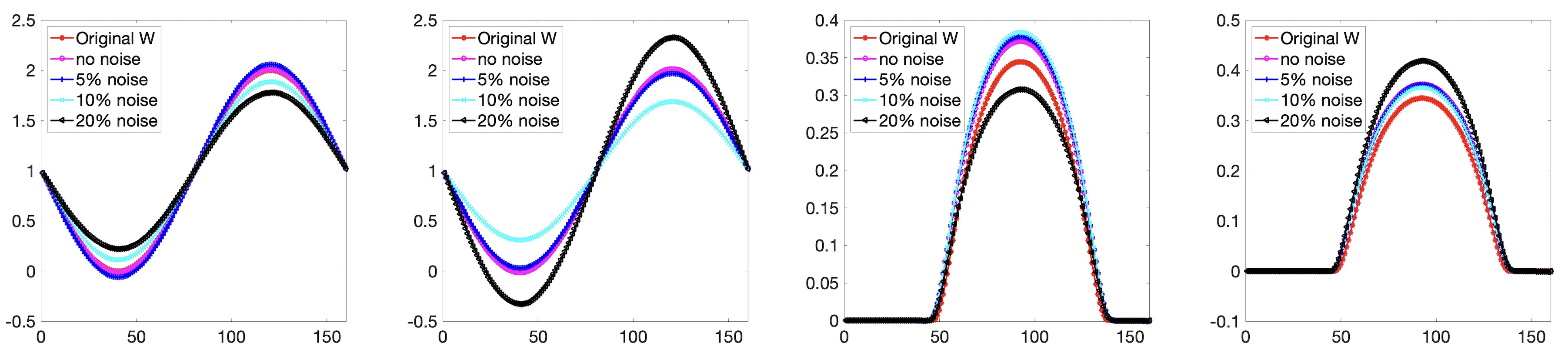}}
     \caption{Profile plots of $W$ reconstructed from $\Tc(\nabla W)$ and $\Lc(\nabla^\perp W)$ with $0\%$, $5\%$, $10 \%$, and $20 \%$ noise. The first pair of images corresponds to Phantom 1, the second pair to Phantom 2.}
   \label{fig:prof1}  
  \end{figure}
\begin{table}[h!]
\centering 
\begin{tabular}{ |p{1.6cm}||p{1cm}||p{1.6cm}|p{1.6cm}|p{1.7cm}|p{1.7cm}|}
 \hline
Phantoms & $\vf$ & No noise & 5\% noise  & 10\% noise & 20\% noise\\
 \hline
PH1   & $\nabla W$    & 0.48\% & 2.97\%  &	6.24\%	& 7.56\%\\
 \hline
PH1   &   $\nabla^\perp W$  & 0.48\% & 3.58\% & 4.74\% &	7.67\% \\
 \hline
 PH2  &  $\nabla W$ & 1.17\%	& 2.81\% &	12.11\%	& 21.51\%\\
  \hline
 PH2    & $\nabla^\perp W$ & 1.17\%	& 2.15\%	& 3.76\%& 17.95\%\\
 \hline
\end{tabular}
\caption{Relative errors of the reconstruction of $W$ from $\Tc (\nabla W)$ or $\Lc (\nabla^\perp W)$.}
\label{table:1}
\end{table}


\subsection{Recovery of a vector field from its LVT and TVT}
In this subsection, we use the combination of $\Lc\vf$ and $\Tc\vf$ to recover the full vector field $\vf = (f_1, f_2)$. We know from Theorem \ref{th: Laplacian of components of f} that $\Delta f_1$ and $\Delta f_2$ can be expressed through $\Lc \vf$ and $\Tc \vf$ as: 
\begin{align*}
   \Delta f_1(\vx) &=   D_{\vv}D_{\vu}\left[\partial_{x}  \Tc \vf(\vx) + \partial_{y}   \Lc \vf(\vx)
   \right],\\
     \Delta f_2(\vx) &=  - D_{\vv}D_{\vu}\left[\partial_{x}  \Lc \vf(\vx) - \partial_{y} \Tc \vf(\vx)
   \right].
\end{align*}
The approach to reconstructing $(f_1,f_2)$ is similar to the technique discussed in Subsection \ref{subsec:Recovery of solenoidal and potential vector field}. We start by computing the gradients $\nabla \Lc \vf = \left(\partial_{x}  \Lc \vf, \partial_{y}  \Lc \vf \right)$ and $\nabla \Tc \vf = \left(\partial_{x}  \Tc \vf, \partial_{y}  \Tc \vf \right)$. 
Using these derivatives, we get the terms $\partial_{x}  \Tc \vf(\vx) + \partial_{y}   \Lc \vf(\vx) $ and $\partial_{x}  \Lc \vf(\vx) - \partial_{y} \Tc \vf(\vx)$ appearing in the right-hand side of the above equations. Now, $\Delta f_1$ and $\Delta f_2$ can be computed by evaluating the directional derivatives as discussed in Subsection \ref{subsec:Recovery of solenoidal and potential vector field}. Then $f_1$ and $f_2$ are recovered by solving numerically the following boundary value problems:
\begin{align*}
   \Delta f_1(\vx) &=   D_{\vv}D_{\vu}\left[\partial_{x}  \Tc \vf(\vx) + \partial_{y}   \Lc \vf(\vx) 
   \right] \mbox{ in } \Omega,\ \  f_1 = 0 \mbox{ on } \partial \Omega; \\
     \Delta f_2(\vx) &=  - D_{\vv}D_{\vu}\left[\partial_{x}  \Lc \vf(\vx) - \partial_{y} \Tc \vf(\vx)
   \right] \mbox{ in } \Omega,\ \  f_2 = 0 \mbox{ on } \partial \Omega.
\end{align*}


\begin{figure}[H] 
\centering
{\includegraphics[width=0.84\textwidth]{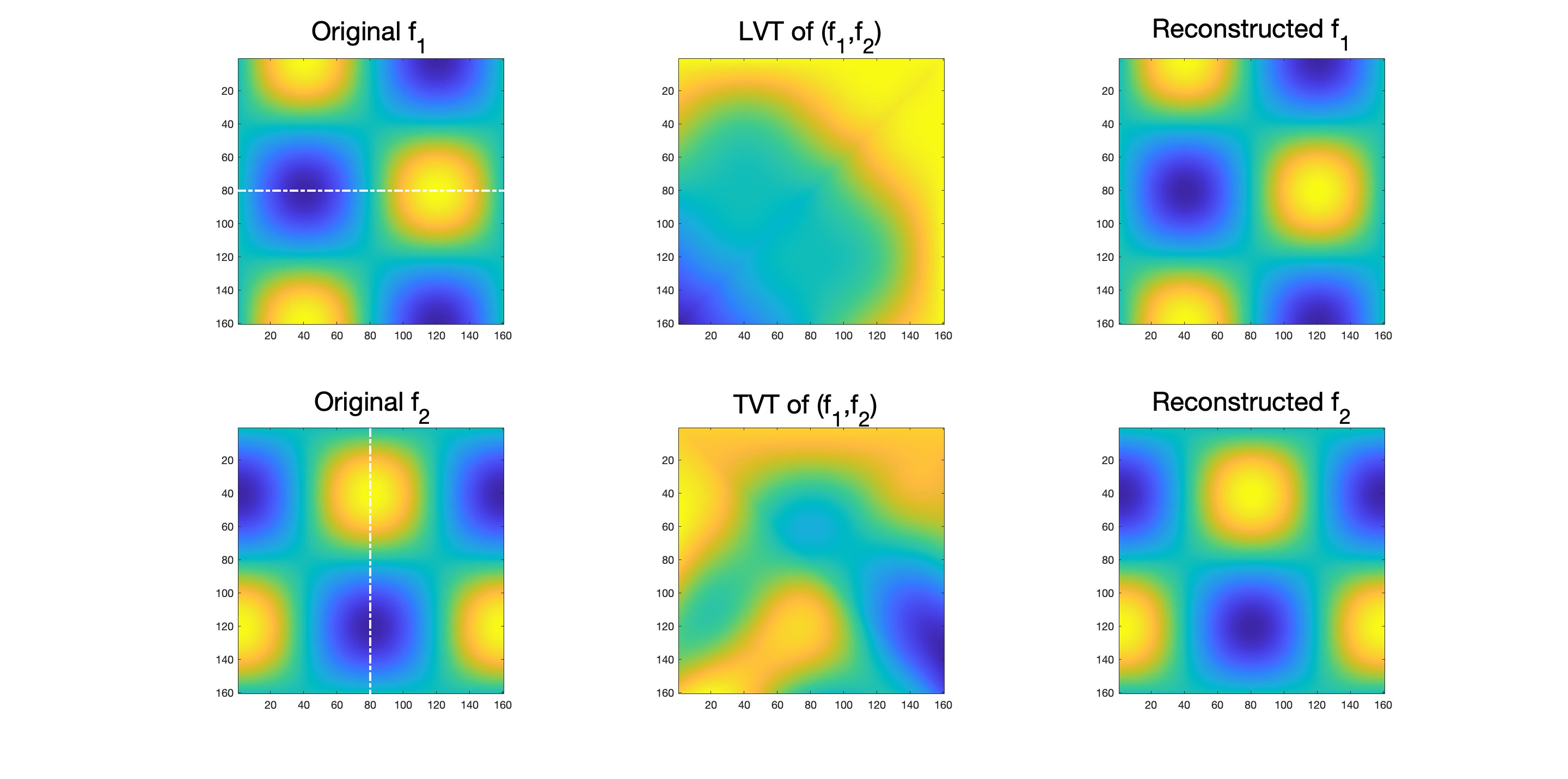}}
    \caption{Components of the original vector field $\vf$ (column 1), its transforms $\Lc\vf$ and $\Tc\vf$ (column 2), and the reconstructed components of $\vf$ (column 3).}
  \end{figure}
\begin{figure}[H]
\centering
{\includegraphics[width=0.91\textwidth]{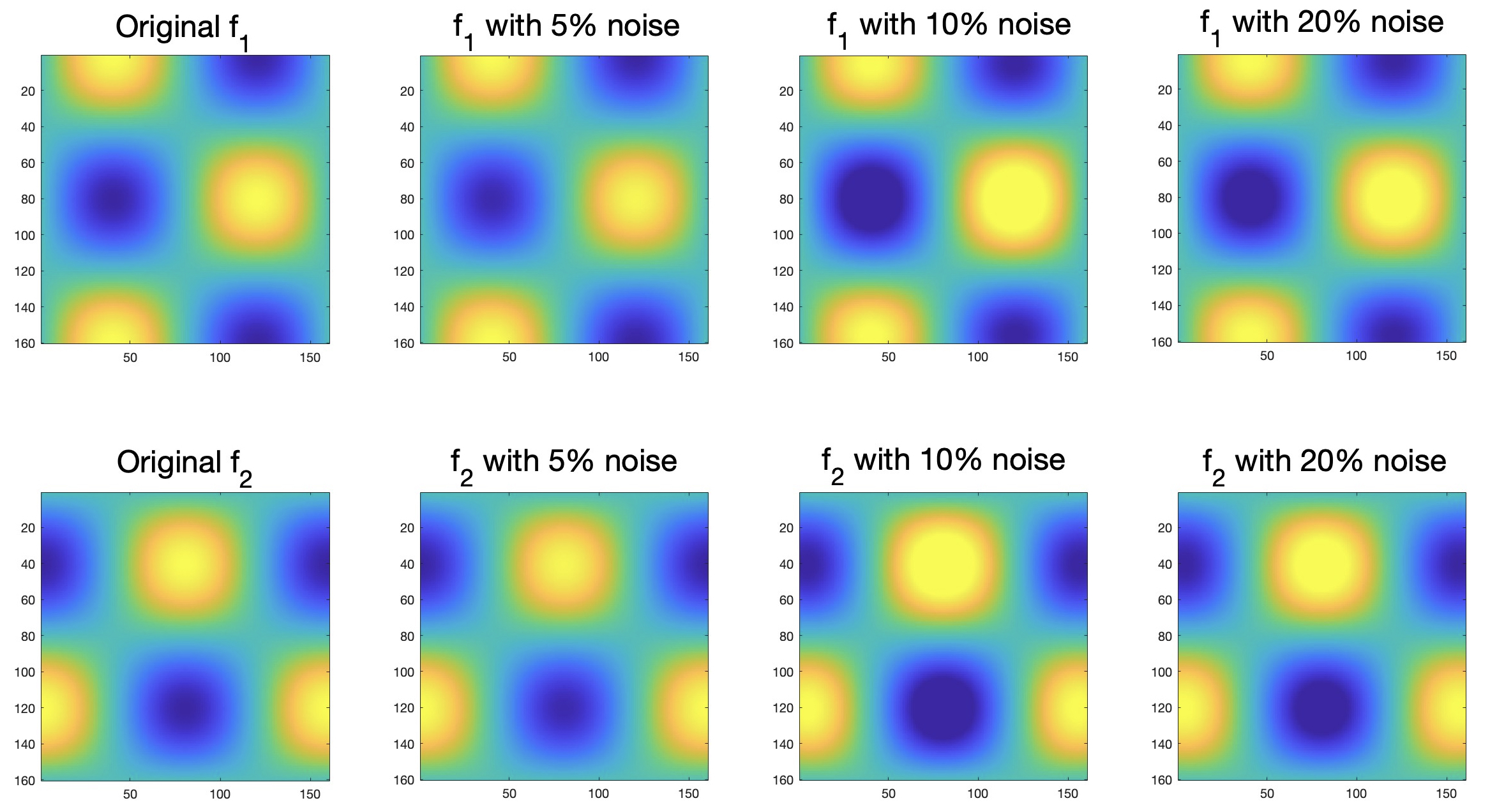}}
     \caption{Reconstructions of components of $\vf$ using $\Lc\vf$ and $\Tc\vf$ with  $5\%$, $10 \%$, and $20 \%$ noise.}
  \end{figure}


\begin{figure}[H]
\centering
{\includegraphics[width=0.84\textwidth]{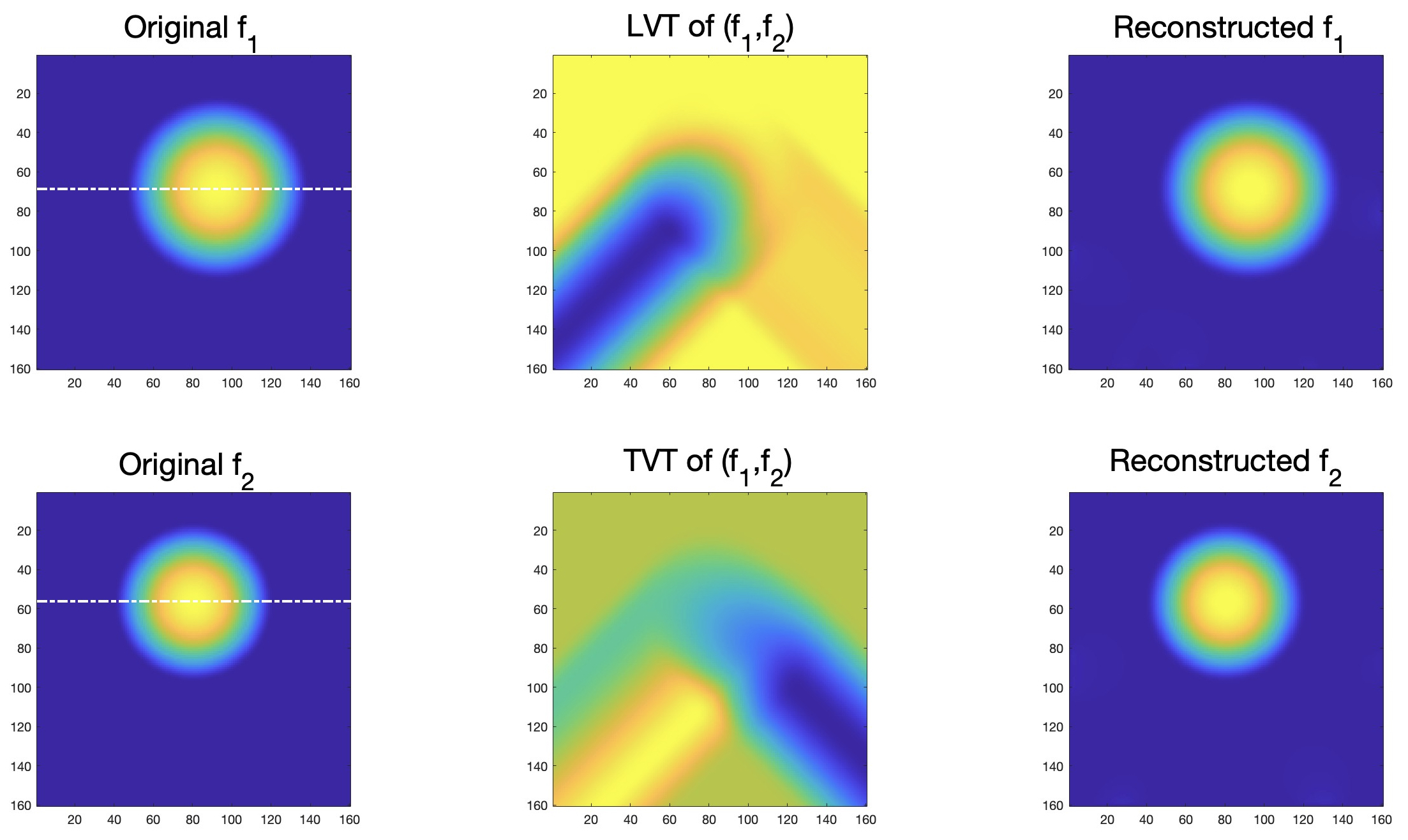}}
      \caption{Components of the original vector field $\vf$ (column 1), its transforms $\Lc\vf$ and $\Tc\vf$ (column 2), and the reconstructed components of $\vf$ (column 3).}
    \label{phantom2}
  \end{figure}
\begin{figure}[H]
\centering
{\includegraphics[width=0.92\textwidth]{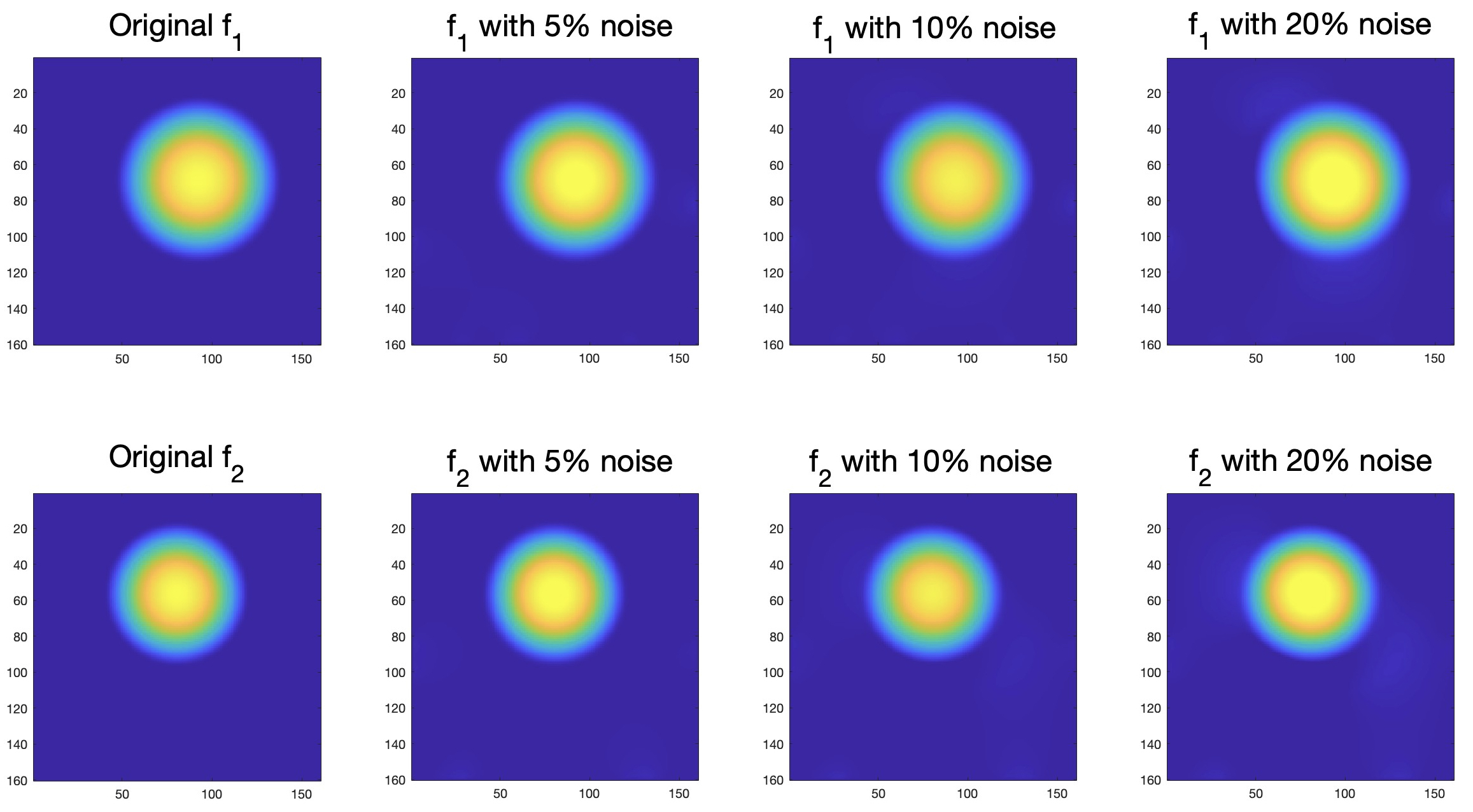}}
     \caption{Reconstructions of components of $\vf$ using $\Lc\vf$ and $\Tc\vf$ with  $5\%$, $10 \%$, and $20 \%$ noise.}
    \label{phantom2_N}
  \end{figure}

\begin{figure}[H]
\centering
{\includegraphics[width=0.84\textwidth]{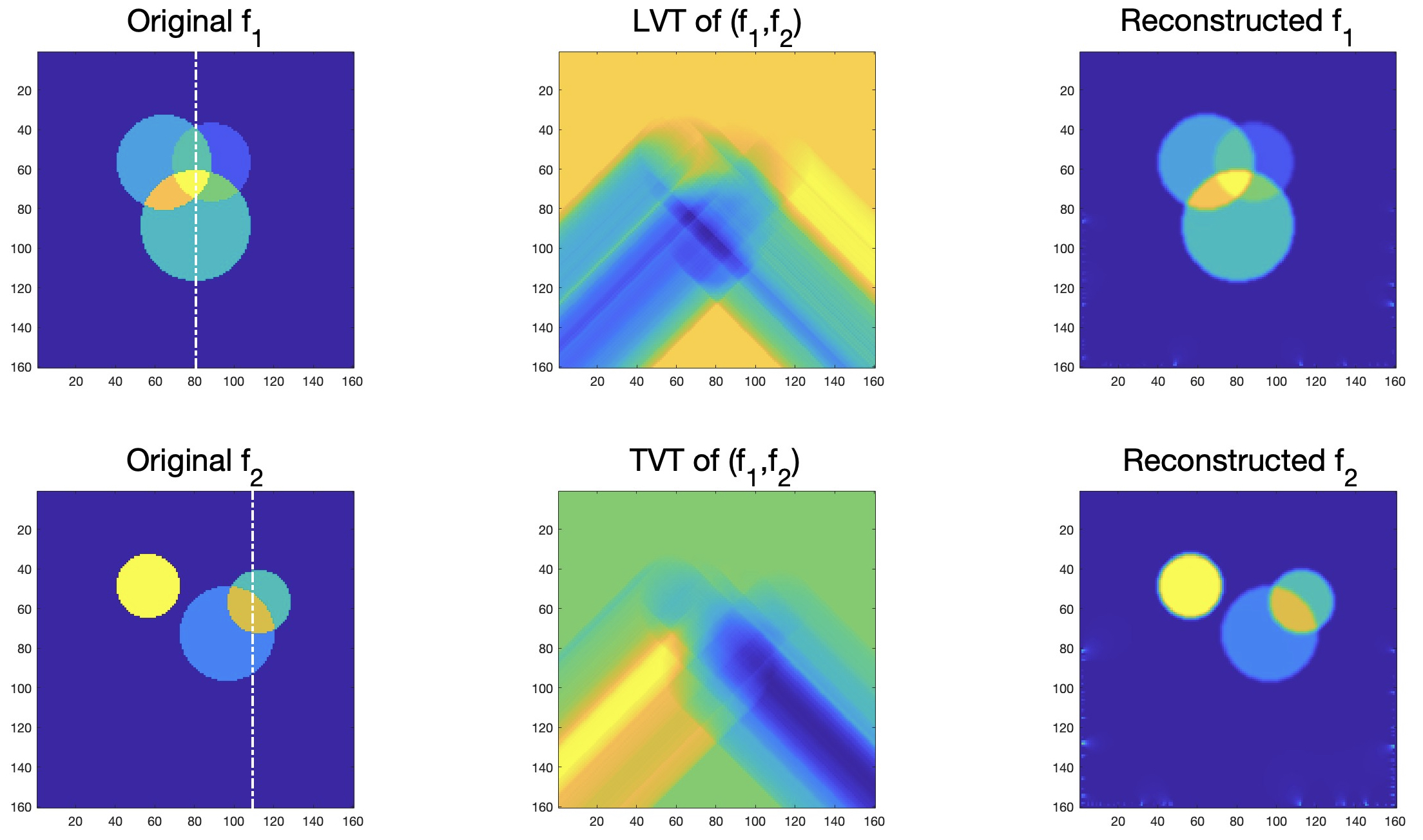}}
      \caption{Components of the original vector field $\vf$ (column 1), its transforms $\Lc\vf$ and $\Tc\vf$ (column 2), and the reconstructed components of $\vf$ (column 3).}
    \label{phantom3}
  \end{figure}
\begin{figure}[H]
\centering
{\includegraphics[width=0.91\textwidth]{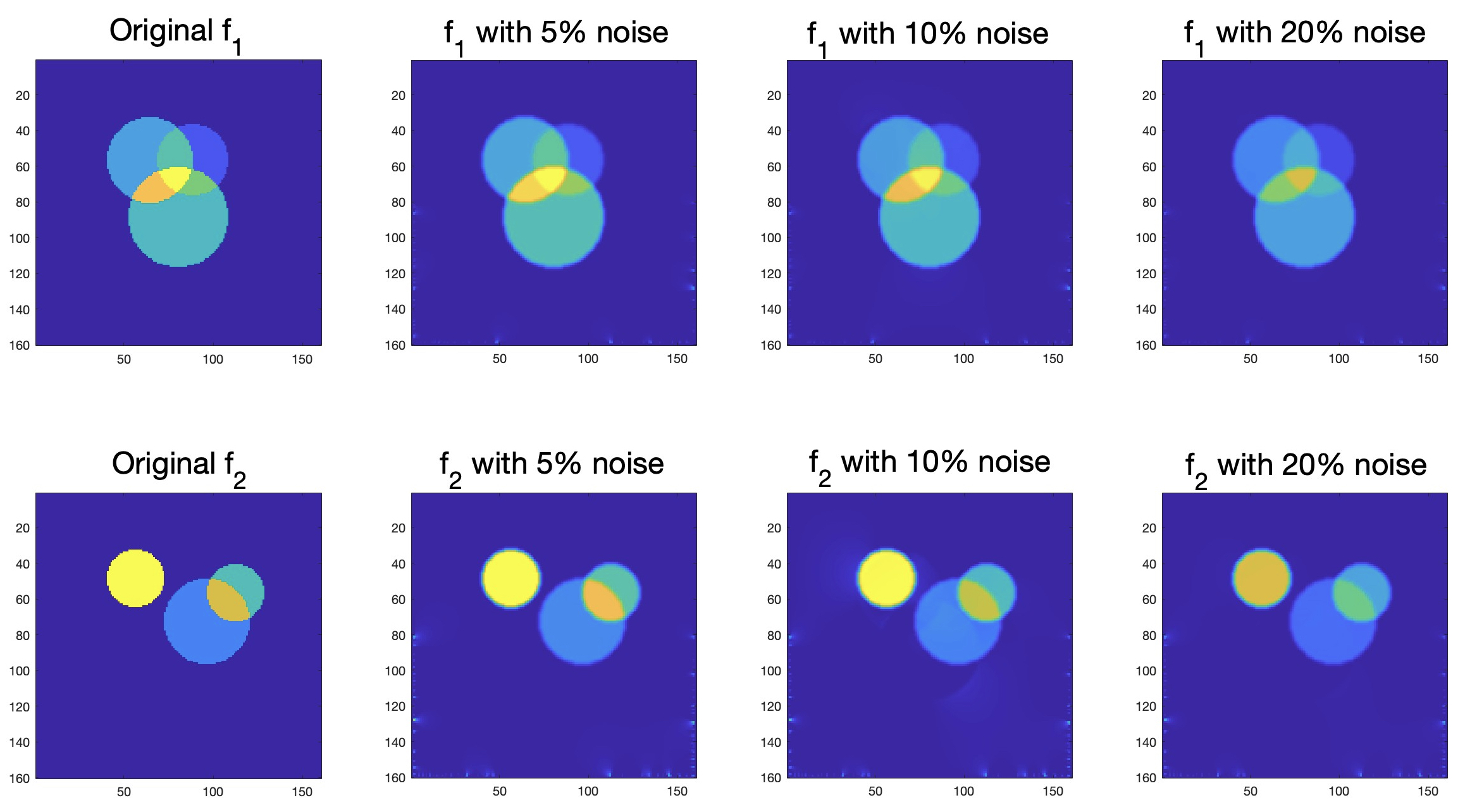}}
     \caption{Reconstructions of components of $\vf$ using $\Lc\vf$ and $\Tc\vf$ with  $5\%$, $10 \%$, and $20 \%$ noise.}
    \label{phantom3_N}
  \end{figure}

\begin{figure}[H]
\centering
{\includegraphics[width=0.95\textwidth]{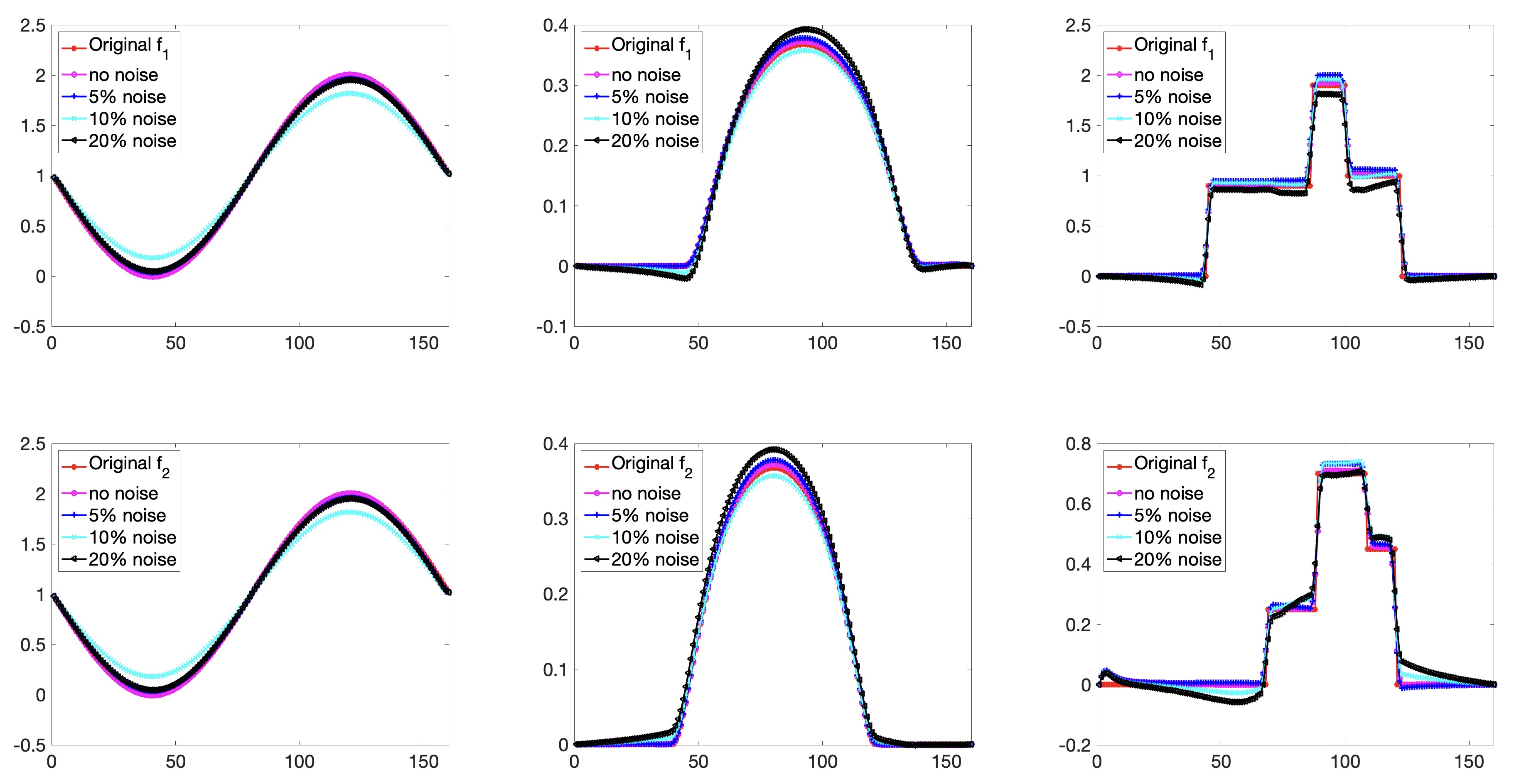}}
     \caption{Profile plots of $f_1$ and $f_2$ reconstructed from $\Lc\vf$ and $\Tc\vf$ with $0\%$, $5\%$, $10 \%$, and $20 \%$ noise. Plots in $j$-th column correspond to Phantom $j$, $j=1,2,3$.}
  \end{figure}
  
 
\begin{table}[h!]
\centering 
\begin{tabular}{ |p{1.6 cm}||p{0.6cm}||p{1.6cm}|p{1.6cm}|p{1.7cm}|p{1.7cm}|  }
 \hline
Phantoms & $\vf$ & No noise & 5\% noise  & 10\% noise & 20\% noise\\
 \hline
PH1   & $f_1$    & 0.96\% &	1.71\%	& 6.26\%	& 9.76\%\\
 \hline
PH1   &   $f_2$  & 0.66\% &	1.58\%	& 6.27\%	& 9.77\%\\
 \hline
 PH2  &  $f_1$ & 1.46\%	& 3.00\%	& 3.78\%	& 8.21\%\\
  \hline
 PH2    & $f_2$ & 1.34\%	& 2.88\%	& 3.92\%	& 8.20\%\\
 \hline 
 PH3  &  $f_1$ & 3.67\%	& 3.86\%	& 6.53\%	& 14.40\%\\
  \hline
 PH3    & $f_2$ & 6.87\% &	7.14\%	& 7.74\%	& 20.3\%\\
 \hline 
\end{tabular}
 \caption{Relative errors of the reconstructions of $f_1$ and $f_2$ from $\Lc\vf$ and $\Tc\vf$.}
\label{table:2}
\end{table}
  

\subsection{Recovery of a vector field from its LVT and LVT1, or TVT and TVT1}\label{sec:rec-moments}
This subsection focuses on combining the V-line transforms (longitudinal and transverse) and their first moments to recover the full vector field $\vf =(f_1, f_2)$. As we described in Theorem \ref{th: reconstruction with integral moments},  this involves the inversion of the signed V-line transform. We discuss below the implementation of the reconstruction process for $f_1$ from $\Lc \vf$ and $\Ic \vf$. The reconstruction of $f_2$ from those transforms, as well as the reconstructions of $f_1$ and $f_2$ from the transverse data ($\Tc \vf$ and $\Jc \vf$) follow similarly. 

Recall from Theorem \ref{th: reconstruction with integral moments} that $f_1$ is given by:
\begin{align*}
   f_1(\vx)=  \frac{1}{\|\vv - \vu\|}D_{\vv}D_{\vu}\int_0^\infty\left\{ \frac{\partial \Ic \vf}{\partial x} + u_2 \Xc_{\vu}^1( \operatorname{curl}  \vf) - v_2 \Xc_{\vv}^1 ( \operatorname{curl} \vf) \right\}(\vx + t \vw) dt.
\end{align*}

As a first step, we compute two directional derivatives (as discussed in Subsection \ref{subsec:Recovery of solenoidal and potential vector field}) of $\Lc \vf$ to generate $\operatorname{curl} \vf$ (see  
formula \eqref{eq:curl-explicit}). Then we apply the procedure discussed in Subsection \ref{subsec:data formulation} to compute the first moment divergent beam transforms $\Xc_{\vu}^1( \operatorname{curl} \vf)$  and $\Xc_{\vv}^1 ( \operatorname{curl} \vf)$. Using the Matlab built-in function \texttt{gradient}, we find the partial derivatives $\left(\frac{\partial \Ic \vf}{\partial x}, \frac{\partial \Ic \vf}{\partial y} \right) $ of the first moment data $\Ic \vf$. By combining these quantities, we evaluate the integrand in the formula for $f_1$ quoted above, i.e.
$$
I\doteq\frac{\partial \Ic \vf}{\partial x} + u_2 \Xc_{\vu}^1( \operatorname{curl}  \vf) - v_2 \Xc_{\vv}^1 ( \operatorname{curl} \vf). 
$$
Notice, that the aforementioned integral itself is nothing but the divergent beam transform of the evaluated function $I$ along the direction $\vw = \vv -\vu$, which we already know how to compute (see Subsection \ref{subsec:data formulation}). Finally, we apply the directional derivatives $D_{\vu}$ and $D_{\vv}$ to the result obtained after integration to generate $ \|\vv - \vu\| f_1 =\sqrt{2} f_1$.   

\begin{rem}
    It was shown in \cite{Gaik_Mohammad_Rohit} that function $I$ coincides with the signed V-line transform (SVL) of $f_1$, i.e.
$$
I=\frac{\partial \Ic \vf}{\partial x} + u_2 \Xc_{\vu}^1( \operatorname{curl}  \vf) - v_2 \Xc_{\vv}^1 ( \operatorname{curl} \vf)=\Xc_{\vu}f_1-\Xc_{\vv}f_1.
$$   
Therefore, as an intermediate step of our procedure we recover SVL of $f_1$, and the follow-up steps are ensuing the inversion of SVL.
\end{rem}

\begin{rem}
    The SVL inversion procedure used here was developed in \cite{amb-lat_2019} and requires
    data along V-lines with vertices in a larger set than the support of the image function (recall Remark \ref{rem:support}). Therefore, the algorithm is tested on a ``truncated'' version of Phantom 1, and the original Phantoms 2 and 3.
\end{rem}
Some of the reconstructed images presented below include artifacts that spread along the divergent beams involved in the associated inversion formulas. Such artifacts are typical for the numerical inversions of various V-line transforms (e.g. see \cite{amb-book, amb-lat_2019, Florescu-Markel-Schotland, Gouia_Amb_V-line, Sherson}) and can be explained by microlocal properties of the divergent beam transform. Of particular importance here is the relation between the wavefront sets of a distribution $h\in\mathcal{D}'(\mathbb{R}^2)$ and its divergent beam transform $\Xc_{\vgamma}h$. It is known (e.g. see \cite{amb-book, Sherson}) that 
\begin{equation}
    \operatorname{WF}(\mathcal{X}_{\vgamma}h)\subseteq
    \operatorname{WF}(h) \cup
   \big\{\big(\vx-t\vgamma,\vxi\big)\,\big|\; \big(\vx,\vxi\big)\in \operatorname{WF}(h),\, \vxi\in \vgamma^\perp,\, t>0\big\}.
\end{equation}
In other words, in addition to the true singularities (e.g. jump discontinuities) of $h$, its divergent beam transform data may also include a set of additional singularities, which start at the points where $h$ has singularities
in the direction $\vgamma^\perp$ and propagate in the direction $-\vgamma$.
This implies that in the images reconstructed from $\Xc_{\vgamma}h$, artifacts may appear along rays in the direction of $-\vgamma$ that originate at and are tangent to the boundary of some feature (non-smoothness) in $h$.

\vspace{2mm}

In the numerical reconstructions presented below, we recover the scalar components $f_1$ and $f_2$ of the vector field $\vf$ using the inversion formulas (\ref{eq:recovery f1}) and (\ref{eq: recovery f2}) from Theorem \ref{th: reconstruction with integral moments}. In both cases we take a divergent beam transform $\Xc_{\vgamma}h$ of some processed data $h$, where $-\vgamma=(\vu-\vv)/||\vu-\vv||=(1,0)$, followed by two directional derivatives. Notice, that the processed data set $h$ is different in (\ref{eq:recovery f1}) and (\ref{eq: recovery f2}), and it includes a different set of singularities. In the case of (\ref{eq:recovery f1}), a portion of the singularities are due to $\frac{\partial \Ic \vf}{\partial x_1}$, while in (\ref{eq: recovery f2}) a similar portion is due to $\frac{\partial \Ic \vf}{\partial x_2}$.  

In the phantoms depicted in Figures \ref{L_Moments_phantom1} and \ref{L_Moments_phantom1_N}, the singularities of $\frac{\partial \Ic \vf}{\partial x_1}$ are vertical (thus, $\vgamma$ is not tangent to them), and the reconstructions of $f_1$ are free of horizontal artifacts. At the same time, the singularities of $\frac{\partial \Ic \vf}{\partial x_2}$ are horizontal (thus, $\vgamma$ is tangent to them), leading to strong horizontal artifacts in the reconstructions of $f_2$. Similar artifacts can also be observed in Figures \ref{L_Moments_phantom3} and \ref{L_Moments_phantom3_N}, which involve another set of piecewise constant images.  

Another portion of singularities (in the processed data set $h$ used in formulas (\ref{eq:recovery f1}) and (\ref{eq: recovery f2})) comes from  $\Xc_{\vu}^1( \operatorname{curl} \vf)$ and $\Xc_{\vv}^1( \operatorname{curl} \vf)$. These data sets, in their own right, include singularities along rays in the directions $-\vu$ and $-\vv$  that originate at and are tangent to the boundary of some feature (non-smoothness) in $\operatorname{curl} \vf$. The ``diagonal'' artifacts corresponding to these singularities can be observed in Figures \ref{L_Moments_phantom1}-\ref{L_Moments_phantom3_N}.



\subsubsection{Recovery of a vector field from its LVT and LVT1}
\begin{figure}[H]
\centering
{\includegraphics[width=0.84\textwidth]{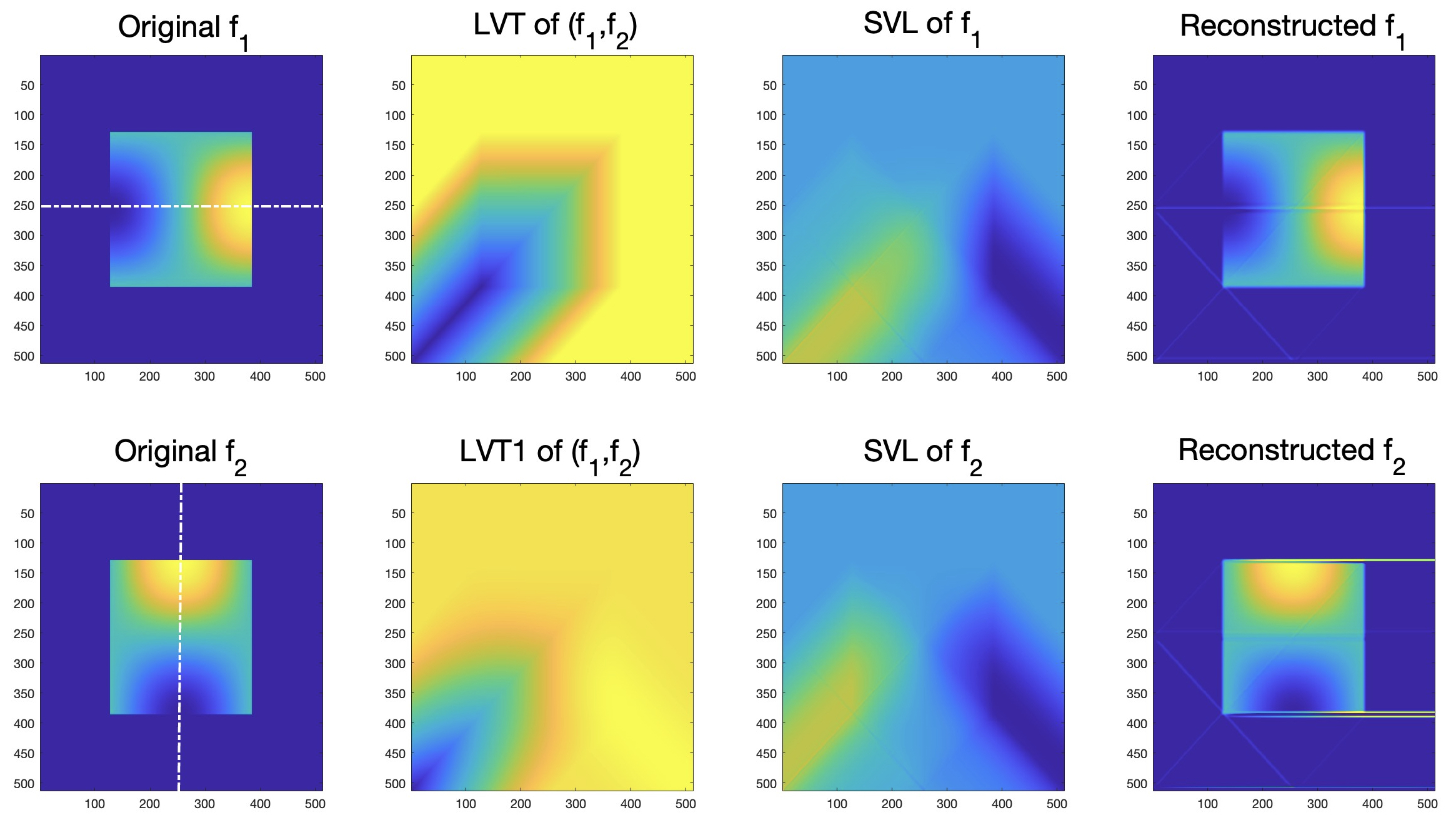}}
  \caption{Components of the original field $\vf$ (column 1), $\Lc \vf$ and $\Ic \vf$ (column 2), signed V-line transform of the components (column 3), and reconstructed components of $\vf$ (column 4).}
    \label{L_Moments_phantom1}
  \end{figure}
\begin{figure}[H]
\centering
{\includegraphics[width=0.84\textwidth]{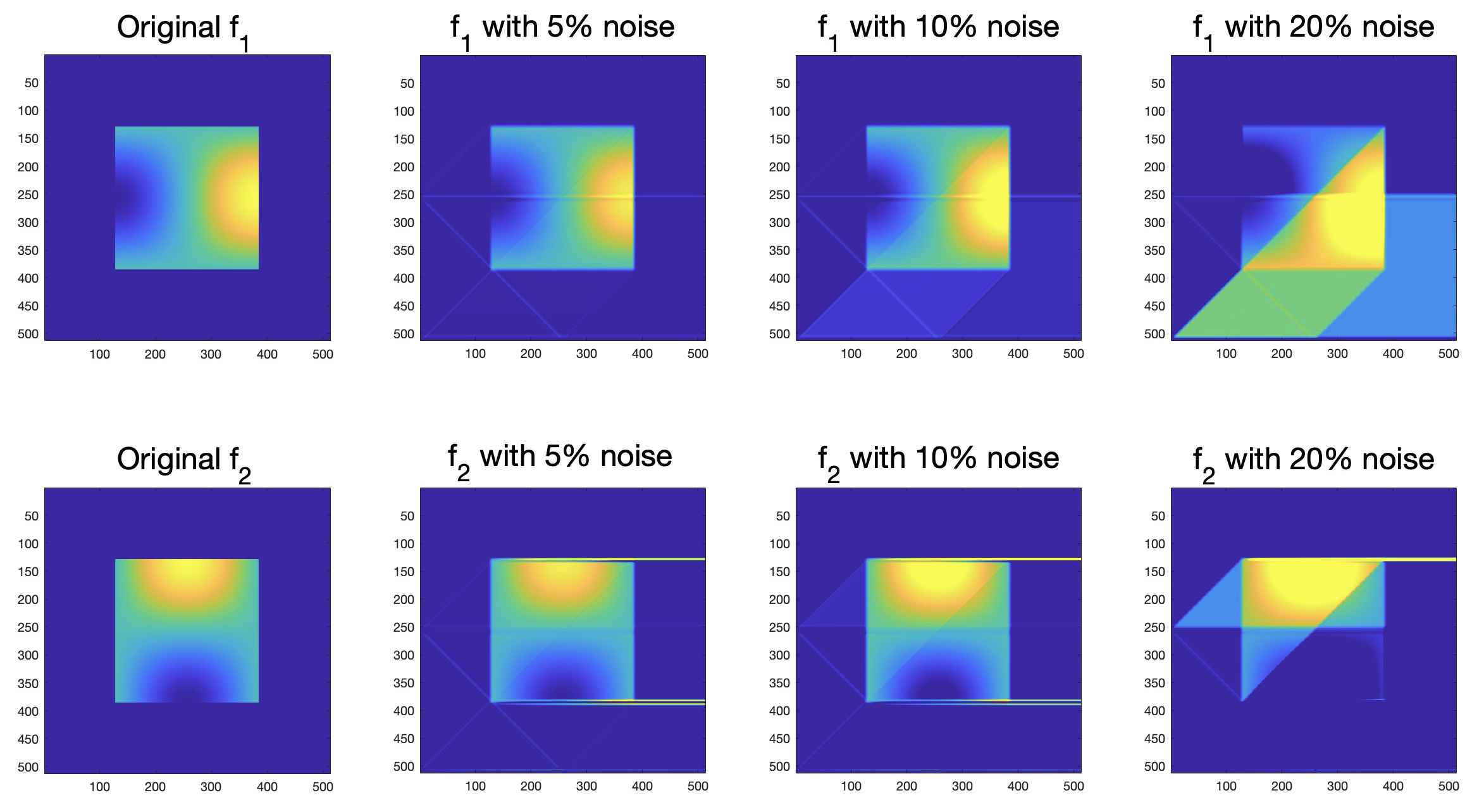}}
      \caption{Reconstructions of components of $\vf$ using $\Lc\vf$ and $\Ic\vf$ with  $5\%$, $10 \%$, and $20 \%$ noise.}
    \label{L_Moments_phantom1_N}
  \end{figure}

\begin{figure}[H]
\centering
{\includegraphics[width=0.9\textwidth]{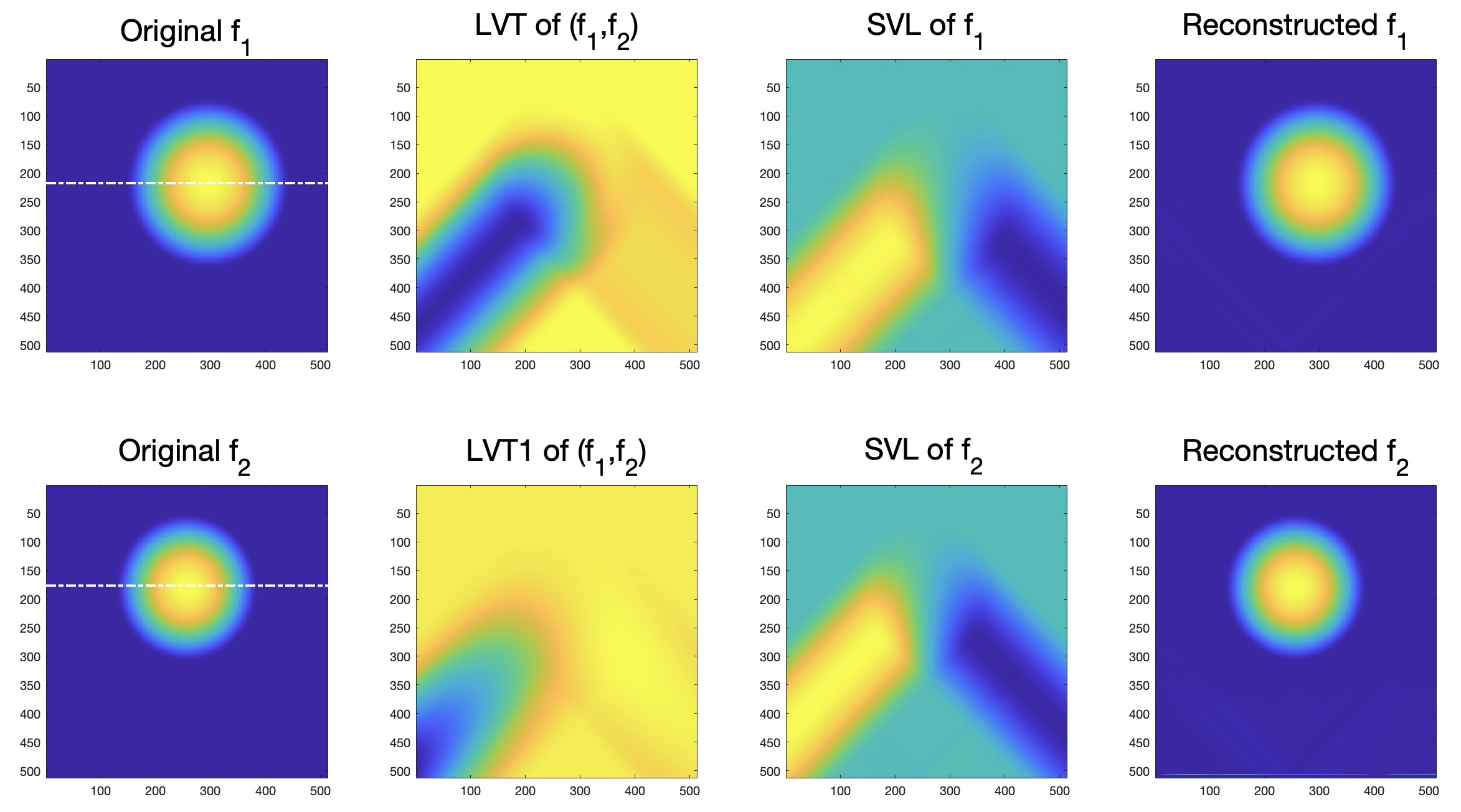}}
   \caption{Components of the original field $\vf$ (column 1), $\Lc \vf$ and $\Ic \vf$ (column 2), signed V-line transform of the components (column 3), and reconstructed components of $\vf$ (column 4).}
    \label{L_Moments_phantom2}
  \end{figure}
\begin{figure}[H]
\centering
{\includegraphics[width=0.9\textwidth]{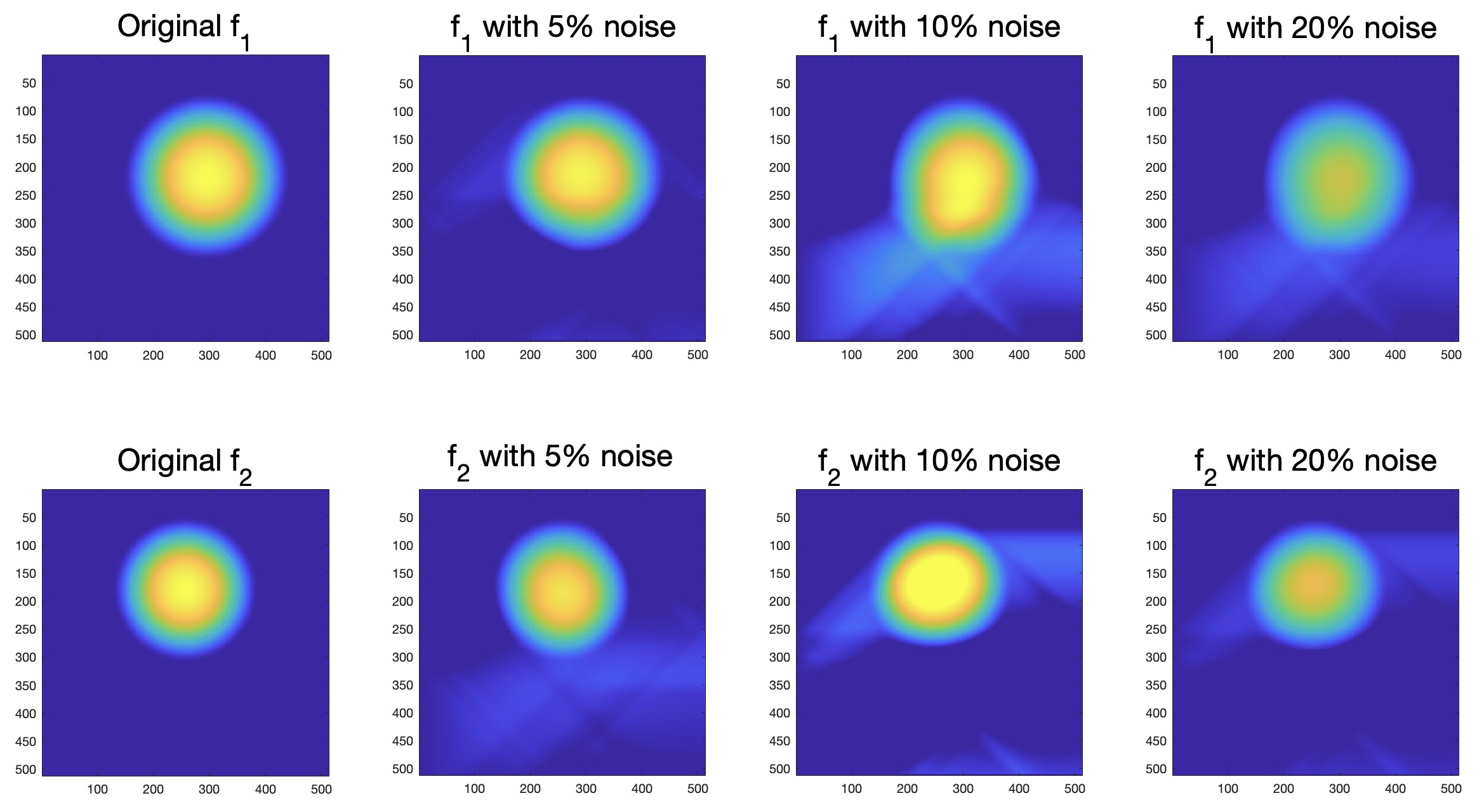}}
 \caption{Reconstructions of components of $\vf$ using $\Lc\vf$ and $\Ic\vf$ with  $5\%$, $10 \%$, and $20 \%$ noise.}
    \label{L_Moments_phantom2_N}
  \end{figure}
\begin{figure}[H]
\centering
{\includegraphics[width=0.9\textwidth]{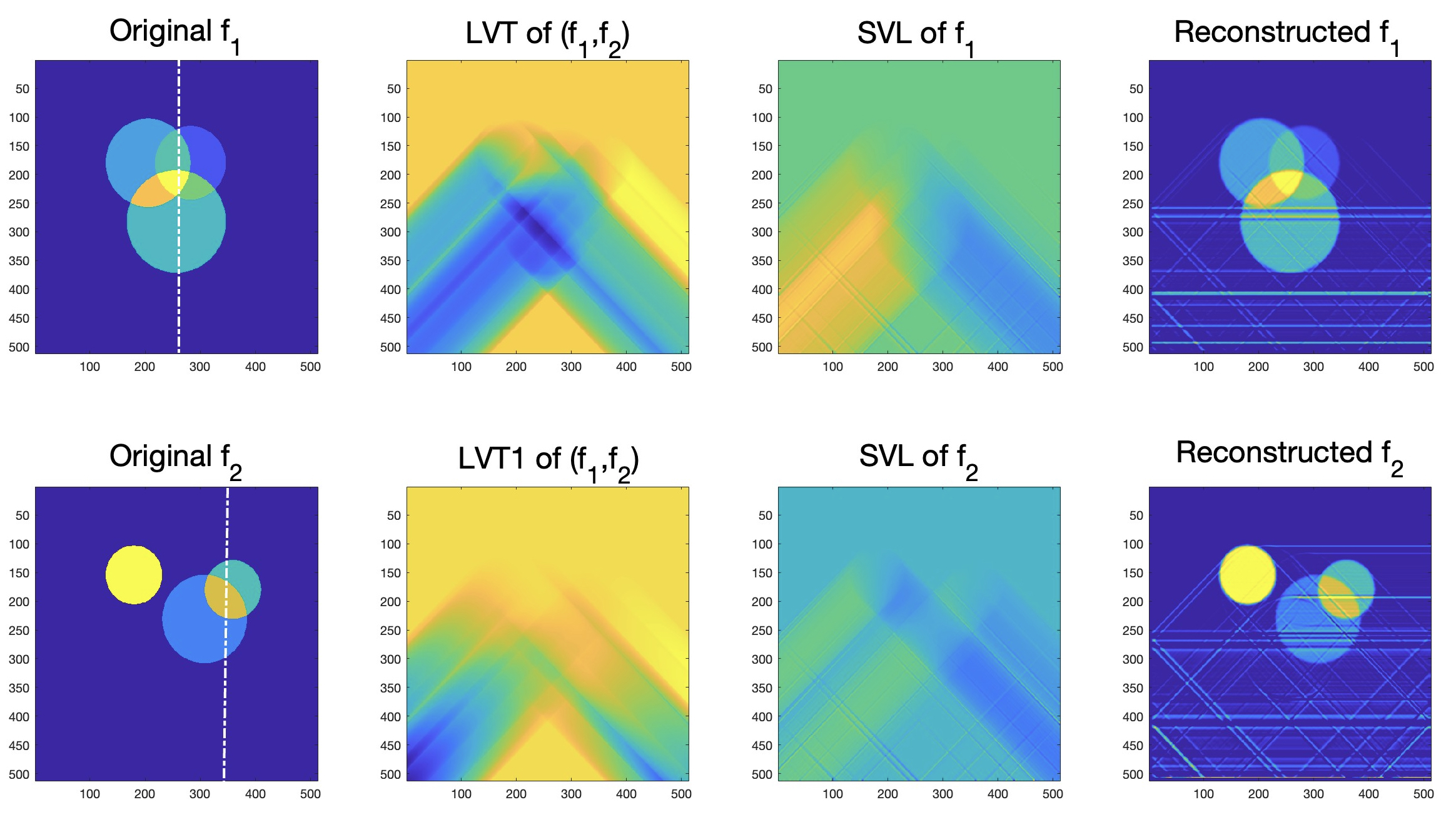}}
     \caption{Components of the original field $\vf$ (column 1), $\Lc \vf$ and $\Ic \vf$ (column 2), signed V-line transform of the components (column 3), and reconstructed components of $\vf$ (column 4).}
    \label{L_Moments_phantom3}
  \end{figure}
\begin{figure}[H]
\centering
{\includegraphics[width=0.9\textwidth]{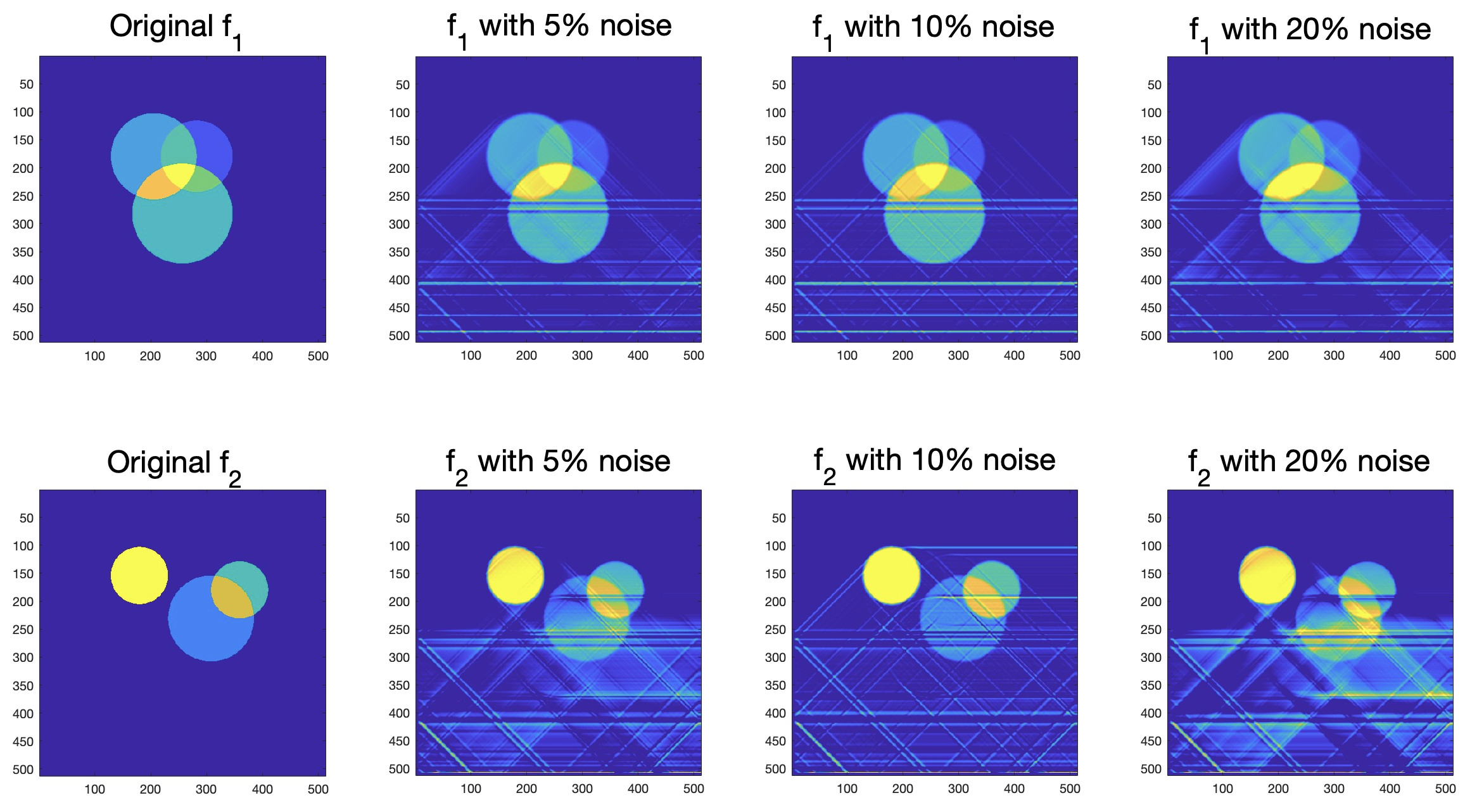}}
   \caption{Reconstructions of components of $\vf$ using $\Lc\vf$ and $\Ic\vf$ with  $5\%$, $10 \%$, and $20 \%$ noise.}
    \label{L_Moments_phantom3_N}
  \end{figure}
  \begin{figure}[H]
\centering
{\includegraphics[width=0.94\textwidth]{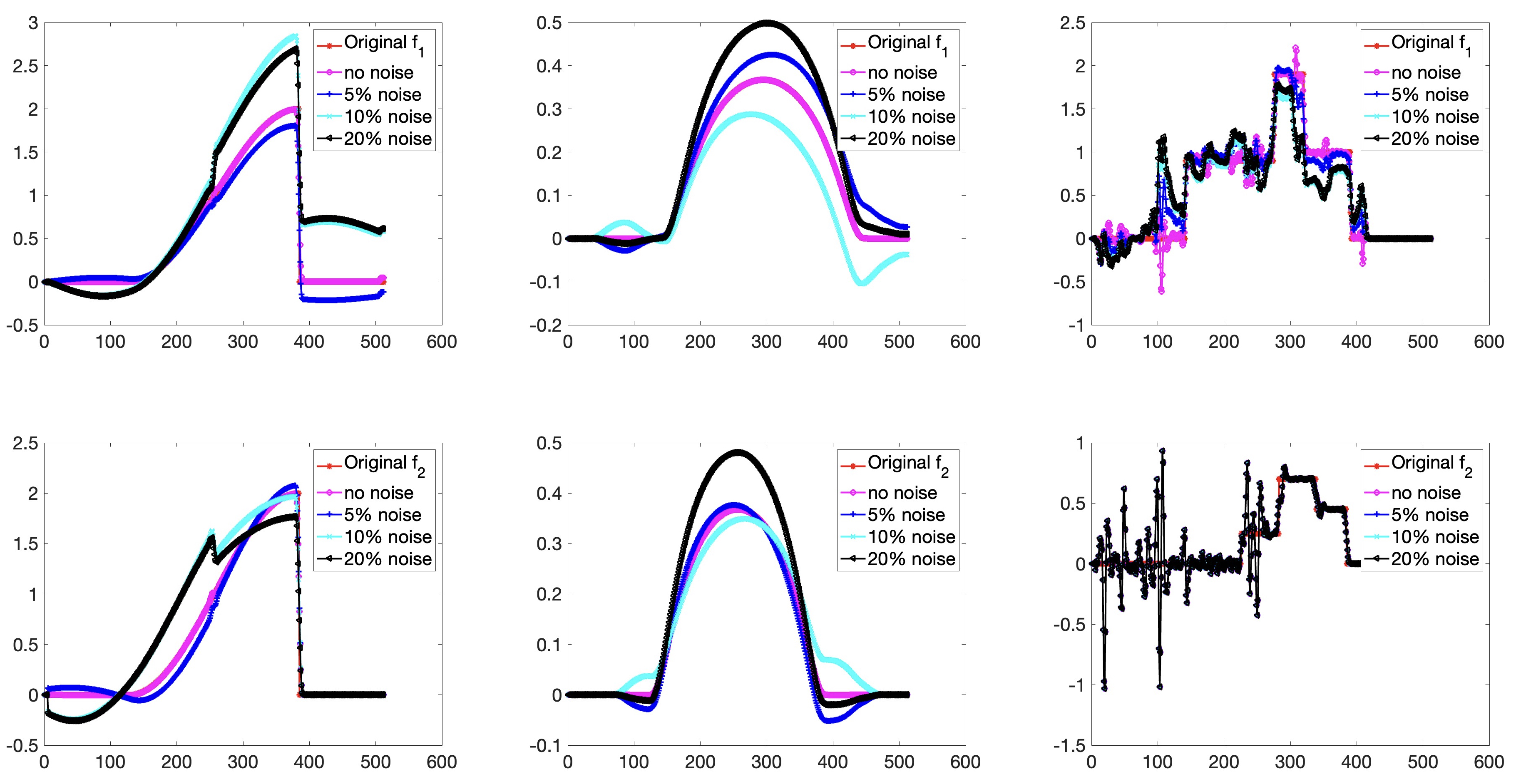}}
     \caption{Profile plots of $f_1$ and $f_2$ reconstructed from $\Lc(\vf)$ and $\Ic(\vf)$ with $0\%$, $5\%$, $10 \%$, and $20 \%$ noise. Plots in $j$-th column correspond to Phantom $j$, $j=1,2,3$.}
  \end{figure}
\begin{table}[h!]
\centering 
\begin{tabular}{ |p{1.6 cm}||p{0.6cm}||p{1.6cm}|p{1.6cm}|p{1.7cm}|p{1.7cm}|  }
 \hline
Phantoms & $\vf$ & No noise & 5\% noise  & 10\% noise & 20\% noise\\
 \hline
PH1   & $f_1$    & 6.52\%	& 7.73\%	& 16.03\% & 170.01\%\\
 \hline
PH1   &   $f_2$  & 52.41\%	& 60.58\% &	45.48\%  & 	491.23\%\\
 \hline
 PH2  &  $f_1$ & 1.06\% & 8.42\%	& 20.84\% &	72.80\%\\
  \hline
 PH2    & $f_2$ & 2.05\%	& 9.08\% & 30.19\% &	108.94\%\\
 \hline 
 PH3  &  $f_1$ & 48.76\% & 45.33\% & 51.85\% & 65.45\%\\
  \hline
 PH3    & $f_2$ & 47.06\% &	51.08\% &	52.94\% &	166.51\%\\
 \hline 
\end{tabular}
\caption{Relative errors of the reconstructions of $f_1$ and $f_2$ from $\Lc\vf$ and $\Ic\vf$.}
\label{table:3}
\end{table}

\subsubsection{Recovery of a vector field from its TVT and TVT1}
\begin{figure}[H]
\centering
{\includegraphics[width=0.84\textwidth]{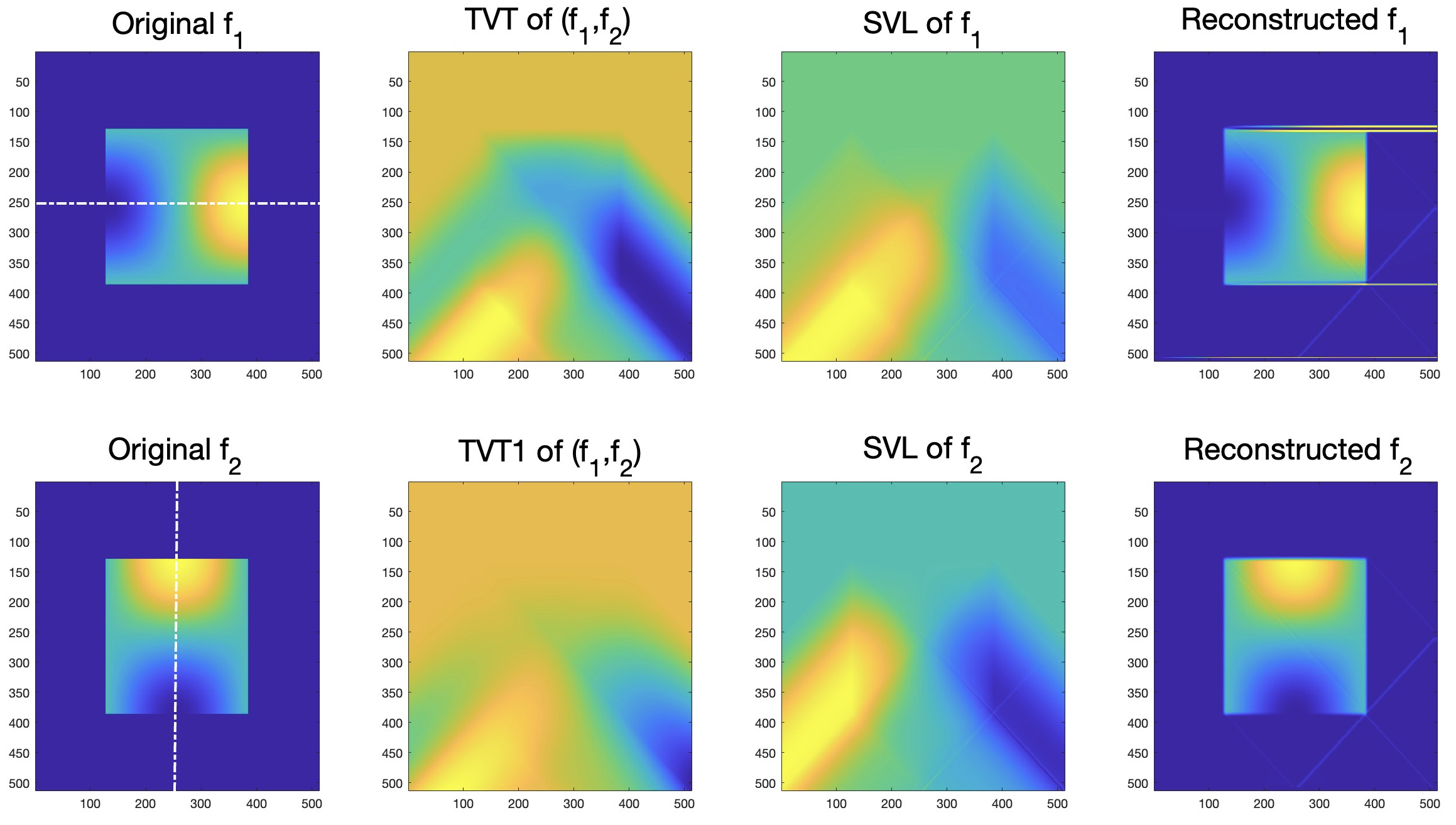}}
    \caption{Components of the original field $\vf$ (column 1), $\Tc \vf$ and $\Jc \vf$ (column 2), signed V-line transform of the components (column 3), and reconstructed components of $\vf$ (column 4).}
    \label{T_Moments_phantom1}
  \end{figure}
\begin{figure}[H]
\centering
{\includegraphics[width=0.84\textwidth]{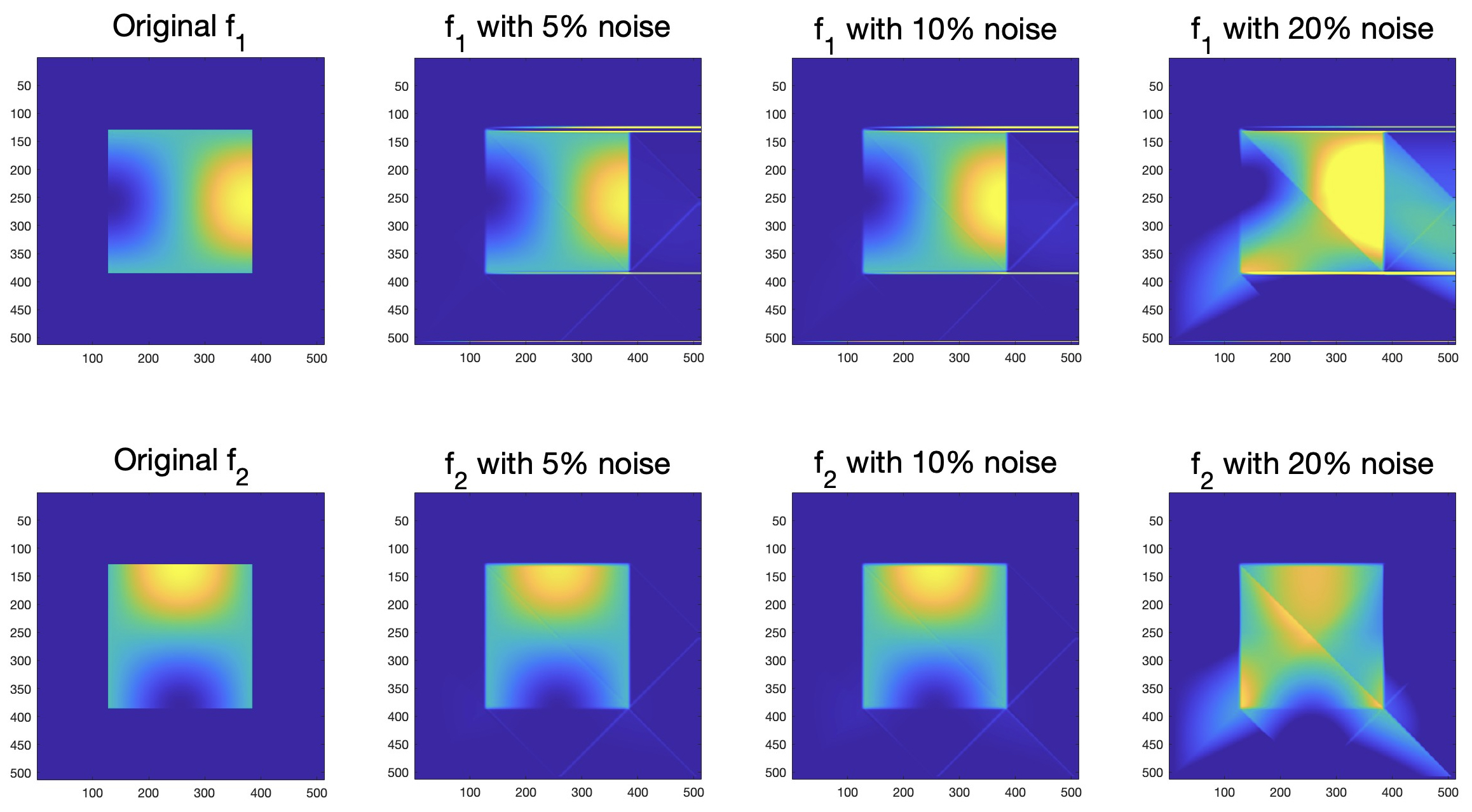}}
     \caption{Reconstructions of components of $\vf$ using $\Tc\vf$ and $\Jc\vf$ with  $5\%$, $10 \%$, and $20 \%$ noise.}
    \label{T_Moments_phantom1_N}
  \end{figure}
\begin{figure}[H]
\centering
{\includegraphics[width=0.9\textwidth]{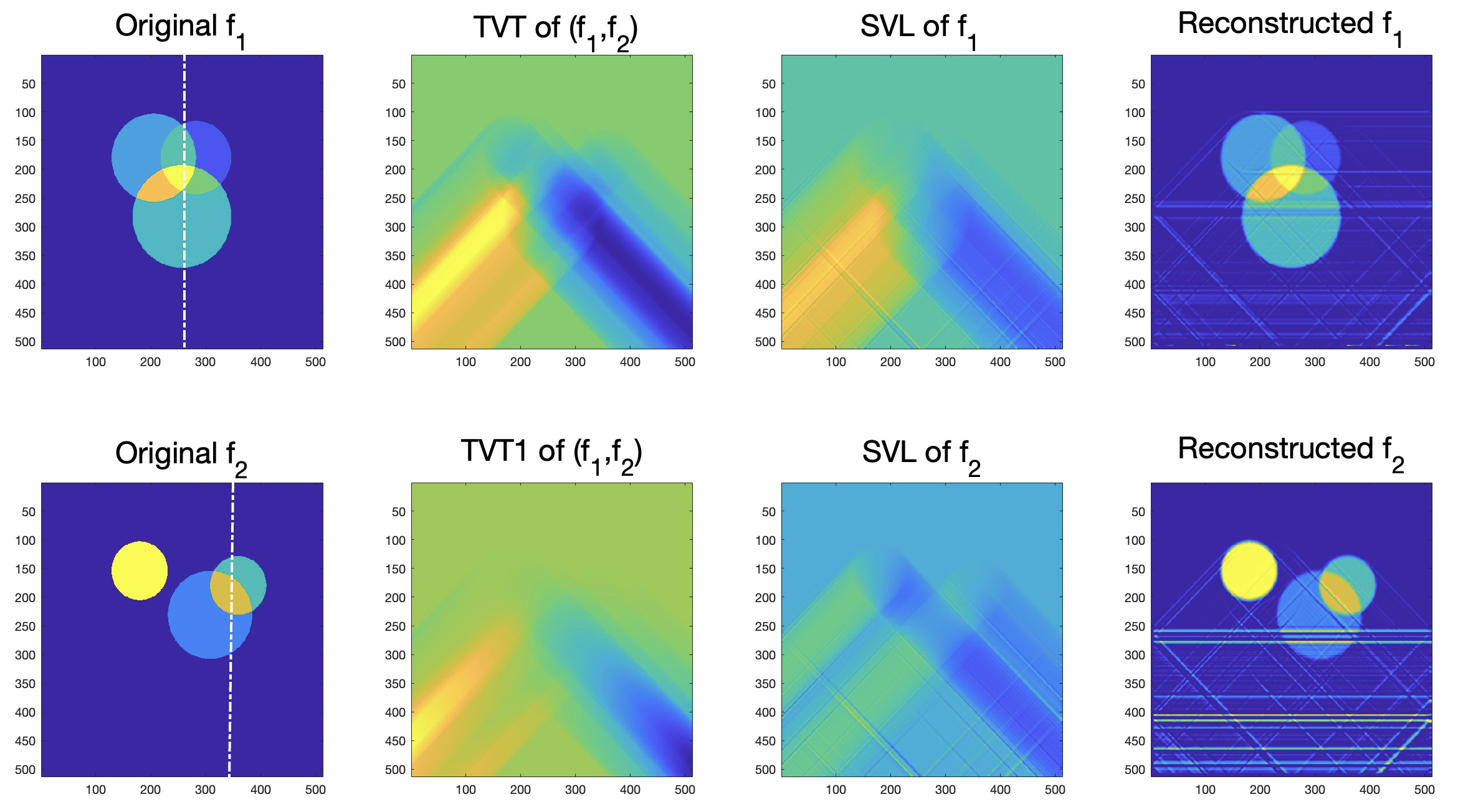}}
\caption{Components of the original field $\vf$ (column 1), $\Tc \vf$ and $\Jc \vf$ (column 2), signed V-line transform of the components (column 3), and reconstructed components of $\vf$ (column 4).}
    \label{T_Moments_phantom2}
  \end{figure}
\begin{figure}[H]
\centering
{\includegraphics[width=0.9\textwidth]{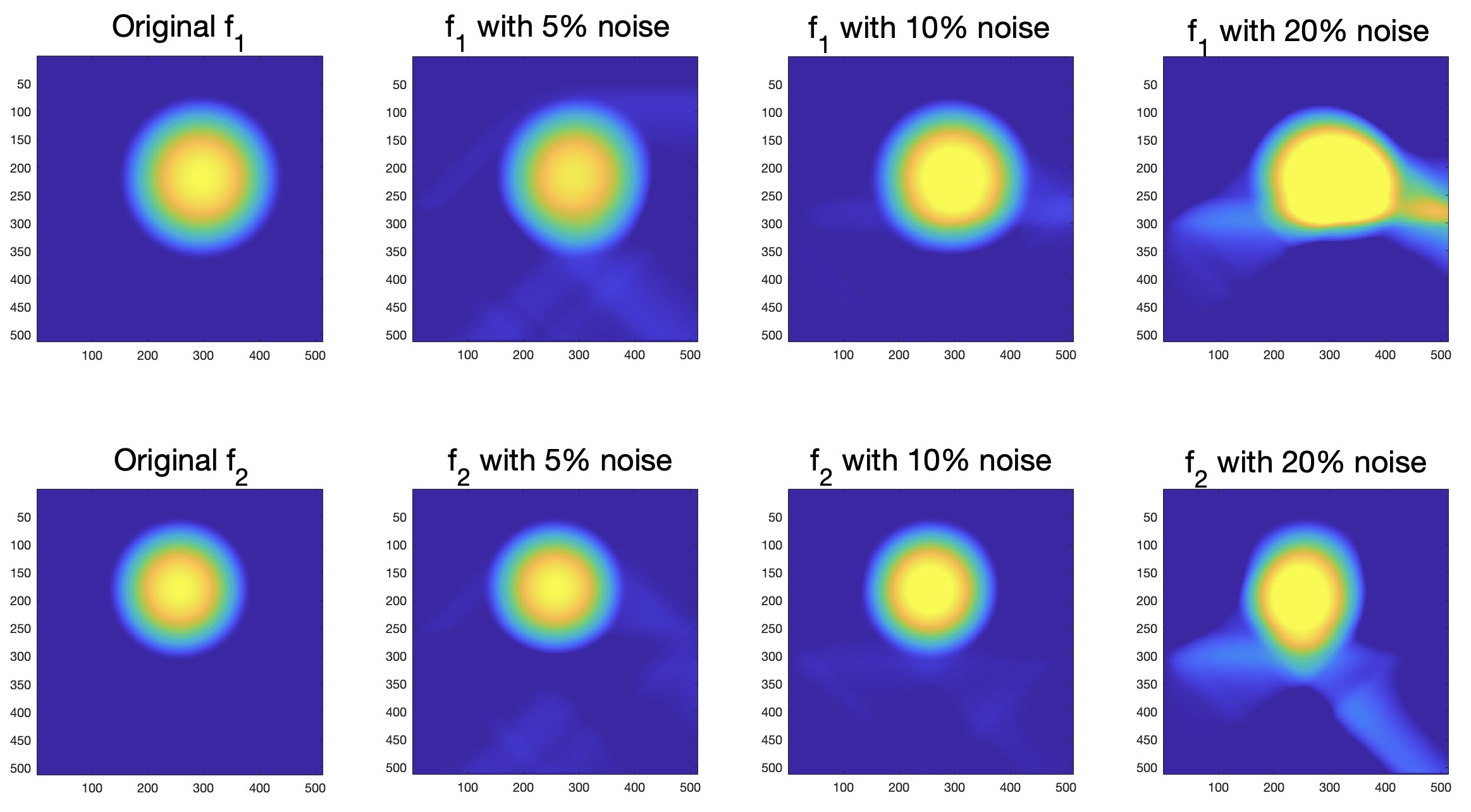}}
     \caption{Reconstructions of components of $\vf$ using $\Tc\vf$ and $\Jc\vf$ with  $5\%$, $10 \%$, and $20 \%$ noise.}
    \label{T_Moments_phantom2_N}
  \end{figure}

\begin{figure}[H]
\centering
{\includegraphics[width=0.85\textwidth]{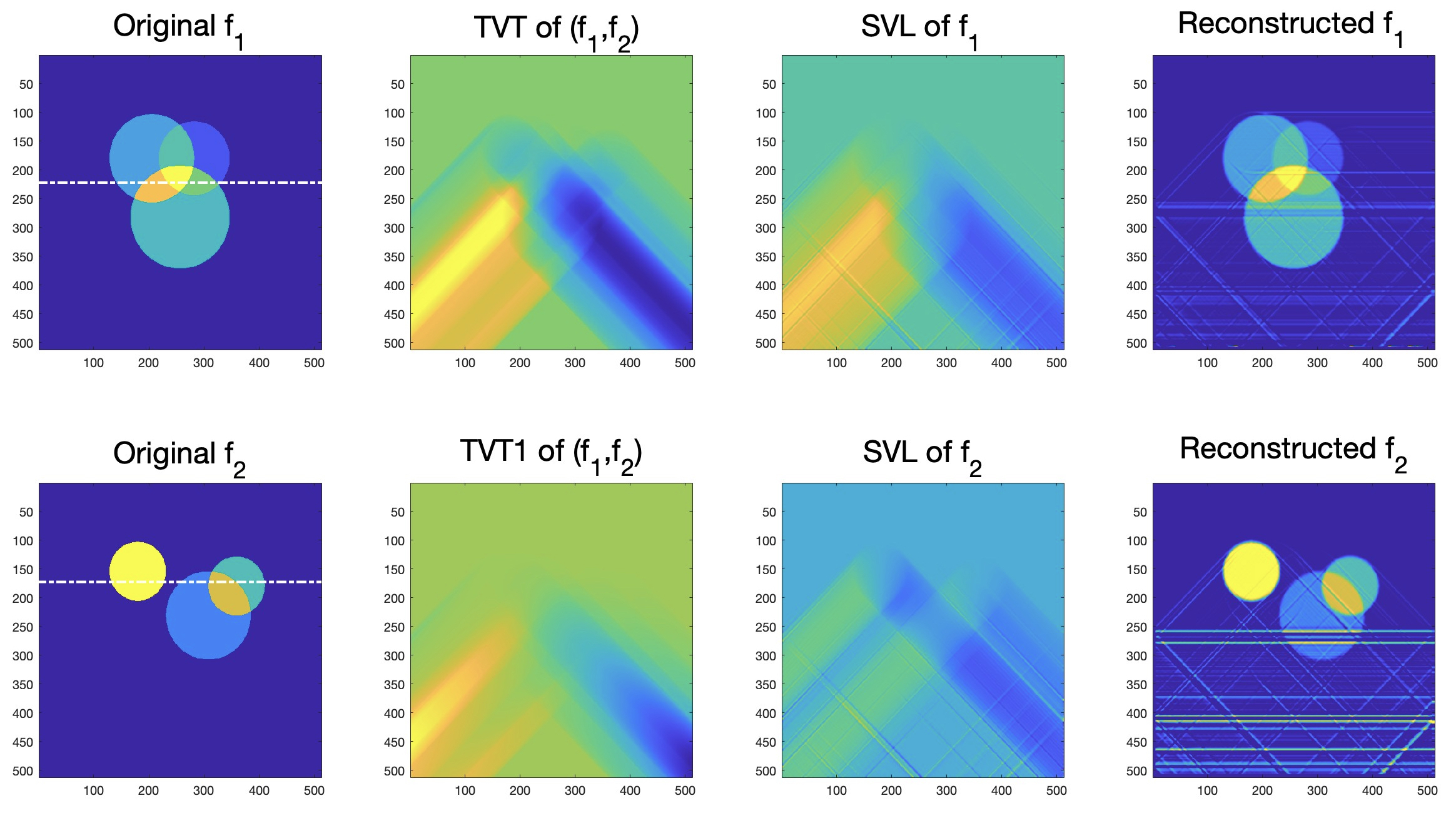}}
\caption{Components of the original field $\vf$ (column 1), $\Tc \vf$ and $\Jc \vf$ (column 2), signed V-line transform of the components (column 3), and reconstructed components of $\vf$ (column 4).}
    \label{T_Moments_phantom3}
  \end{figure}
\begin{figure}[H]
\centering
{\includegraphics[width=0.84\textwidth]{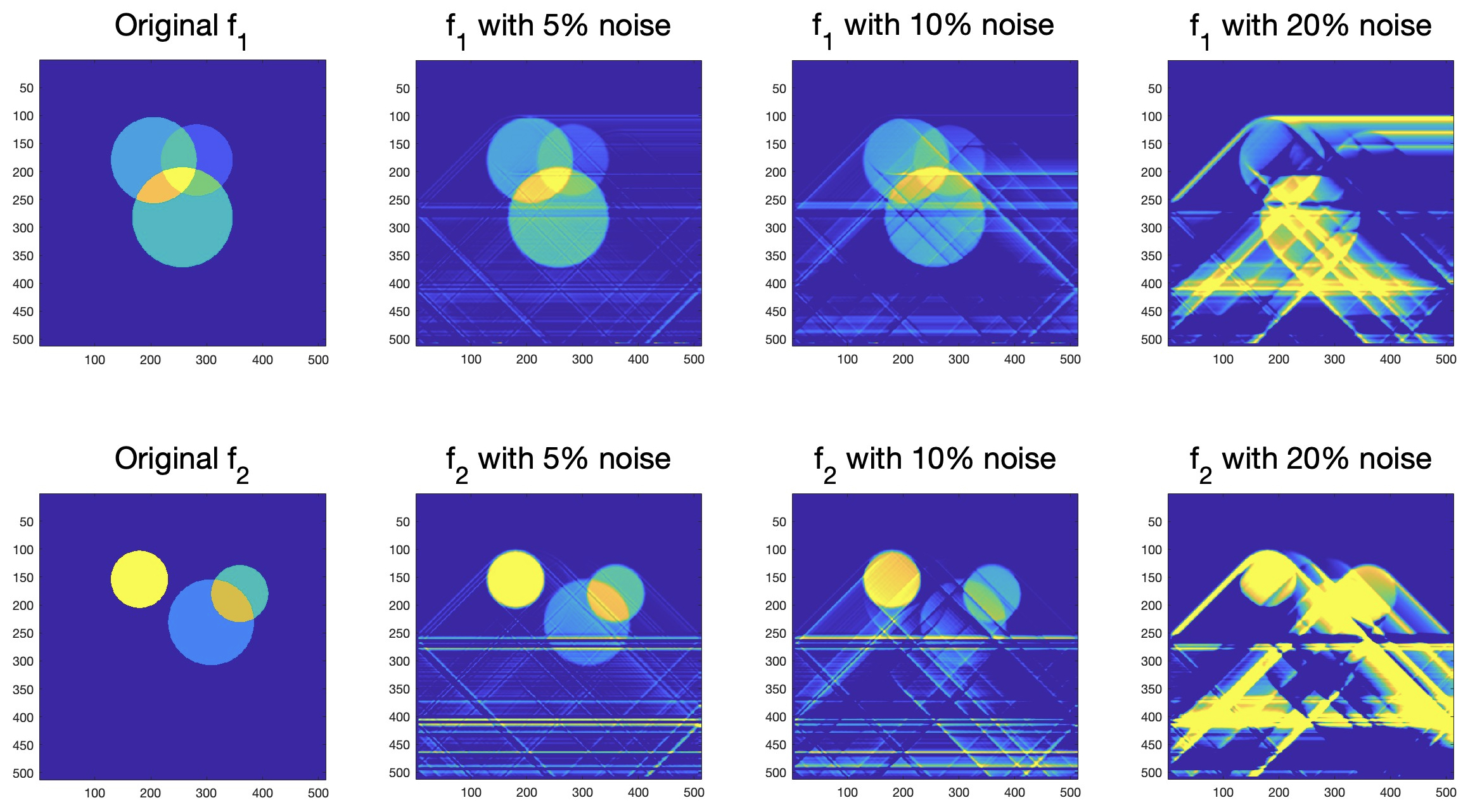}}
        \caption{Reconstructions of components of $\vf$ using $\Tc\vf$ and $\Jc\vf$ with  $5\%$, $10 \%$, and $20 \%$ noise.}
    \label{T_Moments_phantom3_N}
  \end{figure}
  \begin{figure}[H]
\centering
{\includegraphics[width=0.94\textwidth]{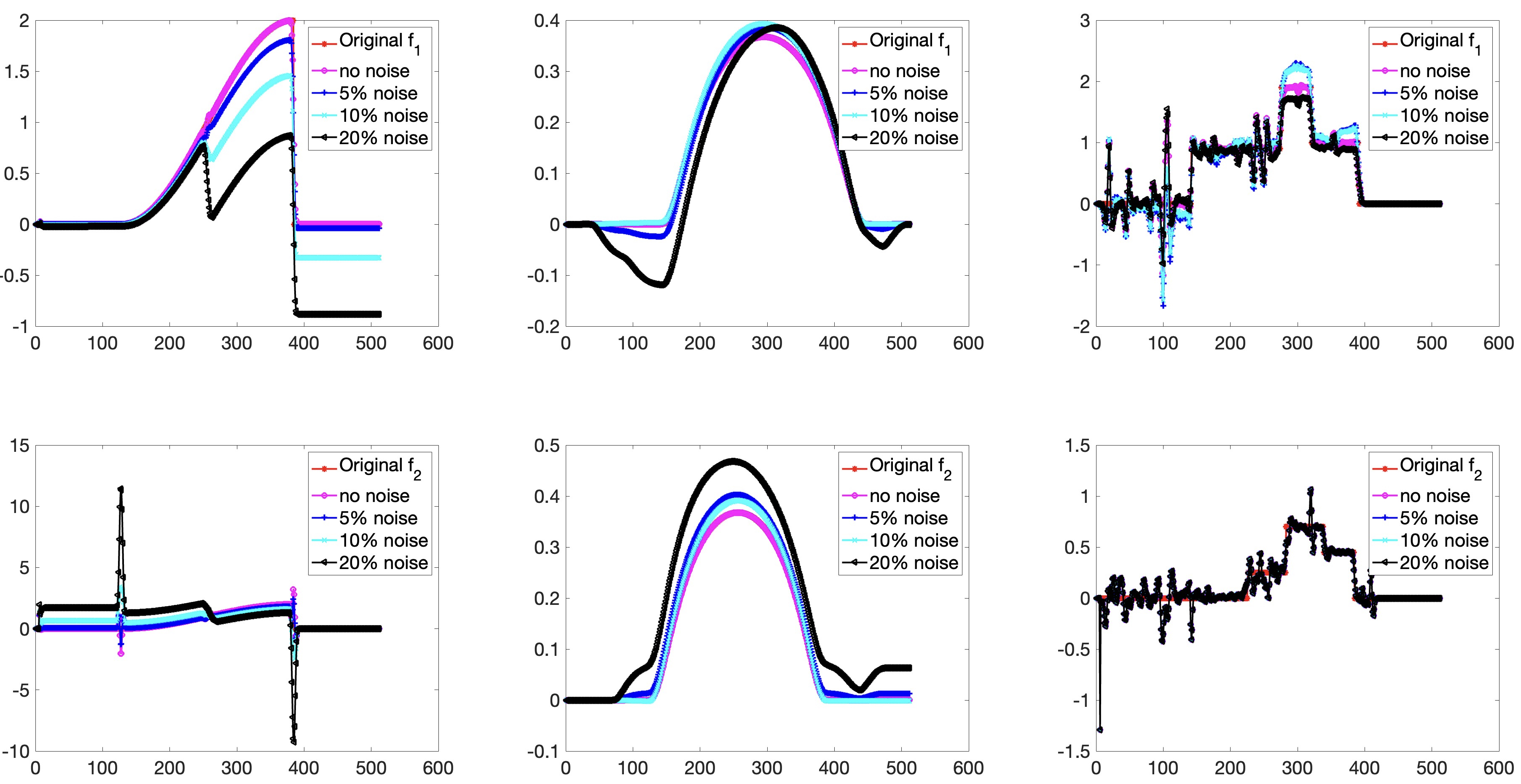}}
     \caption{Profile plots of $f_1$ and $f_2$ reconstructed from $\Tc(\vf)$ and $\Jc(\vf)$ with $0\%$, $5\%$, $10 \%$, and $20 \%$ noise. Plots in $j$-th column correspond to Phantom $j$, $j=1,2,3$.}
  \end{figure}
\begin{table}[h!]
\centering   
\begin{tabular}{ |p{1.6 cm}||p{0.6cm}||p{1.6cm}|p{1.6cm}|p{1.7cm}|p{1.7cm}|  }
 \hline
Phantoms & $\vf$ & Not noise & 5\% noise  & 10\% noise & 20\% noise\\
 \hline
PH1   & $f_1$    & 53.52\%  & 54.89\% &	60.11\% & 412.32\%\\
 \hline
PH1   &   $f_2$  & 6.91\% &	9.88\%	& 16.80\%	 & 110.01\%\\
 \hline
 PH2  &  $f_1$ & 1.49\%	& 16.98\% & 19.49\% & 62.16\%\\
  \hline
 PH2    & $f_2$ & 0.95\%	& 8.18\% & 8.61\% & 38.10\%\\
 \hline 
 PH3  &  $f_1$ & 18.72\%	& 19.10\% & 35.38\% & 82.85\%\\
  \hline
 PH3    & $f_2$ & 94.12\%	& 95.43\% &	95.81\% &158.23\%\\
 \hline 
\end{tabular}
\caption{Relative errors of the reconstructions of $f_1$ and $f_2$ from $\Tc\vf$ and $\Jc\vf$.}
\label{table:4}
\end{table}
\subsection{Effects of the angle between the rays of the V-line on reconstructions}\label{Sec:angles-effect}
In all numerical results presented up to this point, the unit vectors defining the V-lines were $\vu = (\cos \phi, \sin \phi)$ and $\vv = (\cos (\pi-\phi), \sin (\pi-\phi))$, where $\phi=\pi/4$. To test the effects of the V-line opening angle $\pi-2\phi$ on the reconstructions, we have run numerical simulations for various other angles $\phi\in(0,\pi/2)$. The results show that the inversion method using the combination of LVT and TVT data is very robust and works well for all angles and all phantoms. The methods using LVT or TVT with their corresponding moments produce accurate results on smooth phantoms. However, when applied to piecewise constant phantoms, the quality of reconstruction deteriorates as the opening angle moves away from $\pi/2$. The figures below show a representative sample of reconstructions with different V-line opening angles. Since the quality of reconstructions is very similar for $f_1$ and $f_2$, we show only the results for $f_1$.

  \begin{figure}[H]
\centering
{\includegraphics[width=1.0\textwidth]{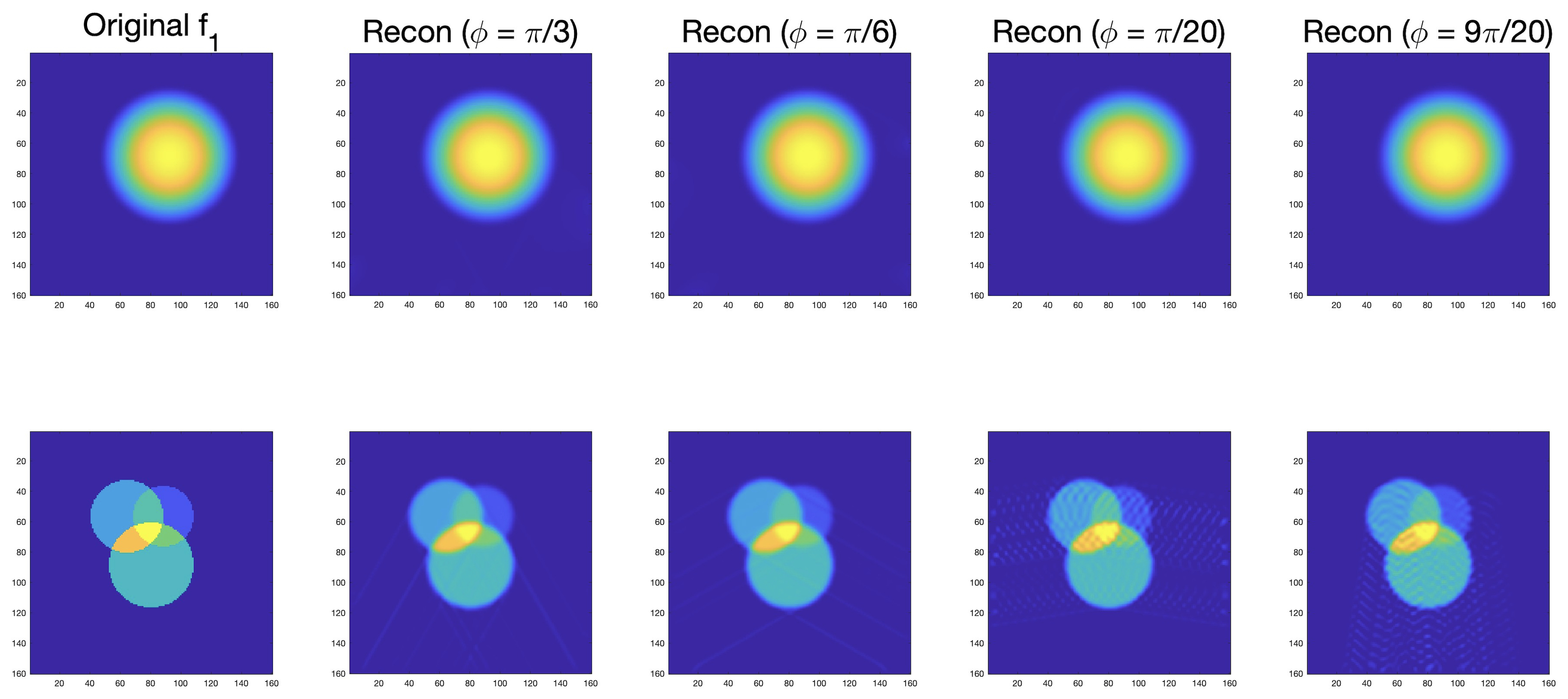}}
     \caption{Reconstructions using a combination of LVT and TVT data}
  \end{figure}
    \begin{figure}[H]
\centering
{\includegraphics[width=1.0\textwidth]{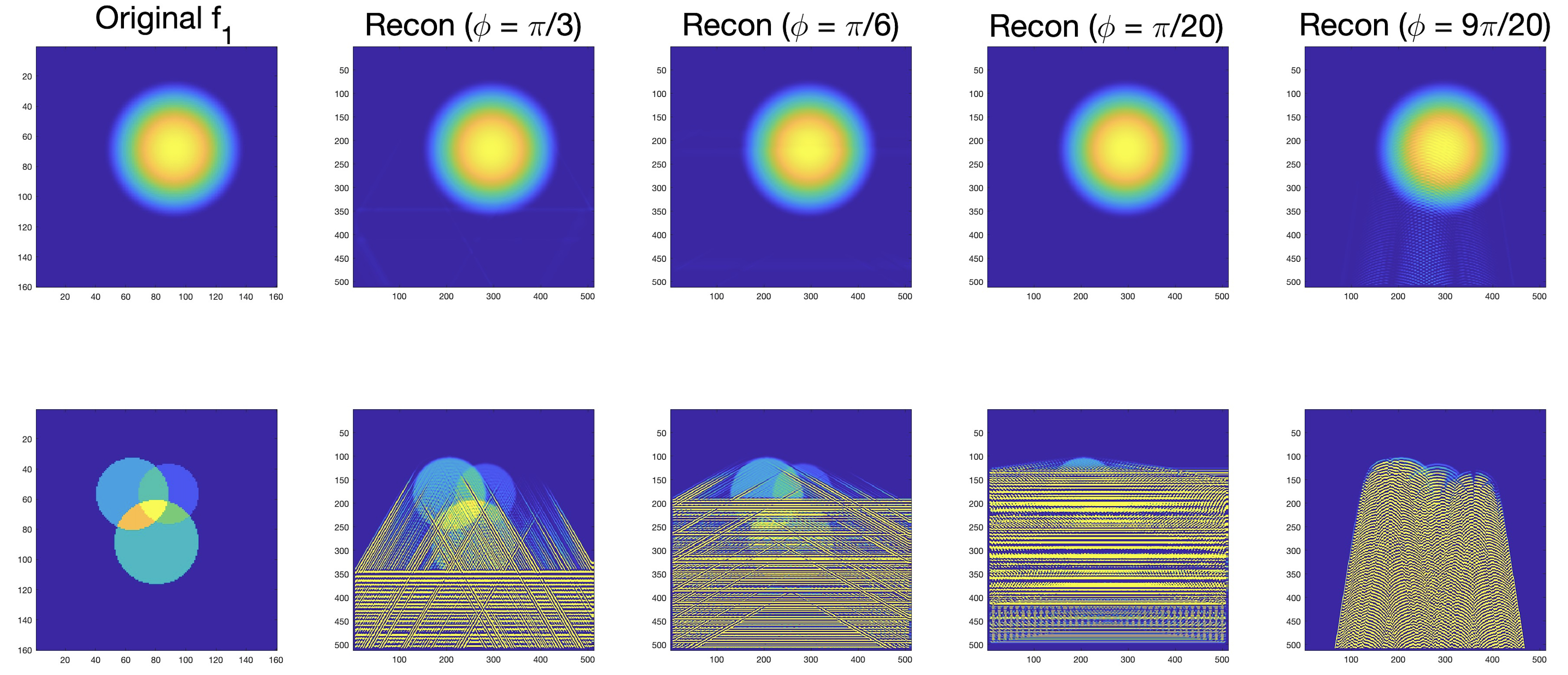}}
     \caption{Reconstructions using a combination of LVT and LVT1 data}
  \end{figure}
The artifacts appearing in these reconstructions have been explained in Section \ref{Sec:angles-effect}, but we would like to emphasize a few things here. The horizontal artifacts start at the locations of an abrupt cut of the (non-compactly supported) data on the edges of the (compactly supported) image domain. Notice that when $\phi=9\pi/20$, these cuts happen outside of the field of view, and there are no horizontal artifacts. The strength of the artifacts gradually increases as the V-line opening angle moves away from $\pi/2$. To demonstrate that, we present below a few additional simulations with various angles $\phi$ close to $\pi/4$. The analysis of the strengths of such singularities and the development of various techniques for their reduction are interesting and non-trivial topics of research in microlocal analysis, which are beyond the scope of this article. We refer the reader interested in this subject to \cite{Krishnan2015} and the references therein. 
  
\begin{figure}[H]
\centering
{\includegraphics[width=1.0\textwidth]{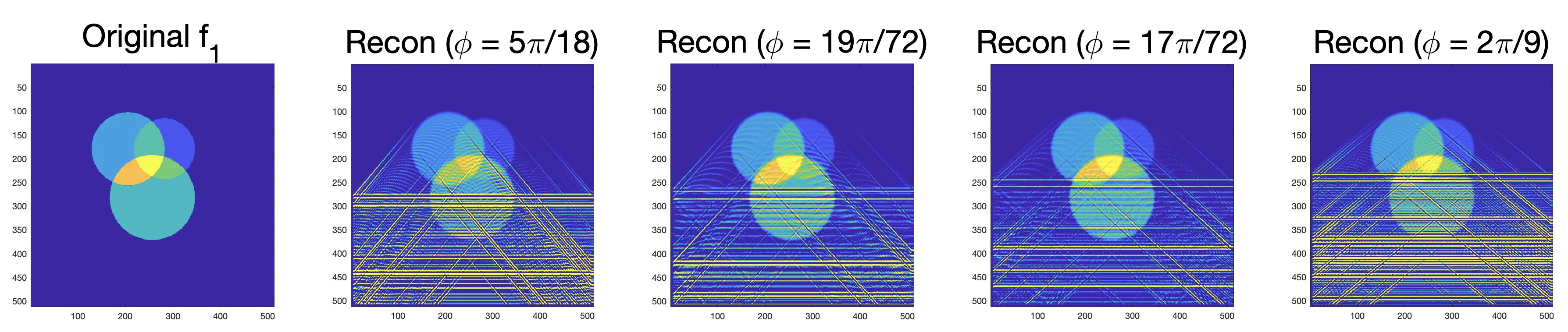}}
     \caption{Reconstructions using a combination of LVT and LVT1 data}
  \end{figure}

The reconstructions from TVT and TVT1 demonstrate the same type of behavior as those from LVT and LVT1 and are not presented here for the benefit of space and to avoid redundancy.
 \subsection{Recovery of a vector field from its vector star transform}
 This subsection is devoted to the reconstruction of a vector field $\vf$ from its vector star transform $\Sc \vf$. In our numerical simulations(see Figures \ref{Star_phantom1} - \ref{Star_phantom3_N}) we consider the stars with a variable location of the vertex and three branches directed along $\vgamma_i = (\cos \phi_i, \sin \phi_i)$, where $\phi_1 = 0$, $\phi_2 =  2\pi/3$, and $\phi_3 = 4\pi/3$. Recall from Theorem \ref{th:star} that the (component-wise) Radon transform of the unknown vector field $\vf$ is expressed in terms of its vector star transform as follows:
\begin{equation}\label{eq:radon of f}
\mathcal{R}\vf\,{(\vpsi,s)} = \begin{bmatrix}
 \vgamma(\vpsi) \\
 \vgamma(\vpsi)^{\perp}
\end{bmatrix}^{-1} \frac{d}{ds}\mathcal{R}(\mathcal{S}\vf) (\vpsi,s),
\end{equation}
where 
\begin{equation}\label{eq: def gamma}
\vgamma(\vpsi) := -\sum_{i=1}^3 \frac{\vgamma_i}{ \vpsi \cdot \vgamma_i } \; \in \mathbb{R}^2, \quad c_i = 1\  (\mbox{for } i = 1, 2, 3).  
\end{equation}
By applying (component-wise) $\mathcal{R}^{-1}$ to the above identity, we recover $\vf$. Numerically, the Radon transform and its inverse are carried out through the Matlab in-built functions \texttt{radon} and \texttt{iradon}. The actions of these functions can be briefly described as follows.

\vspace{2mm}
\noindent\texttt{radon}: takes as an input an $n \times n$ pixelized image $F$ and generates its Radon transform $\Rc F(\psi, s)$ for angles $ \psi  = 0, 1, \dots, 179$ (in degrees), and $s=-m,\ldots,m$, where $m = \lfloor n/\sqrt{2}\rfloor + 2$. 

\vspace{2mm}
\noindent\texttt{iradon}: is used to invert the Radon transform and get back the image $F$ from  $\Rc F$.\\

\noindent We break down our procedure of inverting the vector start transform $\mathcal{S}$ into the following steps:
\begin{itemize}
    \item The star transform data $\Sc \vf$ is represented by a pair of $512\times512$ matrices, one corresponding to its longitudinal part and the other to the transverse part (recall formula (\ref{def:star})). We generate them by numerically evaluating the divergent beam transforms $\Xc_{\vgamma_i}\vf$, $i=1,2,3$ and adding them up. Since $\Sc \vf$ has an unbounded support even for a compactly supported vector field $\vf$, the matrices described above represent a truncated approximation of $\Sc \vf$.
\item   We use the function \texttt{radon} to generate $\Rc (\Sc \vf)(\psi_i, s_j)$, which is represented by a pair of $180 \times 729$ matrices. Here $\vpsi = (\psi_i)_{i = 0}^{179}$ is the vector of projection angles in degrees and $(s_j)_{j = 1}^{729}$ is the discretization of the radial variable used for parameterization of the Radon transform.  The truncation of $\Sc\vf$ described above results in numerical errors in the evaluation of $\Rc (\Sc \vf)(\psi_i, s_j)$ along the lines that pass through the truncated ``tails'' of $\Sc\vf$.
\item In the third step we apply the Matlab built-in function \texttt{gradient} to compute $\displaystyle \frac{d}{ds}\mathcal{R}(\mathcal{S}\vf) (\vpsi,s)$. 
   \item Next, for each value of discretized angle $\psi$ we multiply $\displaystyle \frac{d}{ds}\mathcal{R}(\mathcal{S}\vf) (\vpsi,s)$ by the $2\times2$ matrix $ \displaystyle \begin{bmatrix}
 \vgamma(\vpsi) \\
 \vgamma(\vpsi)^{\perp}
\end{bmatrix}^{-1}$, where $\vgamma(\vpsi)$ is given by equation \eqref{eq: def gamma}. This multiplication generates the Radon transforms $\Rc f_1$ and $\Rc f_2$.  
\item Finally, we apply the Matlab built-in function \texttt{iradon} to $\Rc f_1$ and $\Rc f_2$ to get $f_1$ and $f_2$.
\end{itemize}

\begin{rem}
    The errors in data described in the second step of the above list, spread further by the follow-up steps of differentiation and matrix multiplication, resulting in artifacts at the edges of the unit square in reconstructed images. 
    Similar artifacts also appear in the numerical inversions of the star transform on scalar fields (e.g. see \cite{Amb_Lat_star}).
\end{rem}

\begin{figure}[H]
\centering
{\includegraphics[width=0.9\textwidth]{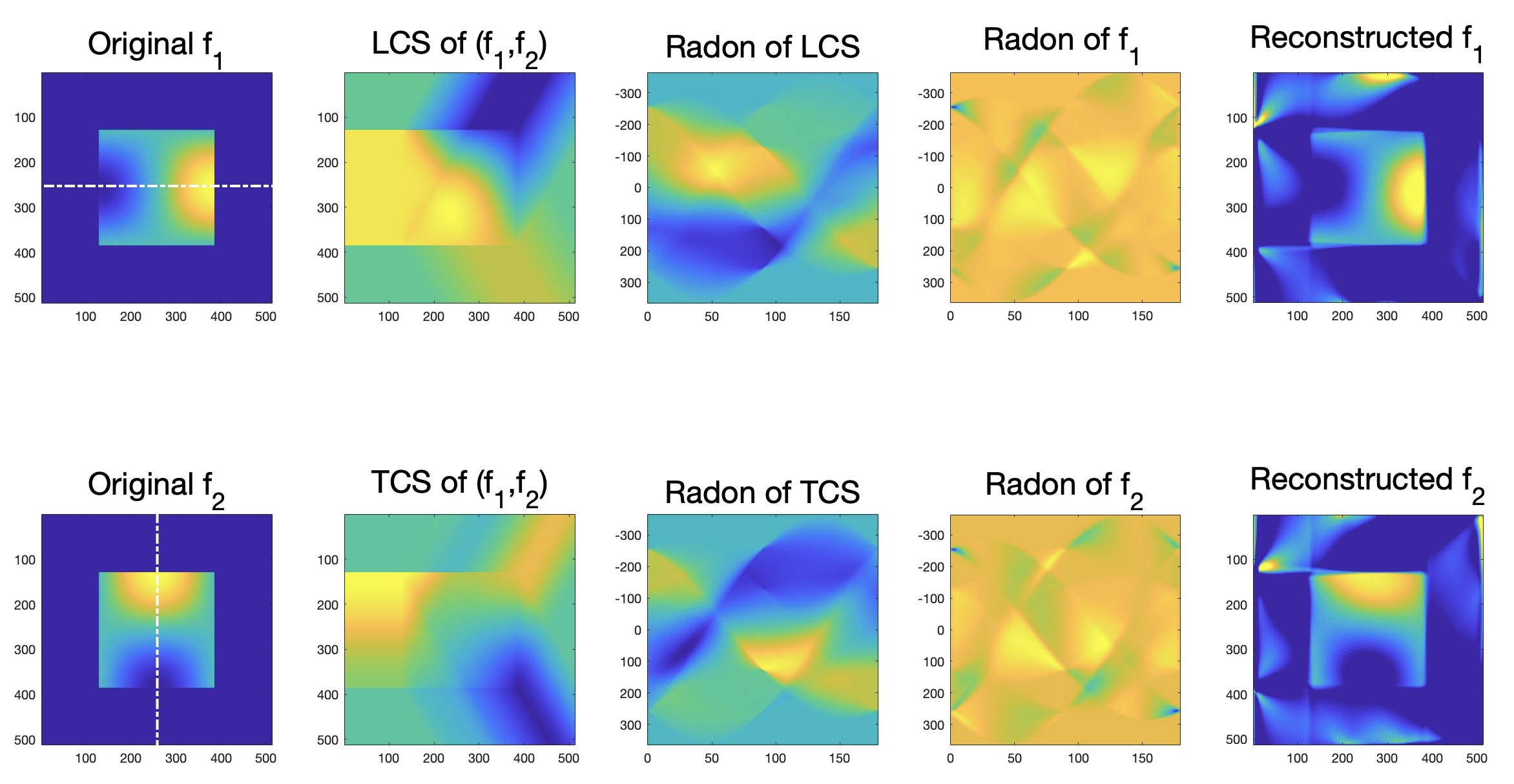}}
\caption{Components of  $\vf$ (column 1), longitudinal (LCS) and transversal (TCS) components of $\Sc\vf$ (column 2), corresponding Radon transforms (column 3), $s$-derivative (column 4), Radon transform of components of $\vf$ (column 5), and reconstructed components of $\vf$ (column 6).}
    \label{Star_phantom1}
  \end{figure}
\begin{figure}[H]
\centering
{\includegraphics[width=0.9\textwidth]{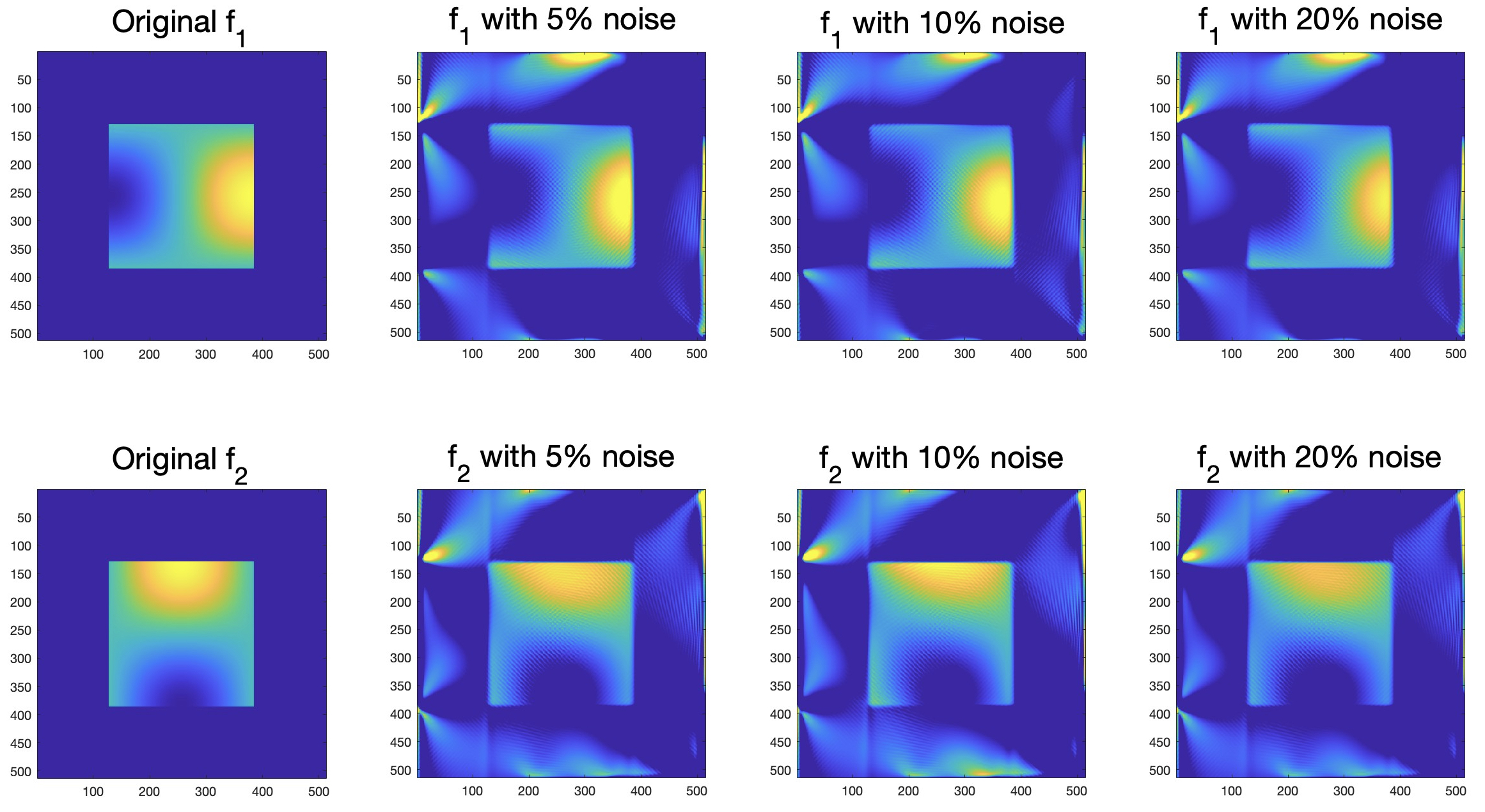}}
     \caption{Reconstructions with  $5\%$, $10 \%$, and $20 \%$ noise.}
    \label{Star_phantom1_N}
  \end{figure}

\begin{figure}[H]
\centering
{\includegraphics[width=0.9\textwidth]{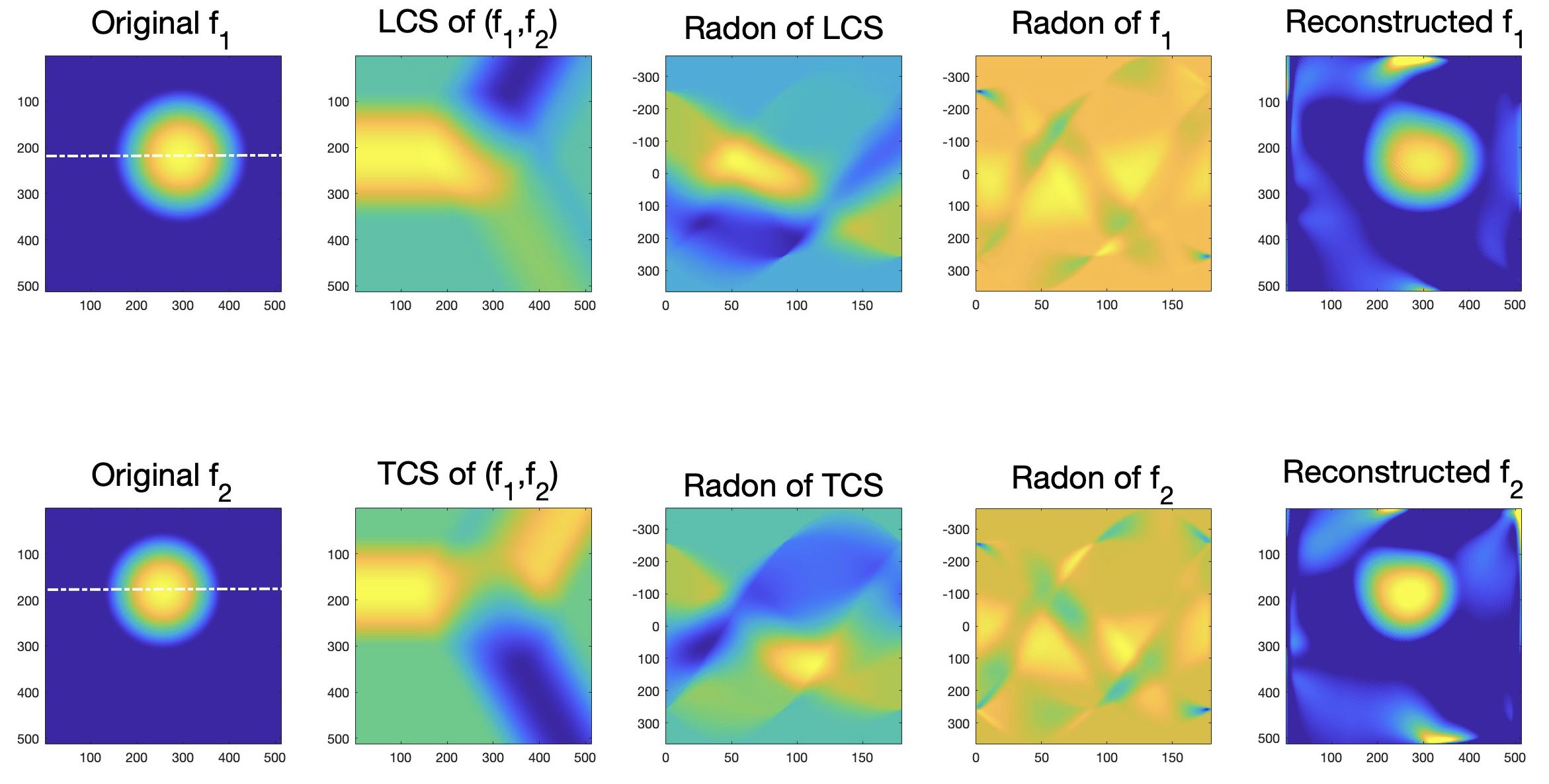}}
   \caption{Components of  $\vf$ (column 1), longitudinal (LCS) and transversal (TCS) components of $\Sc\vf$ (column 2), corresponding Radon transforms (column 3), $s$-derivative (column 4), Radon transform of components of $\vf$ (column 5), and reconstructed components of $\vf$ (column 6).}
    \label{Star_phantom2}
  \end{figure}
\begin{figure}[H]
\centering
{\includegraphics[width=0.9\textwidth]{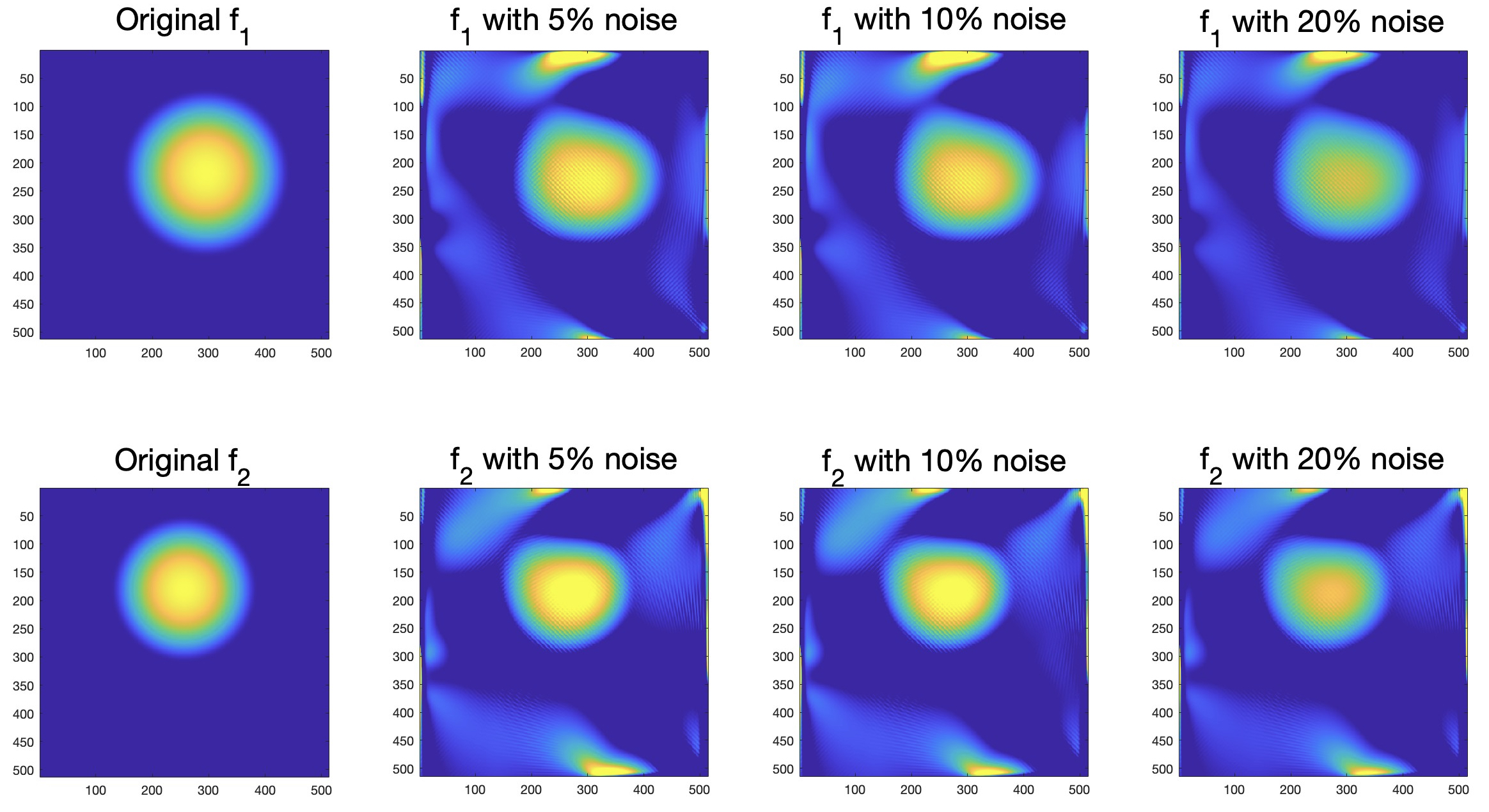}}
    \caption{Reconstructions with  $5\%$, $10 \%$, and $20 \%$ noise.}
    \label{Star_phantom2_N}
  \end{figure}
  
\begin{figure}[H]
\centering
{\includegraphics[width=0.9\textwidth]{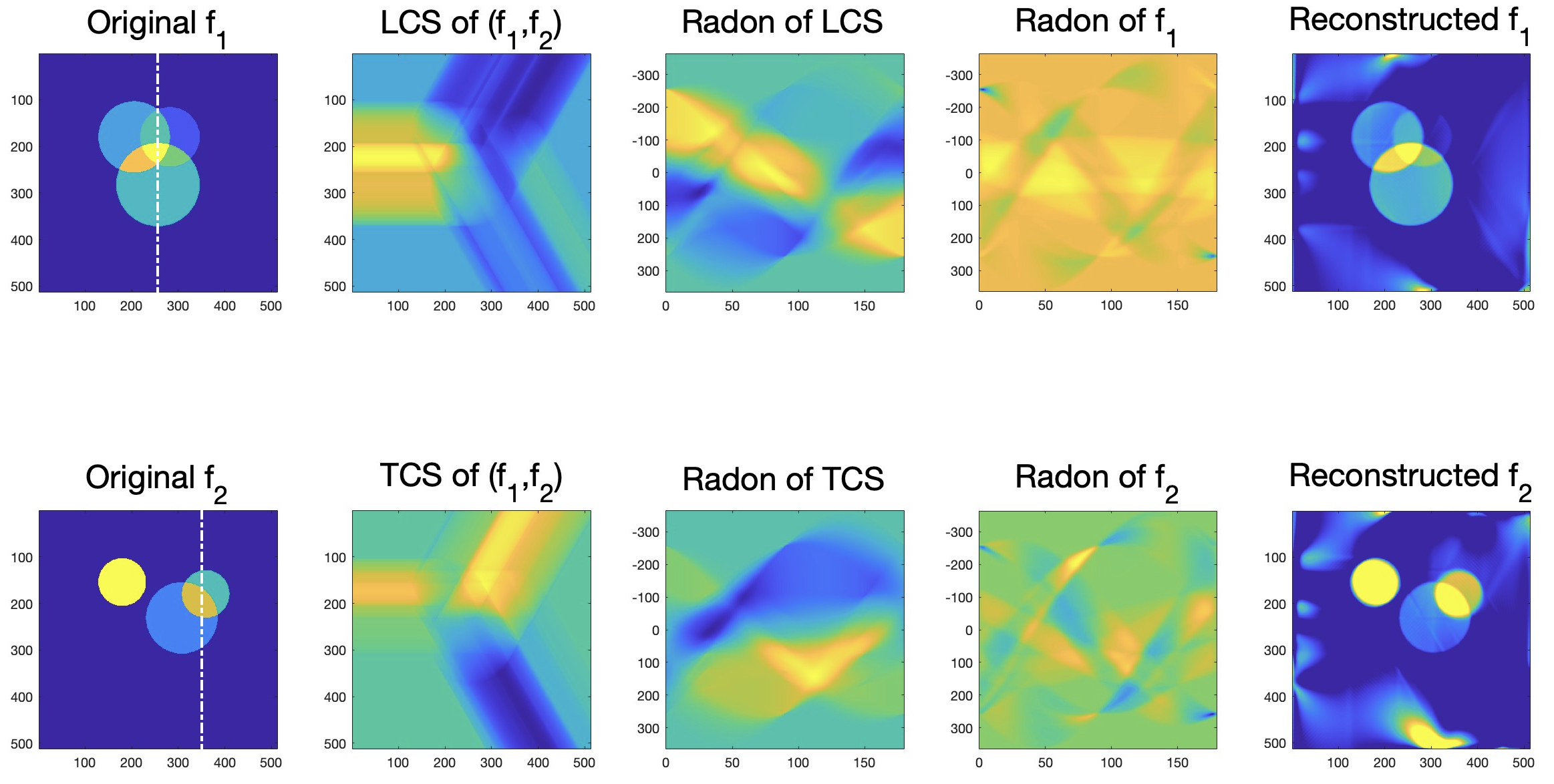}}
      \caption{Components of  $\vf$ (column 1), longitudinal (LCS) and transversal (TCS) components of $\Sc\vf$ (column 2), corresponding Radon transforms (column 3), $s$-derivative (column 4), Radon transform of components of $\vf$ (column 5), and reconstructed components of $\vf$ (column 6).} 
    \label{Star_phantom3}
  \end{figure}
\begin{figure}[H]
\centering
{\includegraphics[width=0.9\textwidth]{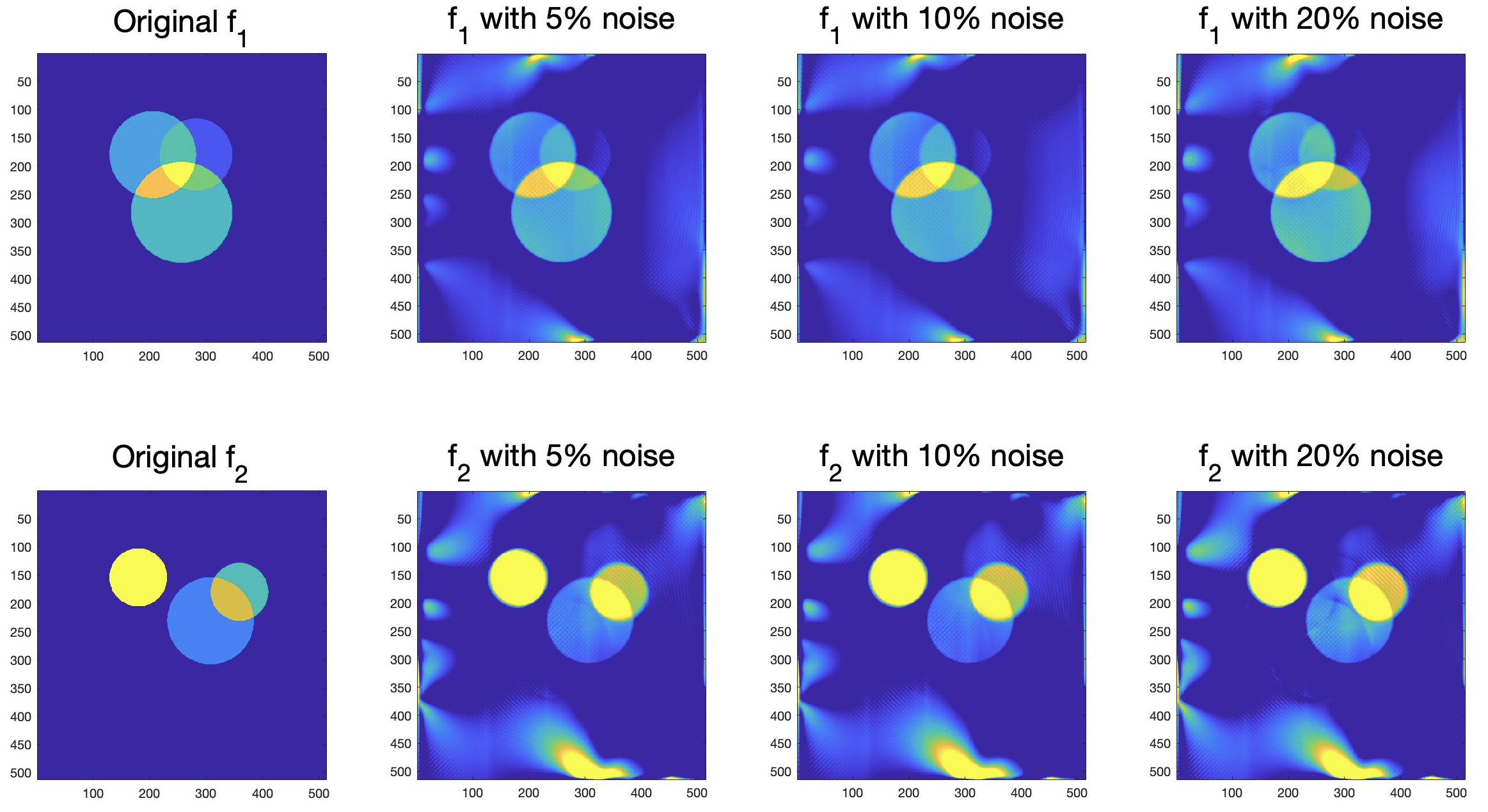}}
 \caption{Reconstructions with  $5\%$, $10 \%$, and $20 \%$ noise.}
    \label{Star_phantom3_N}
  \end{figure}

  \begin{figure}[H]
\centering
{\includegraphics[width=0.94\textwidth]{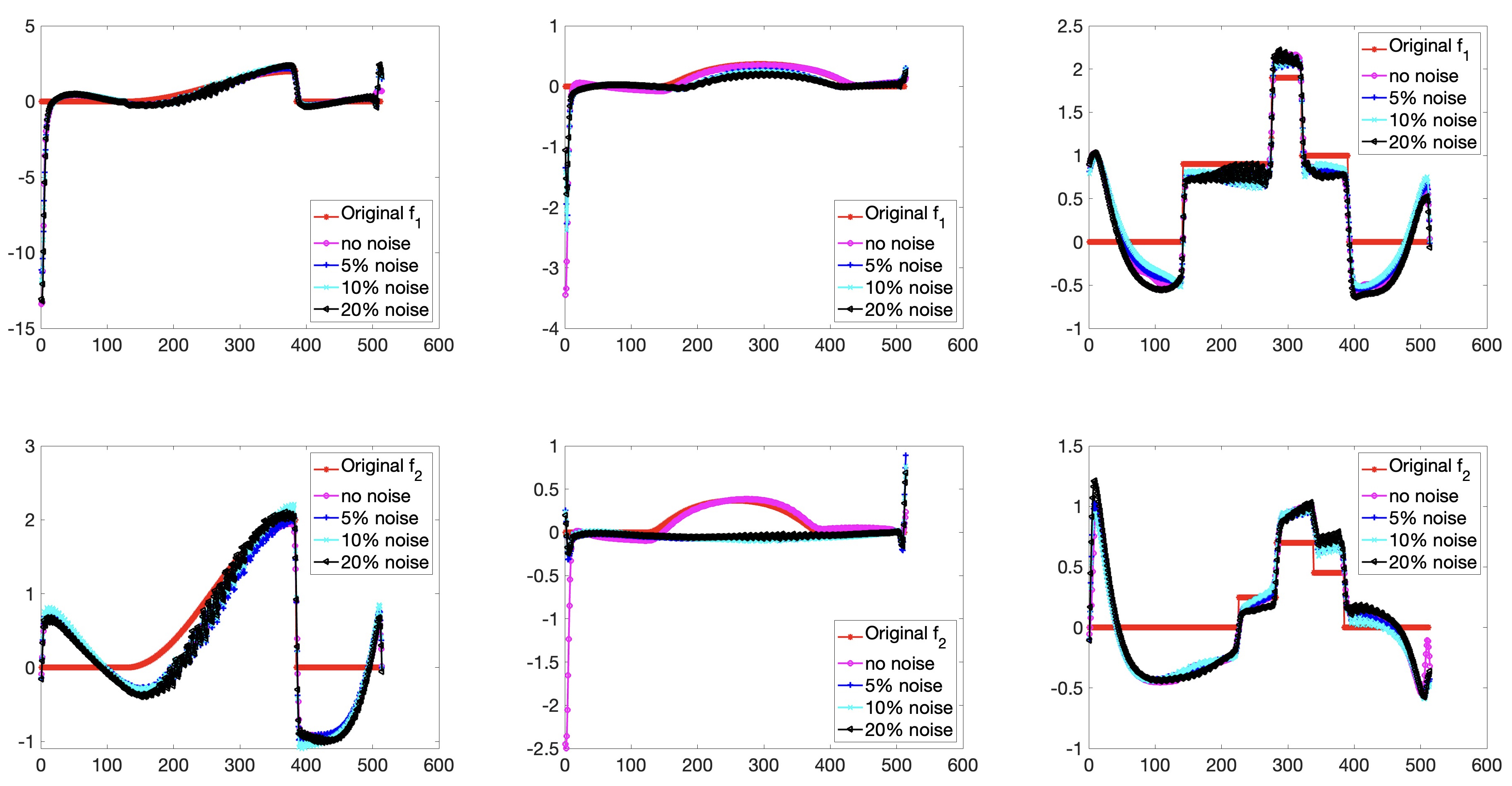}}
     \caption{Profile plots of $f_1$ and $f_2$ reconstructed from $\Sc(\vf)$ with $0\%$, $5\%$, $10 \%$, and $20 \%$ noise.  Plots in $j$-th column correspond to Phantom $j$, $j=1,2,3$.}
  \end{figure}

\begin{table}[h!]
\centering
\begin{tabular}{ |p{1.6 cm}||p{0.6cm}||p{1.6cm}|p{1.6cm}|p{1.7cm}|p{1.7cm}|  }
 \hline
 Phantoms & $\vf$ & No noise & 5\% noise  & 10\% noise & 20\% noise\\
 \hline
PH1   & $f_1$    & 105.01\% & 113.48\% &	122.81\% & 147.59\%\\
 \hline
PH1   &   $f_2$  & 111.48\% & 112.09\% &	123.21\% & 135.12\%\\
 \hline
 PH2  &  $f_1$ & 113.75\% &	128.9\% &	130.57\% & 232.02\%\\
  \hline
 PH2    & $f_2$ & 104.76\% & 163.91\% &	182.58\% & 264.27\%\\
 \hline 
 PH3  &  $f_1$ & 101.60\% & 102.88\% & 111.82\% & 138.78\%\\
  \hline
 PH3    & $f_2$ & 129.59\% &	199.81\% &	202.73\% & 218.80\%\\
 \hline 
\end{tabular}
\caption{Relative errors of the reconstructions of $f_1$ and $f_2$ from $\Sc\vf$.}
\label{table:5}
\end{table}

\section{Vector Star Transform Reconstructions of RGB Images}\label{sec:Star-Python}

In this section, we consider the vector fields on $\Omega \subset \mathbb{R}^2$ as RGB images, where at each point (pixel) $\vx\in \Omega$  the two components of the vector field $\vf(\vx) = (f_1(\vx), f_2(\vx))$  represent the intensity of red and green colors. More explicitly, the discretized versions of the vector field components $R \approx f_1$, $G \approx f_2$ are $N\times N$ matrices representing the red and green layers of the image, i.e. the $(i,j)$-th entries of $R$ and $G$ have values between zero and one, respectively associated with the red and green intensity of the corresponding pixel. The values of the blue layer, $B$, have been ignored in our experiments. In our numerical demonstration, we use three different phantoms and generate their vector star transform data. Then, we apply to the data our formula (\ref{eq:radon of f}), followed by the component-wise inverse Radon transform to reconstruct each image. 

\vspace{2mm}

The Python code for these procedures is available as a notebook in the Google Colab: \href{https://colab.research.google.com/drive/1S2B8KUFpAGBURWTZ4L023fNh_1WRPBBB?usp=sharing}{Star Transform Reconstruction Library}. This notebook allows the user to customize the experiments either by providing as an input a layered NumPy tensor, or by uploading a colored image. The user can also define the branch directions of the star transform and run the experiment, generating the vector star transform followed by the reconstruction step.   

\vspace{2mm}


\textbf{Computing the vector star transform of an image.} Given an input image $\texttt{img}$ of dimensions $N\times N\times 3$, we compute the dot product of the branch direction $\vgamma_1$ and the vector given by the red and blue components of $\texttt{img}$ (we use the red, $R=$\texttt{img[:,:,0]}, and green, $G=$\texttt{img[:,:,1]}, components and ignore the blue, $B=$\texttt{img[:,:,2]})    to obtain a $N\times N $ matrix $L \approx \vf\cdot \vgamma_1$. Then, we use the $\texttt{divergent-beam}$ procedure to obtain the divergent beam transform of $L$, resulting in the \textbf{longitudinal} transform of $\vf$ along the branch $\vgamma_1$. Note that the divergent beam transform of a scalar function $g(\vx)$ at the vertex $\vx_0 \in \mathbb{R}^2$ can be obtained by applying the standard Radon transform to the function $\chi_H (\vx)g(\vx)$, where $H$ is the appropriate half plane with $\vx_0$ on the boundary $\partial H$ and 
$$
\chi_H(\vx):= \begin{cases}1, & \vx \in H, \\ 0, & \vx \notin H.\end{cases}
$$

The above observation reduces the computation of the divergent beam transform to the standard Radon transform, for which we use the Python function $\texttt{skimage.transform.radon}$. To obtain the vector star transform, we add up the contributions from all branches. The computation of the transversal component is done in the same fashion, by taking the (truncated) Radon transform of the transversal component $T \approx \vf^{\perp}\cdot \vgamma_1$ of the vector field $\vf$. In our experiments below, we use three branches along vectors with polar angles  $0$, $3\pi/4$  and $3\pi/2$. \\

\textbf{Reconstructing an image from its vector star transform.} The input for the reconstruction has two components (longitudinal and transversal). Our algorithm follows the inversion formula (\ref{eq:radon of f}). 

\vspace{2mm}

Steps of the reconstruction:  
\begin{itemize}
    \item Apply the standard Radon transform to both components of $\mathcal{S}\vf$.

    \item Use $\texttt{numpy.diff}$ function to compute the $d/ds$ derivative in the Radon domain.

    \item Multiply the resulting vector data by the matrix function 
    $\begin{bmatrix}
         \vgamma(\vpsi) \\
         \vgamma(\vpsi)^{\perp}
    \end{bmatrix}^{-1}$. 
This is the only step in the reconstruction where we have a mixture of the two components. 

    \item Apply the inverse Radon transform procedure ($\texttt{skimage.transform.iradon}$) in Python to reconstruct the image. 
\end{itemize}



  



\subsection*{Sample Reconstructions}

We apply the reconstruction algorithm described above to three different images presented below. The green and red components of the images are taken as components of the vector field. 

\begin{figure}[H]
    \centering
    \subfloat[\centering Red and green balls ($100\times 100$ pixels)]{{\includegraphics[width=4.6cm]{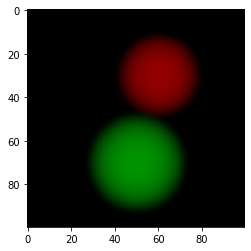} }}%
    \qquad
    \subfloat[\centering Overlapping rectangles ($100\times 100$ pixels)]{{\includegraphics[width=4.6cm]{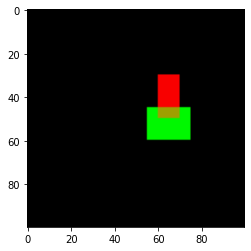} }}%
    \qquad
    \subfloat[\centering Rose ($300\times 300$ pixels)]{{\includegraphics[width=4.6cm]{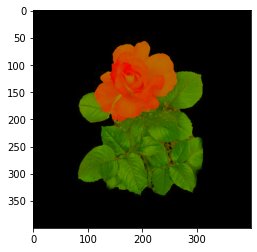} }}
    \caption{Original images (phantoms).}
\end{figure}

\begin{figure}[H]
    \centering
    \subfloat[\centering Red and green balls ($100\times 100$ pixels)]{{\includegraphics[width=4.5cm]{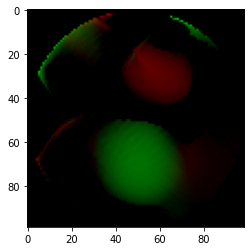} }}%
    \qquad
    \subfloat[\centering Overlapping rectangles ($100\times 100$ pixels)]{{\includegraphics[width=4.5cm]{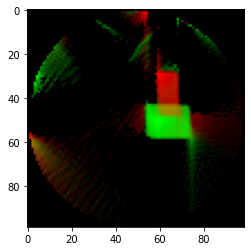} }}%
    \qquad
    \subfloat[\centering Rose ($300\times 300$ pixels)]{{\includegraphics[width=4.5cm]{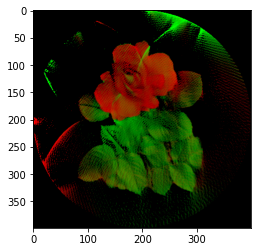} }}%
    \caption{Reconstructions from the vector star transforms of the original images.}%
\end{figure}

\section{Conclusions and Directions of Future Work}\label{sec:conclusions}

In this paper we discussed numerical implementations of various inversion schemes for generalized V-line transforms on vector fields introduced in \cite{Gaik_Mohammad_Rohit}. We demonstrated the possibility of efficient recovery of an unknown vector field from its longitudinal and transverse V-line transforms, their corresponding first moments, and the vector star transform. We examined the performance of our algorithms in a variety of setups with and without noise. 

The technique using a combination of LVT and TVT data proved to have the best characteristics, with the least amount of artifacts in reconstructions and no additional requirements on the support of the vector field. Moreover, this method preserved the high quality of reconstructed images for a wide range of V-line opening angles. The reconstructions using moment transforms and the vector star transform data had artifacts similar to those appearing in numerical inversions of generalized VLTs on scalar fields. In addition to that, in these cases the transform data were required to be known in a larger domain than the support of the vector field. The technique using moment transforms proved to be sensitive to the V-line opening angle when applied to non-smooth phantoms, producing the best reconstructions when that angle is close to $\pi/2$.

The vector tomography problems studied here and in the article \cite{Gaik_Mohammad_Rohit} can also be considered in higher dimensions. In dimensions $n\geq 3$, the longitudinal V-line transform can be defined in a similar fashion as for $n = 2$, while to define the transverse V-line transforms one needs to make a choice for $n-1$ linearly independent transverse directions. Once a choice for the transverse directions is made, questions like injectivity, exact inversion formulae, and numerical inversion algorithms can be asked for these transforms as well. Notice that the family of V-lines in $\mathbb{R}^n$ has $3n-2$ degrees of freedom. Therefore, to have a formally determined inverse problem one would have to consider a judiciously chosen subset of V-lines. Another approach is to consider conical transforms (generalized Radon transforms integrating over conical surfaces) for vector fields in $\mathbb{R}^n$, $n \geq 3$. In this case, there is a natural way to define the transverse transform, but one has to choose $n-1$ directions for longitudinal transforms. Here too one can study problems similar to those mentioned above. 
It is also natural to ask these questions for higher-order tensor fields in $\mathbb{R}^2$ as well as in higher dimensions. In a recent article \cite{Gaik_Rohit_Indrani}, we have studied V-line tensor tomography problem for symmetric $2$-tensor fields in $\mathbb{R}^2$ and addressed questions like kernel description, injectivity, and exact inversion formulae for longitudinal, transverse, and mixed V-line transforms. We feel that similar results should be possible for higher-order tensor fields as well, at least in $\mathbb{R}^2$. The authors plan to address several of the aforementioned problems in future research work and upcoming publications.


\section{Acknowledgements}\label{sec:acknowledge}
GA was partially supported by the NSF grant DMS 1616564 and the NIH grant U01-EB029826. RM was partially supported by SERB SRG grant No. SRG/2022/000947.


\bibliography{references-num}
\bibliographystyle{plain}






\end{document}